\documentclass{amsart}

\usepackage[arrow, matrix, curve]{xy}
\usepackage{amsmath,amsthm,amssymb}
\usepackage{amscd,graphics}
\usepackage[dvips]{graphicx}

\DeclareMathOperator{\spane}{span}
\DeclareMathOperator{\clsp}{\overline{span}}

\DeclareMathOperator{\Aut}{Aut}

\DeclareMathOperator{\Lat}{Lat}

\newcommand{\LL}{\mathcal L}
\newcommand{\M}{\mathcal M}
\newcommand{\DR}{\mathcal DR}
\newcommand{\TT}{\mathcal T}
\newcommand{\SSS}{\mathcal S}

\newcommand{\NN}{\mathcal N}

\newcommand{\JJ}{\mathcal J}
\newcommand{\KK}{\mathcal K}
\newcommand{\RR}{\mathcal R}

\renewcommand{\a}{\widetilde a}

\newcommand{\OO}{\mathcal{O}}

\newcommand{\X}{\widetilde{X}}

\newcommand{\al}{\alpha}

\newcommand{\FF}{\mathcal F}
\newcommand{\A}{\mathcal A}

\newcommand{\C}{\mathbb C}

\newcommand{\Z}{\mathbb Z}
\newcommand{\N}{\mathbb N}

\newtheorem{theorem}{Theorem}[section]
\newtheorem{lemma}[theorem]{Lemma}
\newtheorem{proposition}[theorem]{Proposition}
\newtheorem{corollary}[theorem]{Corollary}

\theoremstyle{definition}
\newtheorem{definition}[theorem]{Definition}
\newtheorem{example}[theorem]{Example}

\theoremstyle{remark}
\newtheorem{remark}[theorem]{Remark}
\newtheorem{notation}[theorem]{Notational conventions}

\numberwithin{equation}{section}
 
\begin{document}
   \title[Generalizing Cuntz-Pimsner and Doplicher-Roberts algebras]{$C^*$-algebras generalizing both relative Cuntz-Pimsner and Doplicher-Roberts algebras}
\author{B. K.  Kwa\'sniewski}
\address{Institute of Mathematics,  University  of Bialystok\\
 ul. Akademicka 2,  PL-15-267  Bialystok,
  Poland}
 \email{bartoszk@math.uwb.edu.pl}  
\urladdr{http://math.uwb.edu.pl/~zaf/kwasniewski}

  \thanks{  This work was in part supported by Polish Ministry of Science and High Education grant number N N201 382634. }  
\begin{abstract} We introduce and analyse the structure of $C^*$-algebras arising from ideals in right tensor $C^*$-precategories, which naturally unify the  approaches  based on Hilbert $C^*$-modules and  $C^*$-categories with tensor structure. We establish   an explicit intrinsic construction of the  algebras considered,  prove a number of key results  such as   structure theorem,  gauge-invariant uniqueness theorem, and describe the gauge-invariant ideal structure.  These  results give a new insight into the corresponding statements for relative Cuntz-Pimsner algebras and are  applied to Doplicher-Roberts algebras associated with $C^*$-correspondences.
\end{abstract}   
\maketitle 

\setcounter{tocdepth}{1}
 \tableofcontents
 
\section*{Introduction}
$C^*$-categories, the categorial analogues of (unital) $C^*$-algebras, arise quite naturally  in different problems of  representation theory, harmonic analysis, or cohomology theory, cf. \cite{glr}, \cite{v} and sources cited there. The recent interest in $C^*$-categories however is essentially due to a series of papers by S. Doplicher and J. E. Roberts where, motivated by  questions arising in  quantum field theory, they developed an abstract duality for compact groups, cf. \cite{dr}. In their  scenario an object of the dual group is represented by a certain tensor $C^*$-category $\TT$, and  
a machinery performing  this duality rests on a construction of a  $C^*$-algebra $\OO_\rho$,   influenced by Cuntz algebras  \cite{Cuntz},  associated functorially to each object $\rho$ of  $\TT$.  This association of $\OO_\rho$ can be  applied, with no substantial modifications, to the case where $\TT$ is just  a  \emph{ right tensor} $C^*$-category, i.e.  a $C^*$-category for which the set of objects is a unital semigroup, with identity $\iota$, and  for any object $\tau \in \TT$ there is a $^*$-functor $\otimes 1_\tau: \TT \to \TT$  (which intuitively  should be thought of as a tensoring on the right with  identity $1_\tau$ in the space of morphisms $\TT(\tau,\tau)$) such that  
$$
\otimes 1_\tau: \TT(\rho,\sigma)\to \TT(\rho\tau,\sigma \tau),\quad \textrm{and} \quad  \otimes 1_\iota= id ,\qquad ((a\otimes 1_\tau)\otimes 1_\omega)=  a \otimes 1_{\tau\omega},
$$
 where $\omega,\rho,\sigma\in \TT$, $a\in \TT(\rho,\sigma)$, and  we write $a\otimes 1_\tau$ for an element customary denoted by $\otimes 1_\tau(a)$. In this paper we adopt a \emph{``contravariant'' convention} that $\TT(\rho,\sigma)$ stands for the space of morphisms from $\sigma$ to $\rho$. %-- this makes multiplication formulas  more wieldy). 
 \par
 The construction of the single algebra $\OO_\rho$, with  $\rho$ fixed, relies only on the semigroup  $\{\rho^n\}_{n\in \N}$ generated by $\rho$ (we include $0$ in $\N$). In such a situation  our framework becomes even more transparent by  assuming that $\TT$ is simply a $C^*$-category with $\N$ as the set of objects, equipped with a   $^*$-functor $\otimes 1: \TT \to \TT$ such that $\otimes 1: \TT(n,m)\to \TT(n+1,m+1)$, $n,m\in \N$. Then it makes sense to denote the \emph{Doplicher-Roberts algebra} (associated to the object $1\in \N$) by $\DR(\TT)$. %
 Such algebras are sometimes also called DR-algebras \cite{v}.
 \par
A somewhat different but closely related and  important class of algebras form the $C^*$-algebras associated with  $C^*$-correspondences, the study of which   was initiated by Pimsner \cite{p}. More specifically, a  $C^*$-correspondence $X$ over a $C^*$-algebra $A$ (sometimes called a Hilbert bimodule)  is a right Hilbert $A$-module equipped with a left action $\phi$ of $A$ by adjointable operators, and  \emph{Cuntz-Pimsner algebra} $\OO_X$ is  constructed  as a  quotient of the \emph{ Toeplitz algebra} of $X$ generated by the \emph{Fock representation} of $X$ on the Fock module $\FF(X)=\bigotimes_{n=0}^\infty X^{\otimes n}$. Algebras arising in this way are known to comprise various  $C^*$-algebras found in the literature:  crossed products by  automorphisms, partial crossed products, crossed products by endomorphisms,  $C^*$-algebras of graphs (in particular Cuntz-Krieger algebras), Exel-Laca algebras, $C^*$-algebras of topological quivers, and many more. It has to be noted that originally Pimsner in his analysis  assumed that the \emph{left action $\phi$  on $X$ is injective}. However,  this seemingly only technical assumption turned out to be crucial. In particular, the efforts to remove this restriction  resulted in a variety of  approaches \cite{ms}, \cite{aee}, \cite{fmr},   \cite{fr},   \cite{katsura1}. We single out  two of these. Firstly, the $C^*$-algebra $\OO_X$ introduced by Katsura \cite{katsura1}  seems to be the most natural candidate for $\OO_X$ in the general case. It is  the smallest $C^*$-algebra among $C^*$-algebras generated by injective representations of $X$  admitting gauge actions, cf. \cite[Prop. 7.14]{katsura2}. Secondly, the so-called \emph{relative  Cuntz-Pimsner algebras} $\OO(J,X)$ of Muhly and Solel \cite{ms} possess  traits of being the most general, since by  particular choices of an ideal $J$ in $A$, one can cover  all the aforementioned constructions. Moreover,  $\OO(J,X)$ arise quite naturally,  when one tries to understand the ideal structure of Cuntz-Pimsner algebras \cite{fmr}, \cite{ms}, and when dealing with  certain concrete problems of description of $C^*$-algebras generated by irreversible dynamical systems \cite{Brow-Rae}, \cite{kwa} \cite{kwa-leb2}.
\par
 The relationship between  Pimsner's algebras and Doplicher-Roberts  algebras $\DR(X)$ associated with a $C^*$-correspondence $X$ was investigated 
in \cite{dpz}, \cite{fmr} where it was assumed, as in Pimsner's paper \cite{p}, that the \emph{left action $\phi$ on  $X$ is  injective}. It was noticed that  $\DR(X)$ is closely related but tends to be larger than $\OO_X$. Namely,   there are natural embeddings $\OO_X \subset \DR(X)\subset \OO_X^{**}$ and  the equality  $\OO_X=\DR(X)$ holds, for instance, if   $X$ is  finite projective,  cf. \cite[Prop. 3.2]{dpz}, \cite[Cor. 6.3]{fmr}.  In particular, if $\OO_X=\DR(X)$  the Doplicher-Roberts construction makes   the analysis of the Cuntz-Pimsner algebra $\OO_X$ very accessible and proves to be very useful, cf. \cite{kpw}. It should be stressed, however, that when $\phi$ is not injective
the relation between $\DR(X)$ and any of the algebras $\OO(J,X)$ is far more elusive and  remains practically untouched. In this article we develop a general approach to overcome these problems. %In particular, this may concern physically oriented mathematicians  interested in superselection structures.
\par
The main  observation is that even though the categorial language of  Doplicher and Roberts does not exactly fit into the formalism of Pimsner, in many aspects it is more natural  and being \emph{properly adapted} it  clarifies the relationship between  all the above-named constructions as well as shed much more  light  on their structures. To support this point of view let us  note that the algebras $\OO(J,X)$, $\DR(X)$,  $\DR(\TT)$ are spanned respectively by images of the spaces of "compact" operators, adjointable operators, and abstract morphisms (arrows):  
$$
\KK( X^{\otimes m},X^{\otimes n}), \qquad \LL( X^{\otimes m},X^{\otimes n}), \qquad \TT(n,m),\qquad n,m \in \N,
$$
where the above building blocks are "glued together" to form the corresponding algebras  in procedures based essentially on the "tensoring" -- either the natural tensoring on $X$,   or an abstract tensoring on $\TT$. 
Namely, if $A$ is unital,  the family  $\TT_X:=\{\LL( X^{\otimes m},X^{\otimes n})\}_{n,m\in \N}$, where $X^{\otimes 0}:=A$, with the tensoring on the right  by the identity operator on $X$  form a natural right tensor $C^*$-category. By  definition $\DR(X):=\DR(\TT_X)$ and the general   Doplicher-Roberts algebra $\DR(\TT)$ is a $C^*$-algebra  endowed with an action of the unit circle for which the $k$-spectral subspace  is an inductive limit of the  inductive sequence
$$
\TT(r+k,r) \stackrel{\otimes 1}{\longrightarrow}\TT(r+k+1,r+1) \stackrel{\otimes 1}{\longrightarrow} \TT(r+k+2,r+2) \stackrel{\otimes 1}{\longrightarrow} ...\,\,\, .
$$
In particular, we have natural homomorphisms $\iota_{n,m}:\TT(n,m)\to \DR(\TT)$ that form a representation  of the $C^*$-category $\TT$,  \cite{glr}, and clearly this representation determines the structure of $\DR(\TT)$. Similarly, the Fock representation (or any other universal representation) of $X$  give rise  to  homomorphisms  $\iota_{n,m}:\KK( X^{\otimes m},X^{\otimes n})\to\OO(J,X)$ which posses the same properties as the aforesaid representation of $\TT$.  It is evident and almost symptomatic that any analysis of $\OO(J,X)$ leads 
 to  an  analysis of the family  $\{\iota_{n,m}\}_{n,m\in \N}$ of such  representations. 
This traces the fact that within our unified approach  we may    transfer a very rich and well developed representation theory of relative Cuntz-Pimsner algebras \cite{ms}, \cite{fmr}, \cite{katsura}, \cite{katsura2}, onto the ground of Doplicher  and Roberts. In particular, we may "improve" the construction of $\DR(\TT)$ so that  the universal representation of $\TT$ in $\DR(\TT)$ is injective, even  when the right tensoring is not. On the other hand,   generalizing the inductive limit construction of $\DR(\TT)$, which makes its structure very accessible, allow us to clear up the description of  ideal structure of $\OO(J,X)$ obtained in \cite{fmr},  \cite{katsura}. 
\par
There are two plain but important new principles  that shine through  our development: 
\begin{itemize}
\item[1)] It is more natural to work with $C^*$-precategories, the categorical analogues of (not necessarily unital) $C^*$-algebras, rather than  with $C^*$-categories.
\item[2)] Not only right tensor $C^*$-precategories but also their ideals naturally give rise to universal $C^*$-algebras.
\end{itemize}
To support 1) note  that the unit-existence-requirement embedded into the notion of a category   causes  dispensable technicalities and lingual inconsequence like that  an ideal in a $C^*$-category, \cite[Def. 1.6]{glr}, may not be a $C^*$-category. Moreover, dealing with a $C^*$-correspondence $X$ over a non-unital $C^*$-algebra $A$ leads to  the following problem: If $\phi$ is not non-degenerate there is no obvious right tensoring on the $C^*$-category $\{\LL( X^{\otimes m},X^{\otimes n})\}_{n,m\in \N}$, as there is no such extension of $\phi:A \to \LL(X)$ up to the multiplier algebra $M(A)=\LL(A)$ of $A$. A natural solution to that problem is  considering  a  smaller $C^*$-precategory  $\TT_X$ where $\TT_X(0,0):=A$ (then $\TT_X$ is a $C^*$-category iff $A$ is unital).
\par
An important remark on account of 2) is that the $C^*$-precategory $\KK_X$ where $\KK_X(n,m):=\KK( X^{\otimes m},X^{\otimes n})$, $n,m\in \N$,  forms an ideal in the right tensor $C^*$-precategory $\TT_X$, but as a rule it is not  a right tensor $C^*$-precategory itself. Indeed, a typical situation is that the tensor product of a "compact" operator with the identity is   no longer "compact", so $\KK_X$ is not preserved under the tensoring $\otimes 1$. As we will show this  turns out not to be an  obstacle. This is due to the new feature of our  construction -- it applies well not only   to right tensor $C^*$-precategories but also to their ideals.
\\

The article is organized as follows.
We begin in Section \ref{preliminary monty} with a  review of basic facts and objects related to $C^*$-correspondences that will be of interest to us. In particular, we  introduce   examples of $C^*$-cor\-res\-pon\-den\-ces associated with  partial morphisms and directed graphs which we will  use throughout the paper to give a dynamical and  combinatorial interpretation of the theory presented. In Section \ref{section categories}, slightly modifying  and extending the terminology of \cite{glr}, we establish the rudiments of the theory of $C^*$-precategories.  A very useful statement here is Theorem  \ref{diagonal of ideals} which characterizes  ideals in $C^*$-precategories via ideals in $C^*$-algebras. When applied to the  $C^*$-precategory $\KK_X$ this  gives  a one-to-one  correspondence  between the ideals in $A$ and $\KK_X$, see Proposition \ref{Ideals in TT(X)}.
\par
One of the most important structures of our analysis -- the right tensor $C^*$-precategories and their representation theory is undertaken in Section \ref{section three}.  We provide a definition of a \emph{right tensor representation} as a representation of a  $C^*$-precategory which is compatible with the tensoring, and what is very important this notion makes sense not only for right tensor $C^*$-precategories but also for their ideals.  Proposition \ref{proposition:rho} states that such representations may be considered as generalizations of representations of  $C^*$-correspondences. Following this line of thinking  we give   a new meaning to the notion of coisometricity introduced by Muhly and Solel \cite{ms},   see Definitions   \ref{nie wiem jak to nazwac definition}, \ref{definition for ideals in categories}, and also Corollary \ref{coisometricicity of representations}.
\par
In Section \ref{section four} we present three different definitions, and establish their equivalence, of the \emph{$C^*$-algebra $\OO_{\TT}(\KK,\JJ)$ of an ideal $\KK$ in a right tensor $C^*$-precategory $\TT$ relative to an ideal $\JJ$}. We define  $\OO_{\TT}(\KK,\JJ)$  as: a universal $C^*$-algebra with respect to right tensor representations of $\KK$ coisometric on $\JJ$, Definition \ref{definition of main object}; an explicitly constructed algebra with explicit formulas for norm and algebraic operations, Subsection \ref{section construction}; and a $C^*$-algebra obtained via inductive limits formed from a specially constructed right tensor $C^*$-precategory $\KK_\JJ$, page  \pageref{dr-definition of our algebras}. Such a variety of points of view results in a numerous immediate interesting remarks. In particular, it allows us   to reveal the relationships between the algebras  $\OO_\TT(\KK,\JJ)$,  $\DR(\TT)$   and algebras admitting circle action, see  Section \ref{Doplichery i cirlce}. 
\par
 The fundamental tool in our analysis of the ideal structure of $\OO_\TT(\KK,\JJ)$ is  \emph{Structure Theorem}  (Theorem \ref{structure theorem}) which generalizes the main goal of \cite{fmr}.  It states that the ideal  $\OO(\NN)$ in $\OO_\TT(\KK,\JJ)$ generated by an invariant ideal $\NN$ in $\TT$ may be naturally identified as $\OO_\TT(\KK\cap \NN ,\JJ \cap \NN)$, and   the quotient $\OO_\TT(\KK,\JJ)/\OO(\NN)$ identifies as $\OO_{\TT/\NN}(\KK/\NN,\JJ/\NN)$. Moreover, we show that the  ideal $\NN$ may be replaced by its $\JJ$-saturation $\SSS_\JJ(\NN)$ (a concept that generalizes $X$-saturation \cite{mt} and negative invariance \cite{katsura2}), so that our Structure Theorem actually gives an embedding of the lattice of invariant,  $\JJ$-saturated ideals in $\TT$ into the lattice of gauge-invariant ideals in   $\OO_\TT(\KK,\JJ)$. Another application of the Structure Theorem establishes  procedures of reduction of relations  defining $\OO_\TT(\KK,\JJ)$, Definition \ref{reduction ideal definition}, Theorem \ref{reduction of relations theorem}. This  broadens and deepens  a topic started in \cite{kwa-leb1}, and  is indispensable in our further considerations. In   Subsection \ref{structure theorem for relative Cuntz-Pimsner algebras subsection} we discuss how Theorem \ref{structure theorem}
improves  \cite[Thm 3.1]{fmr} and why  a gauge-invariant ideal $\OO(I)$ in a relative Cuntz-Pimsner algebra $\OO(J,X)$ is "merely" Morita equivalent  to the corresponding relative Cuntz-Pimsner $\OO(J\cap I, XI)$.
\par
Our analogue of the \emph{gauge-invariant uniqueness theorem} is Theorem \ref{Gauge invariance theorem for O_T(I,J)1}. It extends the corresponding  theorems for relative Cuntz-Pimsner algebras \cite[Thm. 4.1]{fmr}, \cite[Thm. 5.1]{mt}, \cite[Thm. 6.4]{katsura}, \cite[Cor. 11.7]{katsura2}. %what we discuss in more  detail in Section \ref{Representations of algebras associated with X}. 
The main  novelty  is the use of the right tensor $C^*$-precategory $\KK_\JJ$ constructed in Theorem \ref{norm formulas on 0-spectral subspace}. In particular, we establish (in Theorem \ref{theorem on lattice structure}) a lattice isomorphism  that characterizes  the  \emph{gauge-invariant ideal structure} of  $\OO_\TT(\KK,\JJ)$ in terms of invariant $\KK_\JJ$-saturated ideals in $\KK_\JJ$.   Under certain additional assumptions, this result allows to obtain an analogous description  in terms of invariant  $\JJ$-saturated ideals in $\KK$ (Theorem \ref{theorem on lattice structure2}). The general relationship between invariant saturated ideals in  $\KK$ and $\KK_\JJ$ (and thereby also gauge-invariant ideals in  $\OO_\TT(\KK,\JJ)$) is complex;  on the level of relative Cuntz-Pimsner algebras it is completely revealed  in Theorem  \ref{Lattice structure of gauge invariant ideals in O(J,X) theorem} where   the role of $T$-pairs introduced in \cite{katsura2} is also clarified.
\par
%As an application of the developed theory we investigate representations of relative Cuntz-Pimsner  algebras $\OO(J,X)$ and relative Doplicher-Roberts algebras $\DR(J, X)$ defined in the previous section.
Aiming at a generalization of \cite[Thm. 6.6]{fmr},  \cite[Thm. 4.1]{dpz} %, which, in the case the left action $\varphi$ is injective and hence $\OO_X\subset \DR(X)$, identify the representations of $\OO_X$ that extend to $\DR(X)$,  
in Section \ref{Embedding theorems for O_T(I,J)} we  examine  the  general conditions assuring that algebras of type $\OO_\TT(\KK,\JJ)$  embeds  into one another. In this direction we establish two useful results, Propositions \ref{proposition of hierarchy of algebras2}, \ref{proposition of hierarchy of algebras}, which  are non-trivial, inequivalent generalizations of  \cite[Prop. 3.2]{dpz}, \cite[Cor. 6.3]{fmr}. Motivated by  these considerations we introduce an analogue of relative Cuntz-Pimsner algebras - \emph{relative Doplicher-Roberts algebras $\DR(J, X)$}, Definition \ref{relative Doplicher-Roberts algebra definition}. In particular, we give  necessary and sufficient conditions % (Corollary \ref{corollary about correspondences2}) 
under which the natural embedding 
$
\OO(J_0,X)\subset \DR(J, X)
$
holds. \par
In the final section, we describe representations of $\DR(J, X)$ that extends representations of $\OO(J_0,X)$ (Proposition  \ref{combining of two representations})  and  give criteria under which such a representation is  faithful, see  Theorem \ref{main theorem of this section}. Additionally  we show that  every faithful representation of $\OO(J_0,X)$ extends to faithful representation of  $\DR(J, X)$ for an appropriate $J$, see Theorem \ref{if this is  true theorem}. 
\par

%The author would like to thank the referee for perspicacious comments which helped to significantly improve  the text in several places. 
\begin{notation}
Following \cite{katsura2} we  denote by $\N=\{0,1,2,\ldots\}$ 
the set of natural numbers,  by $\C$ the set of complex numbers, and  by $S^1$ the group of complex numbers 
with absolute value $1$. 
We use a convention that 
$\gamma(A,B)=\{\gamma(a,b)\in D\mid a\in A,b\in B\}$ 
for a map $\gamma\colon A\times B\to D$ 
such as inner products, multiplications or representations. 
We denote by  $\spane\{\cdots\}$ a linear spans of $\{\cdots\}$, and by  $\clsp\{\cdots\}$ 
the closure of $\spane\{\cdots\}$.
\end{notation} 
\section{Preliminaries on $C^*$-correspondences} \label{preliminary monty}
We  adopt the standard notations and definitions of objects related to (right) Hilbert $C^*$-modules, cf. \cite{lance}, \cite{rw}. In particular, we denote by  $X$ and $Y$  Hilbert modules over a $C^*$-algebra $A$;  $\LL(X,Y)$ stands for the space of adjointable operators from $X$ into $Y$; and $\KK(X,Y)$ is  the  space  of "compact" operators in $\LL(X,Y)$, that is $\KK(X,Y)=\clsp\{\Theta_{y,x}: x\in X,y\in Y\}$ where $\Theta_{y,x}(z)=y\langle x,z\rangle_A$, $z\in X$. 
If $I$ is an ideal in $A$ (by which we  always mean a closed two-sided ideal), then $XI$ is  both a Hilbert $A$-submodule of $X$ and  a  Hilbert $I$-module, as  we  have 
\begin{equation}\label{XI equation}
 XI=\{xi:x\in X,\, i\in I\}=\{x\in X: \langle x,y\rangle_A\in I \textrm{ for all }y\in X\}, 
 \end{equation}
cf. \cite[Prop. 1.3]{katsura2}. A natural identification of $\KK(XI)$ as subalgebra of $\KK(X)$ was used in  \cite{katsura2}, \cite{fmr}. Actually we have
\begin{lemma}\label{identifying lemma0}
The equalities 
 \begin{align*}
 \clsp\{\Theta_{y,x}: x\in X,y\in YI\}& =\clsp\{\Theta_{y,x}: x\in XI,y\in Y\} \\
&=\clsp\{\Theta_{y,x}: x\in XI,y\in YI\}
 \end{align*}
establish natural identifications 
$$
\KK(X,YI)=\KK(XI,Y)=\KK(XI,YI)\subset \KK(X,Y).
$$
\end{lemma}
\begin{proof}
Clear by definitions, cf. \cite[p. 9]{lance}, \cite{katsura2}.
\end{proof}
The above identifications do not carry  over to  general adjointable maps  due to a possible non-existence of adjoint operators. In order to clarify this situation, let us note that if $I$ is an ideal in $A$, then (by \eqref{XI equation}) for $a\in \LL(X,Y)$ the two conditions $a(X)\subset YI$ and $a^*(Y)\subset XI$ are equivalent. We denote by $\LL_I(X,Y)$ the space of all adjointable maps $a\in \LL(X,Y)$ satisfying these   equivalent conditions. % is clarified by the following lemma.
\begin{lemma}\label{identifying lemma}
The inclusions of sets of mapping and restriction of mappings yield  six natural linear homomorphisms between the space of morphisms, presented on the diagram
$$\footnotesize 
\xymatrix{             &    \LL(XI,YI)  &    \\
\LL(XI,Y)   \ar[ru]^{incls.}  &   \LL(X,YI)  \ar[u]^{restr.} &  \LL(X,Y) \ar@{.>}[lu]_{restr.} \\
 &  \ar[lu]^{restr.} \LL_I(X,Y) \ar[u]^{incls.}\ar[ru]_{incls.} &
}
$$
where
\begin{itemize}
\item[a)] the three homomorphisms from $\LL_I(X,Y)$ to $\LL(XI,Y)$, $\LL(X,YI)$,  $\LL(X,Y)$ are injective, 
\item[b)] homomorphisms  from  $\LL(XI,Y)$ and $\LL(X,YI)$ to $\LL(XI,YI)$ are injective and the one from
  $\LL(X,Y)$ to $\LL(XI,YI)$  in general is not.
  \end{itemize}
Moreover, 
\begin{itemize}
\item[i)] the three intersections in $\LL(XI,YI)$ of any  two images of the three homomorphisms in b) coincide,
\item[ii)] the three maps from $\LL_I(X,Y)$ to $\LL(XI,YI)$, through $\LL(XI,Y)$, $\LL(X,YI)$, $\LL(X,Y)$, coincide, and this common map is an isomorphism onto the common intersection in i).
  \end{itemize}
%\begin{equation}\label{intersection equality}
%\LL(X,YI)\cap \LL(X,Y)= \LL(XI,Y)\cap \LL(X,YI).
%\end{equation}
\end{lemma}
 %\begin{figure}[htb] \begin{center}\setlength{\unitlength}{1mm} \footnotesize $$ \xymatrix{             &    \LL(XI,YI)  &    \\\LL(XI,Y)   \ar[ru]^{incls.}  &   \LL(X,YI)  \ar[u]^{restr.} &  \LL(X,Y) \ar@{.>}[lu]_{restr.} \\ &  \ar[lu]^{restr.} \LL_I(X,Y) \ar[u]^{incls.}\ar[ru]_{incls.} &}$$ \end{center} \caption{Natural homomorphisms given by an inclusion or restrictions\label{hasse diagram}} \end{figure}
 
\begin{proof} To see item a) note that if $a\in \LL_I(Y,X)$, then $a$ treated as a mapping is an element of $\LL(X,YI)$ with the adjoint given by restriction of $a^*\in \LL(Y,X)$  to $YI$. Moreover,  the map $\LL_I(X,Y)\hookrightarrow \LL(XI,Y)$ given by restriction of mappings is an isometric homomorphism since it can be obtained by passing to adjoints  in the inclusion $\LL_I(X,Y)\subset \LL(X,YI)$. To prove item b) one can argue in a  similar way. 
\\
For item i)  note that $a\in \LL(XI,YI)$ is in the intersection of the images of $\LL(XI,Y)$ and   $\LL(X,YI)$  in  $\LL(XI,YI)$  if and only if 
there are  $\a\in \LL(X,YI)$ and  $\widetilde{a^*}\in \LL(Y,XI)$ such that $\widetilde{a}|_{XI}=a$  and   $\widetilde{a^*}|_{YI}=a^*$, but then it follows that $\widetilde{a}\in \LL_I(X,Y)$ where $\widetilde{a}^*=\widetilde{a^*}$. Conversely, for any  $\widetilde{a}\in \LL_I(X,Y)$ its restriction $a=\widetilde{a}|_{XI}$ yields an element  lying in the intersection of the  images of $\LL(XI,Y)$ and   $\LL(X,YI)$  in  $\LL(XI,YI)$.  Similarly one sees, that the reamining two intersections in item i)  coincide and consist of  morphism from $\LL(XI,YI)$ possessing adjointable extensions to elements of $\LL(X,Y)$ which then necessarily lie in $\LL_I(X,Y)$. This  proves i) and ii).
\end{proof}
\begin{remark} All five  injective homomorphisms in the   above lemma become inclusions when viewing $X$ and $Y$ as Hilbert modules in a $C^*$-algebra $\M$, cf. \cite[Prop. 2.1]{dpz}. In general, the six considered maps are not surjective, and the three images of $\LL(XI,Y)$, $\LL(X,YI)$  and $\LL(X,Y)$  in  $\LL(XI,YI)$  are incomparable.
\end{remark}
 An  important fact  for our purposes  is  

\begin{proposition}\label{corollary on ideals in adjointable maps}
The $C^*$-algebra $\KK(XI)$ is an ideal in $\LL_I(X)$ which  in turn is an ideal  in the $C^*$-algebra $\LL(X)$.
\end{proposition}

For an ideal $I$ in a $C^*$-algebra $A$ we may consider the quotient  space $X/XI$ as a Hilbert $A/I$-module with an $A/I$-valued inner product and right action of $A/I$ given by
$
\langle q(x), q(y)\rangle_{A/I}:=q (\langle x, y\rangle_A)$, $q(x)q(a):=q(xa)
$
where $q$ denotes  both the quotient maps $A\to A/I$ and $X\to X/XI$, cf. \cite[Lem. 2.1]{fmr}.
Moreover, we have a natural map $q:\LL(X)\to \LL(X/XI)$ where  $q(a)q(x)=q(ax)$ for $a\in \LL(X)$ and $x\in X$.
\begin{lemma}
The kernel of the map $q:\LL(X)\to \LL(X/XI)$ is $\LL_I(X)$ and the restriction of $q$ to $\KK(X)$ is a surjection onto $\KK(X/XI)$ whose kernel is $\KK(XI)$.
\end{lemma}
\begin{proof} See \cite[Lem. 1.6]{katsura2} and remarks preceding this statement.
\end{proof}
\begin{corollary}[Lem. 2.6 \cite{fmr}]\label{corollary to quotients}
We have a natural isomorphism  $\KK(X)/\KK(XI)\cong \KK(X/XI)$. 
\end{corollary}
For  Hilbert modules with  left action we use the term  $C^*$-correspondence, and we reserve the term Hilbert bimodule  for an object with an additional structure (cf. Definition \ref{bimodule definition} and Proposition \ref{equivalent characteristics of completeness}  below) --  this convention  seems to become standard.  
\begin{definition}
\emph{A $C^*$-correspondence $X$ over a $C^*$-algebra $A$} is  a (right) Hilbert $A$-module equipped with a $^*$-homomorphism $\phi:A \to \LL(X)$. We refer to $\phi$ as the left action of  the $C^*$-correspondence $X$ and  write $a\cdot x := \phi(a)x$.
\end{definition}
Let us fix  a $C^*$-correspondence $X$ over a $C^*$-algebra $A$ and   a Hilbert $A$-module $Y$.  There is a naturally defined tensor product Hilbert $A$-module  $Y\otimes X$, cf. \cite{lance}, \cite{rw} or \cite{katsura}. %In particular for $n\geq 1$ the $n$-fold tensor product $X^{otimes n}=X\otimes .... X
 An ideal $I$ in $A$  is called \emph{$X$-invariant} if $\varphi(I)X\subset XI$, and for such an ideal the quotient $A/I$-module $X/XI$ with right action $q(a)q(x)=q(\varphi(a)x)$ becomes a $C^*$-correspondence  over $A/I$, cf.  \cite[Lem. 2.3]{fmr}, \cite{katsura2}  \cite{kpw}. In particular,  we may consider  two Hilbert $A/I$-modules  $Y/YI\otimes X/ XI$ and  $(Y\otimes X)/ (Y\otimes XI)$.
\begin{lemma}\label{lemma to quotients}
For an $X$-invariant ideal $I$ in $A$   we have a natural isomorphism of Hilbert modules
$$
(Y/YI)\otimes (X/XI)\cong (Y\otimes X)/ (Y\otimes XI).
$$
\end{lemma}
\begin{proof}
Let $x_1,x_2\in X$, $y_1,y_2\in Y$ and $i,j\in I$. Then 
$$
 (y_1+y_2j)\otimes (x_1+x_2i)=y_1\otimes x_1 + \Big(y_1\otimes x_2i + y_2\otimes\varphi(j)x_1+ y_2j\otimes x_2i\Big)
$$
where by $X$-invariance of $I$ the term in brackets belongs to $Y\otimes XI$. Thereby the mapping
$$
(y+ YI)\otimes (x+XI)\longmapsto (y\otimes x) + (Y\otimes XI)
$$
is well defined. Clearly, it is surjective and $A/I$-linear.  The simple calculation
$$
\langle q(y_1)\otimes q(x_1), q(y_2)\otimes q(x_2)\rangle_{A/I}=\langle q(x_1), q(\varphi\langle y_1, y_2 \rangle_{A}) q(x_2)\rangle_{A/I}
$$
$$
=q(\langle x_1, \varphi(\langle y_1, y_2 \rangle_{A}) x_2\rangle_{A})
=q(\langle y_1\otimes x_1,  y_2 \otimes x_2\rangle_{A})=\langle q(y_1\otimes x_1), q(y_2\otimes x_2)\rangle_{A/I}
$$
shows that  the above mapping preserves $A/I$-valued inner products and hence is isometric.
\end{proof}
We have a  homomorphism $\LL(Y)\ni a \to a\otimes 1 \in \LL(Y\otimes X)$ where 
\begin{equation}\label{right tensoring of operators definition}
 (a\otimes 1)(y\otimes x):=ay\otimes x,\qquad x\in X,\,\, y\in Y.  
\end{equation}
The properties of this homomorphism are related to objects that will play important role throughout this paper. We define  
 $$
 J(X):=\phi^{-1}(\KK(X))
 $$
 which is an ideal in $A$. If $J$ is an ideal in  $A$  we define 
$$
J^\bot=\{a\in A: aJ=\{0\}\}
$$
which is also an ideal in $A$  called an \emph{annihilator} of $J$. It is a unique  ideal in $A$ such that $J^\bot\cap J=\{0\}$ and for any ideal $I$ in $A$ we have  $
 J \cap I =\{0\} \,\, \Longrightarrow I \subset J^\bot
 $.
\begin{lemma}\emph{(cf. \cite[Lem. 4.2]{fmr}).}\label{lemma embedding of tensors}
Suppose that $X$ is a $C^*$-correspondence over $A$ and  $Y$ is a right Hilbert $A$-module. 
\begin{itemize}
\item[i)] 
The map $ a\mapsto a\otimes 1$ restricted to  $\LL_{(\ker\phi)^\bot}(Y)$  is isometric.
\item[ii)] If $\{\mu_\lambda\}_\lambda$ is an approximate unit for $\KK(Y)$, then  $\mu_\lambda\otimes 1$ converges strictly to $1\in \LL(Y\otimes X)$. In particular, $(\KK(Y)\otimes 1)\KK(Y,X)=\KK(Y,X)$.
\item[iii)] If $a\in \LL_{(\ker\phi)^\bot}(Y)$ and $a\otimes 1\in \KK(Y\otimes X)$, then  $a\in \KK(Y)$.
\item[iv)] If $a\otimes 1\in \KK(Y\otimes X)$ and $a\in \KK(Y)$, then  $a\in \KK(Y J(X))$.
\end{itemize}
\end{lemma}
\begin{proof}
i) Let $a\in \LL_{(\ker\phi)^\bot}(Y)$. It suffices to prove $\|a\|\leq \|a\otimes 1\|$, and for that purpose take an arbitrary $y\in Y$. Since $\langle ay,ay\rangle_A \in (\ker\phi)^\bot$ and $\phi$ is isometric on $(\ker\phi)^\bot$ we have $\|\phi(\langle ay,ay\rangle_A)\|=\|ay\|^2$. Positivity of $\phi(\langle ay,ay\rangle_A)$ implies that  for each $\varepsilon >0$,  there exists $x\in X$ such that   $\|x\|=1$ and
$$
\|\langle x,\phi(\langle ay,ay\rangle_A)x\rangle_A\|\geq \|\phi(\langle ay,ay\rangle_A)\| -\varepsilon= \|ay\|^2 -\varepsilon.
$$ 
Thus 
$$
\|(a\otimes 1 )(y\otimes x)\|^2=\|\langle x,\phi(\langle ay,ay\rangle_A)x\rangle_A\|\geq \|ay\|^2 -\varepsilon
$$
and as $\|(y\otimes x)\|\leq \|y\|$ we get $ \|a\otimes 1 \| \geq \|a\|$. 
\\
ii) See the proof of \cite[Lem. 4.2 (2)]{fmr}.
\\
iii) Let $\{\mu_\lambda\}_\lambda$ be an approximate unit for $\KK(Y)$. By item ii) we get
$$
0=\lim_\lambda \|a\otimes 1 - (\mu_\lambda\otimes 1 )(a\otimes 1 )\|=\lim_\lambda \|a - \mu_\lambda a\|.
$$
Hence $a$ is "compact". 
\\
iv) See the second part of the proof of  \cite[Lem. 4.2 (2)]{fmr}.
 \end{proof}
Following  \cite{katsura} we clarify the relationship between $C^*$-correspondences and Hilbert bimodules.
 \begin{definition}\label{bimodule definition}
  We say that $X$ is a \emph{Hilbert $A$-bimodule} if it is at the same time a Hilbert left $A$-module and a Hilbert right $A$-module  with sesqui-linear forms  ${_A\langle} \cdot , \cdot  \rangle$ and $\langle \cdot  , \cdot  \rangle_A$ 
related via the so-called imprimitivity condition:
\begin{equation}\label{complete bimodules2}
 x \cdot \langle y ,z \rangle_A = {_A\langle} x , y  \rangle \cdot z, \qquad \textrm{for all}\,\,\, x,y,z\in X.
 \end{equation}  \end{definition}
It follows from  \eqref{complete bimodules2} that left action in a Hilbert bimodule acts by adjointable maps and hence every Hilbert $C^*$-bimodule is a $C^*$-correspondence, cf. \cite[3.3]{katsura1}. In the converse direction we have
  \begin{proposition}\label{equivalent characteristics of completeness} Let $X$ be a $C^*$-correspondence. The following conditions are equivalent:
 \begin{itemize}
 \item[i)] $X$ is a Hilbert $A$-bimodule (that is there exists a sesqui-linear form  ${_A\langle} \cdot , \cdot  \rangle$ for which $X$ together with left action $\phi$ becomes a left Hilbert $A$-module satisfying \eqref{complete bimodules2}),
\item[ii)] there is a  function  $_A\langle \cdot ,\cdot \rangle:X\times X \to (\ker \phi)^\bot$  such that 
\begin{equation}\label{complete bimodules1}
 \phi({_A\langle} x , y  \rangle)=\Theta_{x,y}, \qquad x,y\in X, 
  \end{equation}  \item[iii)]  the mapping $\phi:(\ker\phi)^\bot \cap J(X) \to \KK(X)$ is onto (hence it is automatically an isomorphism). 
 \end{itemize}
The objects in i) and ii) are determined uniquely,   the function $_A\langle \cdot ,\cdot \rangle$ from  item $ii)$ coincides with the inner product from  item $i)$ and  
\begin{equation}\label{complete bimodules3}
{_A\langle} x , y  \rangle = \phi^{-1}(\Theta_{x,y}), \qquad x,y\in X,
 \end{equation}
where $\phi^{-1}$ is the inverse to the isomorphism $\phi:(\ker\phi)^\bot \cap J(X) \to \KK(X)$.
 Moreover denoting $\overline{{_A\langle} X , X  \rangle}:=\clsp\{{_A\langle} x , y  \rangle: x,y\in X\}$ we have 
 $$
\overline{{_A\langle} X , X  \rangle}=(\ker\phi)^\bot \cap J(X).
 $$
 \end{proposition}
 \begin{proof}  i)$\Rightarrow$ ii) One easily sees that conditions \eqref{complete bimodules1} and \eqref{complete bimodules2} are equivalent. Since $\phi$ is injective on $\overline{{_A\langle} X , X  \rangle}$ and $\overline{{_A\langle} X , X  \rangle}$ is an ideal we have $\overline{{_A\langle} X , X  \rangle}\subset (\ker \phi)^\bot$. Consequentely ii) holds.
 \\
 ii)$\Rightarrow$ iii) By \eqref{complete bimodules1},  $\phi$  maps the set $\overline{{_A\langle} X , X  \rangle}$ onto $\KK(X)$ and $\overline{{_A\langle} X , X  \rangle} \subset (\ker\phi)^\bot \cap J(X)$.  As  $\phi$ is injective on $(\ker\phi)^\bot$ we see that $\phi:(\ker\phi)^\bot \cap J(X) \to \KK(X)$ is an isomorphism. In particular,   $\overline{{_A\langle} X , X  \rangle}=(\ker\phi)^\bot \cap J(X)$ and \eqref{complete bimodules3} holds.\\
 iii)$\Rightarrow$ i) It is straightforward to check that 
 $_A\langle \cdot ,\cdot \rangle$  defined by \eqref{complete bimodules3} is an inner product for the left $A$-module $X$. 
  As we already noted, \eqref{complete bimodules1} and \eqref{complete bimodules2} are equivalent, 
 and thus i) holds.
 \end{proof}

 \begin{example}[$C^*$-correspondence of a partial morphism]\label{partial morphism category1.0} By a \emph{partial morphism of a $C^*$-algebra} $A$ we mean  nondegenerate
 a $^*$-homomorphism $\varphi:A\to M(A_0)$ from $A$ to the multiplier algebra $M(A_0)$ of a hereditary subalgebra $A_0$ of $ A$,  cf. \cite{katsura1}. We recall that $\varphi:A\to M(A_0)$ is said to be \emph{nondegenerate} if $\varphi(A)A_0=A_0$. We construct a $C^*$-correspondence $X_\varphi$ from  $\varphi$ in the following way. We let
 $X_\varphi:=A_0A$
and put  
$$
a \cdot x :=\varphi(a)x,\quad x\cdot a:= xa,\quad \textrm{ and}\quad \langle x,y\rangle_A:=x^*y,  
$$
 where $a\in A$, $x,y\in X_\varphi$. Then $J(X_\varphi)=\varphi^{-1}(A_0)$.
 To assert when $X_\varphi$ is a Hilbert bimodule we slightly extend  R. Exel's  definition \cite[Def. 3.1]{exel1} and by  a  \emph{partial automorphism} of $A$ we shall mean a triple $(\theta,I,A_0)$ consisting of an ideal $I$  in $A$, a hereditary subalgebra $A_0$ of $A$ and an isomorphism  $\theta:I\to A_0$. A partial automorphism $(\theta,I,A_0)$ give rise to a partial morphism $\varphi:A\to M(A_0)$ via the formula
 $ \varphi(a)b:=\theta(a\theta^{-1}(b))$, $a\in A$, $b\in A_0$,  and then  we have 
 $
 I=(\ker\varphi)^\bot\cap  \varphi^{-1}(A_0)
 $. Conversely, if  $\varphi:A\to M(A_0)$ is a partial morphism such that $\varphi$ restricted to $
 I:=(\ker\varphi)^\bot\cap  \varphi^{-1}(A_0)$ is an isomorphism onto $A_0$, then $\varphi$ arises from a partial automorphism $(\theta,I,A_0)$ where $\theta:=\varphi|_I$.
 \\
 Therefore, in view of Proposition \ref{equivalent characteristics of completeness}, a $C^*$-correspondence $X_\varphi$ is a Hilbert bimodule if and only if $\varphi$ arises from a partial automorphism $\theta$, and then the "left" inner product is given by
 $$
 _A\langle x,y\rangle:=\theta^{-1}(xy^*).
 $$ 
   We indicate that  every endomorphism $\al:A\to A$ of a $C^*$-algebra $A$  may be treated as a partial morphism where $A_0=\al(A)A\al(A)$. In this sense a partial morphism $\varphi$ is an endomorphism iff $\varphi^{-1}(A_0)=A$. We shall denote a $C^*$-correspondence arising from an endomorphism $\al$ by $X_\al$. It is shown in \cite[Prop. 1.9]{kwa-trans} that if $A$ is unital, then  $X_\al=\al(1)A$ is a Hilbert bimodule iff there exists  a \emph{complete transfer operator} for  $\al$, which is  a bounded, positive  linear map $\LL:A\to A$ 
 such that 
% \begin{equation}
$$
\LL(\al(a)b) =a\LL(b),\qquad \textrm{ and }\qquad \al (\LL(1)) = \al (1), \qquad a,b\in A,
$$
see \cite{kwa-leb3} ($\LL$ exists iff the kernel of $\al$ is unital and the range is hereditary, see \cite{kwa-trans}). If this is the case  the "left" inner product is given by
 $$
 _A\langle x,y\rangle:=\LL(xy^*).
 $$  \end{example}
 \begin{example}[$C^*$-correspondence of a  graph]\label{precategory of a directed graph ex1} Suppose $E=(E^0,E^1,r,s)$ is a  \emph{directed graph} with vertex set $E^0$, edge set $E^1$, and $r,s:E^1\to E^0$ describing the range and the source of edges. 
A \emph{$C^*$-correspondence $X_E$ of the graph }$E$ is  defined in the following manner, cf. \cite{fr}, \cite{fmr}, \cite{mt}. The space
$X_E$ consists of functions
$x:E^1\to \C$ for which 
$$
v\in E^0 \longmapsto \sum_{\{e\in E^1: r(e) = v\}} |x(e)|^2
$$
belongs to $A := C_0(E^0)$, and $X_E$ is a $C^*$-correspondence over $A$ with the operations
\begin{align*}
(x\cdot a)(e) &:= x(e)a(r(e)) \ \mbox{ for } e\in E^1,\\
\langle x,y\rangle_A (v) &:= \sum_{\{e\in E^1: r(e) = v\}}
\overline{x(e)}y(e) \ \mbox{ for }v\in E^0, \mbox{ and}\\
(a\cdot x)(e) &:= a(s(e))x(e) \ \mbox{ for } e\in E^1.
\end{align*}
We note that $X_E$ and $A$ are respectively spanned by the point masses $\{\delta_f: f\in E^1\}$ and  $\{\delta_v: v\in E^0\}$. In particular 
$$
(\ker\phi)^\bot=\clsp\{\delta_v:s^{-1}(v)\neq \emptyset\},\qquad  J(X_E)=\clsp\{\delta_v:|s^{-1}(v)|< +\infty\},
$$
cf. for instance \cite{mt}. Moreover, if $f,g\in  E^1$ and $v\in E^0$ emits finitely many edges, then
$$
\phi(\delta_v)=\sum_{\{e\in E^1:s(e)=v\}}\Theta_{\delta_e,\delta_e},\qquad \textrm{and}\qquad \Theta_{\delta_f,\delta_g}\neq 0\, \Longleftrightarrow \, r(f)=r(g).
$$
It follows that $X_E$ is a Hilbert bimodule iff every vertex of the graph $E$ emits and receives at most one edge (equivalently maps $r,s$ are injective). If this is the case   
the left and right inner products are zero on the complement of  $r(E^1)$ and  $s(E^1)$, respectively, and 
$$
\langle x,y\rangle_A (r(e))=\overline{x(e)}y(e),\qquad _A\langle x,y\rangle (s(e))=x(e)\overline{y(e)}.
$$ 
The above formulas and objects naturally generalize to  \emph{topological graphs} \cite{katsura4}, \cite{katsura1}, that is quadruples
$E= (E^0,E^1, d, r)$ where $E^0$, $E^1$ are locally compact spaces, $d:E^1 \to E^0$ is a local
homeomorphism and $r: E^1 \to  E^0$ is a continuous map. The $C^*$-correspondence  $X_E$ of the topological graph $E$ is defined by the same formulas as above, with $r$ replaced with $d$ and $s$ with $r$. It is a Hilbert bimodule iff $d$, $r$ are homeomorphism onto open sets and then  the bimodule $X_E$ may treated as a bimodule $X_\varphi$ arising from a partial morphism $\varphi:C_0(E^0)\to C_b(r(E^1))$ given by the partial  homeomorphism $r \circ d^{-1}:d(E^1)\to r(E^1)$.   
\end{example}

\section{$C^*$-precategories and their ideals}\label{section categories}
For a notion of a precategory we adopt  the standard definition of  a category with the only exception that we drop the assumption of existence of identity morphisms. In the present paper we shall be interested in precategories with the class of objects being the set of  natural numbers $\N$. However, for future reference and to underscore the role of categorical language,   in this  section (and only in this section) we shall deal with general precategories. 
\begin{definition}
A \emph{precategory }$\TT$ consists of 
a class  of objects, denoted here by $\rho,\sigma,  \tau$, etc.; 
a class $\{\TT(\sigma,\rho)\}_{\sigma,\rho\in\TT}$ of disjoint classes of morphisms (arrows)   where  $\TT(\sigma,\rho)$ stands for the space of morphisms from $\rho$ to $\sigma$; and  a
composition of morphisms  
$
\TT(\tau,\sigma)\times \TT(\sigma,\rho)\ni(a, b)\to ab\in \TT(\tau,\rho)
$,  which is associative, in the sense that 
$
(a b) c=a (b c)
$  whenever the compositions of morphisms $a$, $b$, $c$ are allowable. 
\end{definition}
One  can always equip (if necessary) the sets of morphisms $\TT(\sigma, \sigma)$ with identity morphisms  in such a way that a given  precategory $\TT$ becomes a category. We  generalize the notion of a $C^*$-category, cf. \cite{glr}, \cite{dr}, in an obvious fashion.
\begin{definition}
A precategory $\TT$  is a \emph{$C^*$-precategory} if each set of morphisms $\TT(\sigma,\rho)$ is a complex Banach space, the composition of morphisms gives us a bilinear map 
$$
\TT(\tau,\sigma)\times \TT(\sigma,\rho)\ni(a, b)\to ab\in \TT(\tau,\rho)
$$ 
with $\|ab\|\leq \|a\|\cdot \|b\|$, and  there is an antilinear involutive contravariant functor $^*:\TT\to\TT$ such that if $a \in   \TT(\sigma,\rho)$, then $a^* \in   \TT(\rho,\sigma)$ and the $C^*$-equality $\|a^* a\|=\|a\|^2$ holds. A $C^*$-precategory $\TT$ where  $\TT$ is a category is a \emph{$C^*$-category}.
\end{definition}
\begin{notation}
If $\SSS$ is a \emph{sub-$C^*$-precategory} of $\TT$, that is   if  $\SSS$ and $\TT$ are two $C^*$-precategories  such that each space $\SSS(\sigma,\rho)$ is a closed subspace of $\TT(\sigma,\rho)$, we shall briefly write $\SSS\subset \TT$. If $\SSS$ and $\TT$ are two sub-$C^*$-precategories of another $C^*$-precategory, we denote by $\TT\cap \SSS$ the $C^*$-precategory where $(\SSS\cap \TT)(\sigma,\rho):=\SSS(\sigma,\rho)\cap \TT(\sigma,\rho)$. 
\end{notation}

Each space of morphisms $\TT(\rho,\rho)$ in a  $C^*$-precategory $\TT$ is a $C^*$-algebra, and by the $C^*$-equality  the functor "$^*$" is isometric on every space $\TT(\sigma,\rho)$. %using the $C^*$-equality 
A $C^*$-precategory $\TT$ is a $C^*$-category if and only if every $C^*$-algebra $\TT(\rho,\rho)$, $\rho \in \TT$, is unital.  In general, by adjoining  units to  $C^*$-algebras  $\TT(\rho,\rho)$, $\sigma \in \TT$, one may obtain  "unitization" of $\TT$, which is a $C^*$-category that contains $\TT$ as an ideal  in the sense of the following definition, cf. \cite[Def. 1.6]{glr}.
\begin{definition}
By a (closed two-sided) \emph{ideal $\KK$ in a $C^*$-precategory} $\TT$ we shall mean a collection  of Banach subspaces $\KK(\sigma,\rho)$ in  $\TT(\sigma,\rho)$,  $ \rho,\sigma \in \TT$, such that 
$$
  \TT(\tau,\sigma)\KK(\sigma,\rho) \subset \KK(\tau,\rho)\quad \textrm{ and } \quad \KK(\tau,\sigma)\TT(\sigma,\rho)\, \subset \KK(\tau,\rho),
  $$
  for any $\sigma,\rho,\tau\in \TT.
$ 
\end{definition}

Arguing as for $C^*$-algebras, cf. \cite[Prop. 1.7]{glr}, one shows that an ideal $\KK$ in a $C^*$-precategory $\TT$ is "self-adjoint" (hence it is a $C^*$-precategory). Obviously, each space $\KK(\rho,\rho)$   is a closed two-sided  ideal in the  $C^*$-algebra $\TT(\rho,\rho)$. A useful fact is that    $\KK$  is uniquely determined by these "diagonal ideals". In order to prove this, we  apply the following simple lemma.
\begin{lemma}\label{lemma about approximate units}
Let $\TT$ be a $C^*$-precategory and   let $
a\in \TT(\sigma,\rho)$. Let $\KK(\rho,\rho)$, $\KK(\sigma,\sigma)$ be arbitrary ideals in  $\TT(\rho,\rho)$ and $\TT(\sigma,\sigma)$ respectively, and  let $\{\mu_\lambda\}_\lambda$, $\{\nu_\lambda\}_\lambda$ be approximate units in  $\KK(\rho,\rho)$ and $\KK(\sigma,\sigma)$ respectively. Then 
$$
a^*a\in \KK(\rho,\rho)\,\,\, \Longleftrightarrow \,\,\, \lim_\lambda  a\mu_\lambda =a,
$$
$$
aa^*\in \KK(\sigma,\sigma)\,\,\, \Longleftrightarrow \,\,\, \lim_\lambda  \nu_\lambda a=a.
$$
\end{lemma}%\label{lemma taka tuka}
\begin{proof} If $\lim_\lambda  a\mu_\lambda =a$, then  $\lim_\lambda  a^*a\mu_\lambda =a^*a$ is in $\KK(\rho,\rho)$. Conversely, if $a^*a\in \KK(\rho,\rho)$ then 
 $$
\|a -a\mu_\lambda \|^2= \|(a -a\mu_\lambda)^*(a-a\mu_\lambda)\|\leq \|a^*a-a^*a\mu_\lambda\|+ \|\mu_\lambda\|\cdot  \|a^*a -a^*a \mu_\lambda\|,   
$$
which implies that $ a\mu_\lambda$ converges to $a$. The second equivalence can be proved analogously.
\end{proof}
\begin{theorem}[Characterization of ideals in $C^*$-precategories]\label{diagonal of ideals}
If $\KK$ is an ideal in a $C^*$-precategory $\TT$, then for every objects $\sigma$ and $\rho$
\begin{equation}\label{ideal defining formula}
\KK(\sigma,\rho)=\{a\in \TT(\sigma,\rho):a^*a \in \KK(\rho,\rho)\}=\{a\in \TT(\sigma,\rho):aa^* \in \KK(\sigma,\sigma)\}.
\end{equation}
 Conversely, if   $\{\KK(\rho,\rho)\}_{\rho\in \TT}$ is  a class of ideals $\KK(\rho,\rho)$  in $\TT(\rho,\rho)$, $\sigma\in \TT$, such that the equality  
\begin{equation}\label{ideal defining condition}
\{a\in \TT(\sigma,\rho):a^*a \in \KK(\rho,\rho)\}=\{a\in \TT(\sigma,\rho):aa^* \in \KK(\sigma,\sigma)\}
\end{equation}
holds for every $\sigma, \rho\in \TT$, then relations  \eqref{ideal defining formula} define an ideal $\KK$ in $\TT$.
\end{theorem}
\begin{proof}
Let $\KK$ be an ideal in $\TT$ and let  $a\in \TT(\sigma, \rho)$. If  $a \in \KK(\sigma,\rho)$ then  $a^*a \in \KK(\rho,\rho)$. Conversely, if $a^*a \in \KK(\rho,\rho)$, then using  Lemma \ref{lemma about approximate units} one gets $a \in \KK(\sigma,\rho)$. Thus   $\KK(\sigma,\rho)=\{a\in \TT(\sigma,\rho):a^*a \in \KK(\rho,\rho)\}$, and the equality $\KK(\sigma,\rho)=\{a\in \TT(\sigma,\rho):aa^* \in \KK(\sigma,\sigma) \}$ can be proved analogously.
\\
To prove the second part of assertion we fix a class of ideals  $\{\KK(\rho,\rho)\}_{\sigma\in \TT}$ such that \eqref{ideal defining condition} holds and use \eqref{ideal defining formula} to define $\KK=\{\KK(\sigma,\rho)\}_{\sigma,\rho\in \TT}$. If  $a\in \KK(\sigma,\rho)$, then   $a^*a\in \KK(\rho,\rho)$ and by Lemma \ref{lemma about approximate units} for arbitrary  $b\in \TT(\tau,\sigma)$ we have  
$
(ba)^*(ba)=a^*b^*ba \in \KK(\rho,\rho),
$ that is $ba\in \KK(\tau,\rho)$. This shows that $\TT(\tau,\sigma) \KK(\sigma,\rho)\subset  \KK(\tau,\rho)$. From \eqref{ideal defining formula} it follows that the star functor  preserves $\KK$  and thus 
$\KK(\tau,\sigma)\TT(\sigma,\rho)= (\TT(\rho,\sigma)\KK(\sigma,\tau))^*  \subset (\KK(\rho,\tau))^*=\KK(\tau,\rho)$.
\end{proof}
As a first application of the above statement we construct  an annihilator of an ideal in a $C^*$-precategory. 
\begin{proposition}\label{annihilator preposition}
If  $\JJ$ is an ideal in a $C^*$-precategory $\TT$, then there exists a unique ideal  
 $\JJ^\bot$  in $\TT$ such that
 $$
 \JJ^\bot(\rho,\rho)=\JJ(\rho,\rho)^\bot,\qquad \sigma \in \TT.
 $$
Accordingly, we have $\JJ^\bot\cap \JJ=\{0\}$ and for any ideal $\KK$ in $\TT$  
 $$
 \JJ \cap \KK =\{0\} \,\, \Longrightarrow \KK \subset \JJ^\bot
 $$
 where $\{0\}$ denotes the ideal consisting of zero morphisms.
\end{proposition}
 The ideal $\JJ^\bot$ will be called \emph{an annihilator of }$\JJ$. 
\begin{proof} In view of Theorem  \ref{diagonal of ideals} we need to verify  that the ideals $\JJ(\rho,\rho)^\bot$, $\rho\in \TT$, satisfy  condition \eqref{ideal defining condition}. 
Suppose then that  $a\in \TT(\sigma,\rho)$ is  such that $a^*a\in \JJ(\rho,\rho)^\bot$. By Lemma \ref{lemma about approximate units} 
$a^*b^*ba \in \JJ(\rho,\rho)^\bot$, for all $b\in \TT(\sigma,\sigma)$. In particular, for all $b\in \JJ(\sigma,\sigma)$ we have  
$a^*b^*ba\in (\JJ\cap\JJ^\bot)(\rho,\rho)=\{ 0\}$ and since  
$$
 a^*b^*ba= 0 \Longrightarrow ba=  0  \Longrightarrow baa^*= 0,
$$
it follows that   $aa^*\in \JJ(\sigma,\sigma)^\bot$.
\end{proof}

We now turn to discussion of  maps between $C^*$-precategories.
\begin{definition}
A \emph{homomorphism} $\Phi$ from a $C^*$-precategory $\TT$ to a $C^*$-pre\-ca\-te\-gory  $\SSS$ consists of a mapping $\TT \ni\sigma \mapsto \Phi(\sigma)\in \SSS$ and  linear operators 
$ \TT(\sigma,\rho)\ni a \mapsto \Phi(a) \in \SSS(\Phi(\sigma),\Phi(\rho))$,  $\sigma, \rho\in \TT$, such that 
$$
\Phi(a)\Phi(b)=\Phi(ab),\quad \Phi(a^*)=\Phi(a)^*,\quad  b\in \TT(\sigma,\rho),\,\, a\in \TT(\tau,\sigma),\,\, \sigma, \rho,\tau\in \TT.
$$
The notions such as an  \emph{isomorphism}, \emph{endomorphism},  etc., for $C^*$-precategories are defined in an obvious fashion.  
\end{definition}
A homomorphism $\Phi:\TT\to \SSS$ between $C^*$-categories   is a functor iff it is "unital", that is  for every $\rho \in \TT$, $\Phi$ maps the unit in $\TT(\rho,\rho)$ onto the unit in $\SSS(\Phi(\rho),\Phi(\rho))$. In general,  as the operators $\Phi:\TT(\rho,\rho)\to \SSS(\Phi(\rho),\Phi(\rho))$ are $^*$-homomorphisms of $C^*$-algebras, 
using  $C^*$-equality, 
one easily gets
\begin{proposition}\label{proposition 1.5}
For any homomorphism $\Phi:\TT\to \SSS$  of  $C^*$-precategories   the  operators 
 $$
\Phi:\TT(\sigma,\rho)\to \SSS(\Phi(\sigma),\Phi(\rho)), \qquad \sigma,\rho \in \TT
$$  
 are  contractions, and if they are injective then they are isometries. 
 \end{proposition} 
Clearly,  %for our future purposes 
if $\Phi:\TT\to \SSS$ is a homomorphism and $\KK$ is an ideal in $\SSS$   the collection of sets
 $$
 \Phi^{-1}(\KK)(\sigma,\rho):=\{a\in \TT(\sigma,\rho): \Phi(a)\in \KK(\Phi(\sigma),\Phi(\rho))\},
 $$  
forms an ideal in $\TT$, which we shall   refer to as the \emph{preimage of  $\KK$}.  In particular, preimage of  the zero ideal will be  called a  \emph{kernel} of $\Phi$ and denoted by $\ker \Phi$. 
 \begin{proposition}\label{quotient C*-category}
If $\KK$ is an ideal in a $C^*$-precategory $\TT$, then the precategory $\TT/\KK$ whose morphisms are given by the quotient spaces
$$
(\TT/\KK)\,(\sigma,\rho):=\TT(\sigma,\rho)/ \KK(\sigma,\rho)
$$
is a $C^*$-precategory and the quotient maps $q_\JJ:\TT(\sigma,\rho)\to (\TT/\KK)\,(\sigma,\rho)$ give rise to the quotient homomorphism of $C^*$-precategories $q_\JJ:\TT \to \TT/\KK$.  
\end{proposition}
\begin{proof}
Mimic the argument  that shows  the corresponding fact for $C^*$-algebras.
\end{proof}
%The reader easily sees that within the notations of the foregoing  proposition  of $C^*$-categories in the sense to be made precise below.  
 \begin{definition}\label{representations of categories definition}
By a \emph{representation of a $C^*$-precategory} $\TT$ in a $C^*$-algebra $B$ we mean a homomorphism $\pi:\TT \to B$ where $B$ is considered as a $C^*$-precategory with a single object. Equivalently $\pi$ may be treated as a collection $\{\pi_{\sigma,\rho}\}_{\sigma,\rho\in \TT}$ of linear operators  $\pi_{\sigma,\rho}:\TT(\sigma,\rho)\to B$ such that
$$
\pi_{\sigma,\rho}(a)^*=\pi_{\rho,\sigma}(a^*), \quad \textrm{ and } \quad  \pi_{\tau,\rho}(ba)=\pi_{\tau,\sigma}(b)\pi_{\sigma,\rho}(a),
$$ 
for  $a\in \TT(\sigma,\rho)$, $b\in \TT(\tau,\sigma)$.  
By a \emph{representation} of $\TT$ \emph{in a Hilbert space} $H$  we mean a representation of  $\TT$ in the $C^*$-algebra $L(H)$ of all bounded operators. We say that the representation is faithful if all the mappings $\{\pi_{\sigma,\rho}\}_{\sigma,\rho\in \TT}$ are injective.
 \end{definition}  
\begin{remark}\label{gelfand-naimark remark} The above  definition   differs from  \cite[Def. 1.8]{glr}. However, as in \cite[Prop. 1.14]{glr} (which can  be  easily refined to deal with $C^*$-precategories) one can show that   every $C^*$-precategory $\TT$ may be faithfully  represented in a Hilbert space.
 \end{remark}
 The following statement generalizes an elementary, but frequently used, fact that a representation of an ideal in a $C^*$-algebra extends to a representation of this $C^*$-algebra.  
\begin{proposition} \label{extensions of representations on Hilbert spaces0} Suppose that $\KK$ is an ideal in $\TT$ and  $\{\pi_{\sigma,\rho}\}_{\sigma,\rho\in \TT}$ is a representation of $\KK$ in a Hilbert space $H$. If we denote by    $P_\rho$, $\rho\in \TT$, the orthogonal projection onto the essential subspaces of $\pi_{\rho,\rho}$:   
$$
  P_\rho H=\pi_{\rho,\rho}(\KK(\rho,\rho))H,
  $$ 
 then 
    essential subspaces of $\pi_{\sigma,\rho}$, $n\in\N$, are contained in $P_\rho H$, and there is  a unique extension $\overline{\pi}=\{\overline{\pi}_{\sigma,\rho}\}_{\sigma,\rho\in \TT}$ of  $\pi$  to a representation of a $C^*$-precategory $\TT$   such that 
essential subspaces of $\overline{\pi}_{\sigma,\rho}$ is  contained in $P_\rho H$, $\sigma,\rho\in\N$.  This  extension is   determined by relations
\begin{equation}\label{formula defining extensions of right tensor representations}
\overline{\pi}_{\sigma,\rho}(a) \pi_{\rho,\rho}(b)h = \pi_{\sigma,\rho}(ab)h,
\end{equation}
$a\in \TT(\sigma,\rho), \, b \in \KK(\rho,\rho), h \in H$.
\end{proposition}
\begin{proof}
By Lemma \ref{lemma about approximate units} we have $\KK(\rho,\sigma)=\KK(\rho,\rho)\KK(\rho,\sigma)$ and since the  essential subspace of $\pi_{\sigma,\rho}$ is  $\pi_{\sigma,\rho}(\KK(\sigma,\rho))^*H=\pi_{\rho,\sigma} (\KK(\rho,\sigma))H$  it is contained in $P_\rho H$.
\\
 Since for  $a\in \TT(\sigma,\rho), \, b \in \KK(\rho,\rho), h \in H$, and an approximate unit $\{\mu_\lambda\}$ for $\KK(\rho,\rho)$  we have 
\begin{align*}
	\|\pi_{\sigma,\rho}(ab)h\| =& \lim_\lambda \|\pi_{\sigma,\rho}(a\mu_\lambda)\pi_{\rho,\rho}(b)h\|\leq \lim_\lambda \|\pi_{\sigma,\rho}(a\mu_\lambda)\| \|\pi_{\rho,\rho}(b)h\|
	\\
	 \leq & \|a\|\cdot  \|\pi_{\rho,\rho}(b)h\|,
\end{align*}
formula \eqref{formula defining extensions of right tensor representations}  give rise to the representation  $\overline{\pi}_{\sigma,\rho}$ of $\TT(\sigma,\rho)$ on $P_\rho H$.
 Defining $\overline{\pi}_{\sigma,\rho}$ to be zero on  $(P_\rho H)^{\bot}$  one readily sees that $\overline{\pi}=\{\overline{\pi}_{\sigma,\rho}\}_{\sigma,\rho\in \TT}$ is a  representation of  $\TT$ in $H$. Obviously,  every representation $\overline{\pi}=\{\overline{\pi}_{\sigma,\rho}\}_{\sigma,\rho\in \TT}$  that extends  $\pi$  satisfies  \eqref{formula defining extensions of right tensor representations}, which  together with the requirement that  $\overline{\pi}_{\sigma,\rho}|_{_{(P_\rho H)^\bot}}\equiv 0$, for all $\rho,\sigma \in \TT$,   determines $\overline{\pi}$ uniquely. % Checking that $\overline{\pi}$ is  right tensor representation will  be straightforward once we prove iii).
  \end{proof}
 \begin{example}[$C^*$-category of Hilbert modules]\label{category of hilbert modules example1}
Let $\{X_\rho\}_{\rho\in \TT}$ be the family of right Hilbert modules over a  $C^*$-algebra $A$, indexed by a  collection of objects $\TT$. Then  $\TT$   with morphisms being the adjointable maps between the Hilbert $A$-modules: 
$$
\TT(\sigma,\rho):=\LL(X_\rho,X_\sigma), \qquad \rho,\, \sigma \in \TT, 
$$
becomes  a  $C^*$-category. The collection  of spaces of "compact" operators 
$$
\KK(\sigma,\rho):=\KK(X_\rho,X_\sigma),\qquad \rho,\sigma\in \TT.
$$
give rise to the \emph{ideal $\KK$ of "compact operators"} in $\TT$.
Moreover, using Proposition \ref{corollary on ideals in adjointable maps}, we see that every  ideal $J$ in $A$  naturally   give rise to two ideals in $\TT$. One, denoted by  $\JJ$,   consists of all the adjointable maps with ranges in the spaces  $X_\rho J$, $\rho\in \TT$:
$$
\JJ(\sigma,\rho):=\LL_J(X_\rho,X_\sigma),\qquad \rho,\sigma\in \TT.
$$
The other one is  $\KK\cap \JJ$, that is $
(\KK \cap \JJ)(\sigma,\rho)=\KK(X_\rho,X_\sigma J)$, $\rho,\sigma\in \TT.
$
\end{example}
\subsection{$C^*$-precategory $\TT_X$ of a $C^*$-correspondence}
 \label{paragraph for TT(X)1} 
Here we examine our model example of a $C^*$-precategory, which  shall be equipped with an additional structure of a right tensor $C^*$-precategory  in Example \ref{paragraph for TT(X)}. Throughout this subsection we fix a $C^*$-correspondence $X$ over a $C^*$-algebra $A$.
\begin{definition}\label{definition of T_X}
 We define a \emph{$C^*$-precategory $\TT_X$ of the $C^*$-correspondence} $X$ to be the $C^*$-precategory whose objects are natural numbers $\N=\{0,1,...\}$,    spaces of morphisms are the following adjointable maps
 $$
\TT_X(n,m):=\begin{cases} 
\KK( X^{\otimes m},X^{\otimes n}), & \textrm{ if  }\,\, n=0\textrm{ or } m =0, \\
\LL( X^{\otimes m},X^{\otimes n}), & \textrm{ if  }\,\, n,m\geq 1,
\end{cases}
$$
 where  $X^{\otimes n} := X\otimes\dotsm\otimes X$ is the $n$-fold tensor product  and $X^{\otimes 0}:=A$ is the standard Hilbert $A$-module. Moreover, we shall use the  Riesz-Fr\'echet theorem  \cite[Lem. 2.32]{rw}
to  assume the following identifications 
\begin{equation}\label{Riesz-Frechet}
 \KK(A,X^{\otimes n})= X^{\otimes n}, \qquad  \KK(X^{\otimes n},A)=\X^{\otimes n}, \qquad n\in \N,
\end{equation}
where $\X$ denotes  the dual $C^*$-correspondence to $X$, i.e. a left Hilbert $A$-module equipped with a right action such that there is an anti-linear isomorphism $\flat:X\to \X$ that preserves the corresponding structures ($\X^{\otimes n}:=\widetilde{ X^{\otimes n}}$). Thus we have 
$$
\TT_X(0,0)=A,\qquad  \TT_X(n,0)=X^{\otimes n} ,\qquad  \TT_X(0,m)=\X^{\otimes m}.
$$
\end{definition} 
Obviously, we may consider   $\TT_X$   as a sub-$C^*$-precategory  in the $C^*$-category $\TT=\{\LL( X^{\otimes m},X^{\otimes n})\}_{n,m\in \N}$, and if $A$ is unital, then actually   $\TT_X=\TT$. The reason why we deal with $\TT_X$ rather than  $\TT$ is explained in Remark \ref{remark babel}. 
Likewise in Example \ref{category of hilbert modules example1} we associate with an   ideal $J$ in $A$ two ideals in $\TT_X$.
\begin{definition}\label{notation for ideals definition} 
 We denote  by $\KK_X:=\{\KK( X^{\otimes m},X^{\otimes n})\}_{n,m\in \N}$ the ideal in $\TT_X$ consisting of "compact" operators, and for an ideal  $J$ in $A$ we put 
 $$ 
\KK_X(J):=\{\KK( X^{\otimes m},X^{\otimes n}J)\}_{n,m\in \N},
$$
$$
\TT_X(J):=\{\LL_J( X^{\otimes m},X^{\otimes n})\}_{n,m\in \N} \cap \TT_X.
$$
%In particular, $\KK_X(J)$ and $\TT_X(J)$ are ideals in $\TT_X$,   and  $\KK_X(J)=\TT_X(J)\cap\KK_X$. 
\end{definition}
These ideals give useful estimates for arbitrary ideals in $\TT_X$.
\begin{proposition}\label{Ideals in TT(X)}
Let   $\JJ$ be an ideal in  $\TT_X$ and put $J:=\JJ(0,0)$. Then  
$$
\KK_X(J)\subset  \JJ \subset \TT_X(J).
$$ 
In particular, the relations 
\begin{equation}\label{relations that establish bijection between ideals}
 J=\JJ(0,0),\qquad \JJ=\KK_X(J) 
\end{equation}
 establish a one-to-one correspondence between  ideals  in $A$ and   $\KK_X$. 
 \end{proposition}
\begin{proof}
By Theorem   \ref{diagonal of ideals} an element $x\in X^{\otimes n}=\TT_X(n,0)$ belongs to $\JJ(n,0)$ iff $x^*x=\langle x,x\rangle_A\in J$. Thus, by  Hewitt-Cohen Factorization Theorem, we have   $\JJ(n,0)=X^{\otimes n} J$. Since for $x\in \JJ(n,0)=X^{\otimes n} J$ and $y\in \TT_X(m,0)=X^{\otimes m}$ we identify $xy^*\in \JJ(n,m)$ with the "one-dimensional" operator $\Theta_{x,y}\in \KK(X^{\otimes m},X^{\otimes n}J)$, one sees that $\KK(X^{\otimes m},X^{\otimes n} J)\subset \JJ(n,m)$ and consequently $\KK_X(J)\subset  \JJ$. To prove that $\JJ \subset \TT_X(J)$ let $a\in \JJ(n,m) \subset \LL(X^{\otimes m},X^{\otimes n})$. Since for arbitrary  $x\in \TT(n,0)=X^{\otimes n} $ and $y\in \TT(m,0)=X^{\otimes m}$ we have
$$
x^*a^* y=\langle ax,y\rangle_A \in J,
$$
it follows that $a$ takes values in $X^{\otimes n} J$.
\end{proof}
\begin{corollary}\label{Ideals in II_X(I)}
For any ideal $I$ in $A$ relations  \eqref{relations that establish bijection between ideals}
 establish a one-to-one correspondence between  ideals $J$ in $I$ and ideals $\JJ$ in  $\KK_X(I)$.
\end{corollary}
\begin{proof}
Use the transitivity of relation of being an ideal.  
\end{proof}
 %Taking $X=A$ where $A$ is non-unital commutative $C^*$-algebra, one easily sees that the presented estimates in above statement are sharp. 
\begin{remark}\label{remark on difference between  II_{XI} and II_X(I)}
For an ideal $I$ in $A$, $XI$ is naturally considered as a $C^*$-correspon\-dence over $I$.  Applying Proposition \ref{Ideals in TT(X)} to $XI$ we get a one-to-one correspondence between  ideals $J$ in $I$ and ideals $\JJ$ in $\KK_{XI}$, given by relations
 $$
  J=\JJ(0,0),\qquad \JJ=\KK_{XI}(J). 
 $$  
Thus by Corollary \ref{Ideals in II_X(I)} the ideal structures of $\KK_{XI}$ and $\KK_{X}(I)$ are isomorphic, even though   $\KK_{XI}=\{\KK( XI^{\otimes m},XI^{\otimes n})\}_{n,m\in \N}\neq \KK_X(I)=\{\KK( X^{\otimes m},X^{\otimes n}I)\}_{n,m\in \N}$ unless $\phi(I)X=XI$. This phenomena will be pursued in Subsection  \ref{structure theorem for relative Cuntz-Pimsner algebras subsection}.  
 \end{remark}
 We illustrate the introduced objects in the context of $C^*$-correspondences from Examples \ref{partial morphism category1.0}, \ref{precategory of a directed graph ex1}.
 \begin{example}[$C^*$-precategory of a partial morphism]\label{partial morphism category} Let $X_\varphi$ be the  $C^*$-cor\-res\-pon\-dence of a partial morphism $\varphi:A\to M(A_0)$. %, cf . Example \ref{partial morphism category1.0}. 
 We shall call  $\TT_\varphi:=\TT_{X_\varphi}$  a \emph{$C^*$-precategory of} $\varphi$.  We have the following natural identifications for the ideal    
 $\KK_{\varphi}:=\KK_{X_\varphi}=\{\KK( X_\varphi^{\otimes m},X_\varphi^{\otimes n})\}_{n,m\in \N}$:
 \begin{equation}\label{ideal for partial morhpsims}
\KK_{\varphi}(n,m)=\varphi\Big( \varphi\big(...\varphi( \underbrace{A_0) A_0)...\big)\Big)A_0}_{n}\,A\, \varphi\Big( \varphi\big(...\varphi( \underbrace{A_0) A_0)...\big)\Big)A_0}_{m},
\end{equation}
$n,m\in \N$, and  in view of Proposition \ref{Ideals in TT(X)}  every ideal $\JJ$ in $\KK_{\varphi}$ is of the form %$\KK_{\varphi}(J)$, that is
 $$
 \JJ(n,m)=\varphi\Big( \varphi\big(...\varphi( \underbrace{A_0) A_0)...\big)\Big)A_0}_{n}\,J\, \varphi\Big( \varphi\big(...\varphi( \underbrace{A_0) A_0)...\big)\Big)A_0}_{m},\qquad n,m\in \N,
 $$
 where $J=\JJ(0,0)$ is an ideal in $A$. In particular, if $\varphi$ is a partial morphism  arising from a partial automorphism $(\theta,I,A_0)$, then denoting by $D_n$ the domain of $\theta^{-n}$, cf. \cite{exel1}, the above ideals are given by
 $$
 \KK_{\varphi}(n,m)=D_{\max\{m,n\}},\qquad \JJ(n,m)=J\cap D_{\max\{m,n\}} , \qquad n,m\in \N.
 $$ 
If in turn $\varphi$ arise  from an endomorphism $\al:A\to A$ , we then have
 $$
 \KK_{\varphi}(n,m)=\al^n(A)A\al^m(A),\qquad \JJ(n,m)=\al^n(A) J\al^m(A), \qquad n,m\in \N.
 $$
 In this event we shall  denote the $C^*$-precategory $\TT_\varphi$ by $\TT_\al$.
  \end{example}
 
 \begin{example}[$C^*$-precategory of a directed graph]\label{precategory of a directed graph example} 
We set $\TT_E:=\TT_{X_E}$ where $X_{E}$ is a $C^*$-correspondence  associated with a directed graph $E=(E^0,E^1,r,s)$, %, cf. Example \ref{precategory of a directed graph ex1}, 
and we call it a \emph{$C^*$-precategory of} $E$. For $n \geq 1$ we denote by $E^n$ the set of all paths of length $n$, i.e. the set of sequences $(e_1,e_2,...,e_n)$ where $e_i\in E^1$ and $r(e_i)=s(e_{i+1})$, and we put $r^{(n)}(e_1,e_2,...,e_n):=r(e_n)$, $s^{(n)}(e_1,e_2,...,e_n):=s(e_1)$. The quadruple $E^{(n)}:=(E^0,E^{n},r^{(n)},s^{(n)})$ is a directed graph and the $C^*$-correspondence $X_{E^{(n)}}$  may be naturally identified with $X_E^{\otimes n}$. In particular,   $X_E^{\otimes n}$ is spanned by the point masses $\{\delta_\mu: \mu\in E^n\}$ and the ideal of "compact" operators  $\KK_{E}:=\KK_{X_E}$ 
 is spanned by the "matrix units" 
\begin{equation}\label{matrix units of graphs}
\Theta_{\delta_\mu,\delta_\nu}\qquad \textrm{ such that }\qquad  r^{(n)}(\mu)=r^{(m)}(\nu),\qquad \mu \in E^n, \nu\in E^m.
\end{equation}
Since every ideal in $A=C_0(E^0)$ is determined by its hull contained in $E^0$
, in view of Proposition \ref{Ideals in TT(X)},   the equalities 
\begin{equation}\label{ideals of coimoterictiy vs vertices}
\JJ(n,m)=\clsp\{\Theta_{\delta_\mu,\delta_\nu}: r^{(n)}(\mu)=r^{(m)}(\nu)\notin V,\mu \in E^n, \nu\in E^m\},
\end{equation}
 for $n,m\in \N$, 
establish a one-to-one correspondence between subsets $V$ of $E^0$ and ideals $\JJ$ in $\KK_{E}$.
 \end{example}
\section{Right tensor $C^*$-pre\-categories and their representations}\label{section three}
The notion of a (strict) tensor $C^*$-category in the case the underlying semigroup is  $\N$ could be adapted  to $C^*$-precategories as follows.
\begin{definition}\label{right tensor C-category definition}
A    $C^*$-precategory $\TT$ with the set of objects $\N=\{0,1,2,...\}$ together with designated endomorphism $\otimes 1:\TT\to \TT$ sending $n$ to $n+1$: 
$$
\otimes 1:\TT(n,m)\to \TT(n+1,m+1)
$$
shall be called a \emph{right tensor  $C^*$-precategory}.
We  say that $\otimes 1$  is  a \emph{right tensoring} on $\TT$ and  write $a \otimes 1$ instead of $\otimes 1(a)$ for a morphism $a$ in $\TT$. If $\TT$ is a $C^*$-category we shall refer to it simply as a \emph{right tensor  $C^*$-category}.
 \end{definition}
Iterating  a right tensoring $\otimes 1$ on a $C^*$-precategory $\TT$  one gets the semigroup  $\{\otimes 1^k\}_{k\in \N}$ of  endomorphisms  $\otimes 1^k:\TT\to \TT$
$$ \otimes 1^k:\TT(n,m)\to \TT(n+k,m+k),$$
 where by convention we put $\otimes 1^0:=id$. The model example of a right-tensor $C^*$-precategory is
the $C^*$-precategory  $\TT_X$  with  a right tensoring induced by the homomorphism $\phi:A\to \LL(X)$ (cf. \cite{dpz}, \cite[Def. 1.6, 1.7]{katsura}), which we now describe in detail. %This will be our model example of a right tensor $C^$-category.
\begin{example} [Right tensor $C^*$-precategory $\TT_X$ of a $C^*$-correspondence $X$.] \label{paragraph for TT(X)}
 Let $\TT_X$ be as  in Definition \ref{definition of T_X}.   If $n>0, m>0$, we have a  natural  tensoring $$\LL( X^{\otimes m},X^{\otimes n}) \ni a \longmapsto  a\otimes 1 \in \LL( X^{\otimes (m+1)},X^{\otimes (n+1)}),$$  given by \eqref{right tensoring of operators definition}.
For $a\in \TT_X(0,0)=A$  we put 
 $
 a\otimes 1:= \phi(a)\in \LL(X)
 $.  In order to define the right tensoring on the spaces $\TT_X(n,0)$, $\TT_X(0,n)$, for $n>0$,  we use  the assumed  identifications \eqref{Riesz-Frechet} and the mappings $L^{(n)}:X^{\otimes n}\mapsto \LL(X, X^{\otimes n+1})$ and     $D^{(n)}:\KK(X^{\otimes n},A)\mapsto \LL(X^{\otimes n+1}, X)$   determined by  the formulas
$$
[L^{(n)}(x)](y):=x\otimes y,\qquad    x\in X^{\otimes n}, y\in X,
$$
$$ D^{(n)}(a)(y_1\otimes y_2):=\phi(a(y_1))y_2,\qquad     y_1\in X^{\otimes n}, y_2\in X,
 $$
where $a\in \KK(X^{\otimes n},A)$. Alternatively, when $\KK(X^{\otimes n},A)$ is identified with  $\X^{\otimes n}$ we have  $[D^{(n)}(\flat(x))](y_1\otimes y_2)= \phi(\langle x,y_1\rangle_A )y_2$. We put
$$
\label{following notation}
a \otimes 1 :=\begin{cases} L^{(n)}(a), & a\in \TT_X(n,0), \\
D^{(n)}(a), & a\in \TT_X(0,n).
\end{cases}
$$
In this way $\TT_X$ becomes a right tensor $C^*$-precategory with right tensoring "$\otimes 1$".
 We shall call it a  \emph{right tensor $C^*$-precategory of} $X$. 
\end{example}
\begin{remark}\label{remark babel}
The $^*$-homomorphism $\varphi:A\to \LL(X)$ need not extend to a $^*$-homo\-morphism $M(A)\to  \LL(X)$ unless it is nondegenerate (or $A$ is unital). In other words, in general there is no obvious right tensoring on the $C^*$-category  $\TT:=\{\LL( X^{\otimes m},X^{\otimes n})\}_{n,m\in \N}$ and thus the sub-$C^*$-precategory  $\TT_X\subset \TT$ seems to be a more appropriate object to work with. 
\end{remark}
Applying the above construction to Examples \ref{partial morphism category},  \ref{precategory of a directed graph example}
we get respectively a \emph{right tensor $C^*$-precategory $\TT_\varphi$ of a partial morphism} $\varphi$ and a \emph{right tensor $C^*$-precategory $\TT_E$ of a directed graph} $E$. %In particular cases  right tensoring in these $C^*$-precategories is very simple and natural.
\begin{example}[Right tensor $C^*$-precategory of  a partial morphism]
\label{partial morphism category1}
The $C^*$-pre\-ca\-te\-gory 
$\TT_\varphi$ from Example \ref{partial morphism category} is a right tensor $C^*$-precategory with right tensoring induced by $\varphi$. In particular, if $\varphi$ arises from an endomorphism   $\al:A\to A$ of   a  $C^*$-algebra $A$,   the ideal $\KK_{\varphi}=\{\al^n(A)A\al^m(A)\}_{n,m\in \N}$  is a right tensor $C^*$-precategory itself and the right tensoring on  $\KK_{\varphi}$ assumes the form
$$
\al^n(A)A\al^m(A)\ni a \longrightarrow a \otimes 1 =  \al(a) \in \al^{n+1}(A)A\al^{m+1}(A).
$$
If additionally  $A$ is unital the above formula describes a right tensoring on $\TT_\al=\TT_\varphi$ as then we have $\TT_\al=\KK_{\varphi}=\{\al^n(1)A\al^m(1)\}_{n,m\in \N}$.\end{example}
\begin{example}[Right tensor $C^*$-precategory of a directed graph]
\label{there is no label}
The $C^*$-precate\-gory $\TT_E$ from Example \ref{precategory of a directed graph example} is a right tensor $C^*$-precategory where the right tensoring is induced by "composition" of graphs (equivalently by multiplication of incidence matrices). In the event  the set of edges $E^1$ is finite,  the $C^*$-precategory $\TT_E$ coincides with the ideal of "compact" operators $\KK_{E}$ and  the right tensoring is determined by the formula
$$
\Theta_{\delta_\mu,\delta_\nu} \otimes 1 := \sum_{\{e\in E^1:\,r^{(n)}(\mu)=s(e)\}}\Theta_{\delta_{\mu e},\delta_{\nu e}}\,,\qquad \mu \in E^n, \nu\in E^m,
$$
where $\mu e$ and  $\nu e$ are the paths  obtained by concatenation. 
\end{example}
We now turn to investigation of representations of ideals in right tensor $C^*$-precategories which respect the right tensoring. As we shall see in Proposition \ref{proposition:rho} these representations  generalize  representations of $C^*$-correspondences.
\begin{definition}\label{tensor representation definition} 
Let $\KK$ be an ideal in a right tensor $C^*$-precategory $\TT$. We will say that a representation  $\{\pi_{n,m}\}_{n,m\in \N}$  of $\KK$ %, cf. Definition  \ref{representations of categories definition}, 
is a \emph{right tensor representation} if it satisfies     
\begin{equation}\label{right tensor representation condition}
\pi_{n,m}(a) \pi_{m+k,l}(b)= \pi_{n+k,l}((a\otimes 1^{k})\, b) 
\end{equation}
 for all $a\in \KK(n,m)$ and $b\in \KK(m+k,l)$, $k,l,m,n\in \N$. 
 \end{definition} 
 \begin{remark}
 Since $\KK$ is an ideal  the right hand side  of \eqref{right tensor representation condition}  makes sense. Furthermore,  by  taking adjoints one gets the symmetrized  version of this equation:
  $$
 \pi_{n,m+k}(a)\pi_{m,l}(b)= \pi_{n,l+k}( a (b\otimes 1^{k})),
$$
where $a\in \KK(n,m+k)$ and $b\in \KK(m,l)$, $k,l,m,n\in \N$. 
\end{remark}
We introduce the class of  ideals in a right tensor $C^*$-precategory $\TT$ whose right tensor representations extends naturally to $\TT$. 
\begin{definition}
We shall say that a right tensoring $\otimes 1$ \emph{acts nondegenerately} on an ideal $\KK$ in a right tensor $C^*$-precategory $\TT$, or that $\KK$ is \emph{$\otimes 1$-nondegenerate} if 
$$
(\KK(n,m)\otimes 1)\KK(n+1,m+1)= \KK(n+1,m+1),\qquad n,m\in \N,
$$
briefly $(\KK\otimes 1)\KK=\{\KK(n,m)\}_{n,m>0}$.
 \end{definition}
 We use the term "nondegeneracy" because of the similarity to a condition defining nondegeneracy of homomorphism of $C^*$-algebras, cf Example \ref{partial morphism category1.0}.  By Lemma \ref{lemma embedding of tensors} ii) the ideal $\KK_X:=\{\KK( X^{\otimes m},X^{\otimes n})\}_{n,m\in \N}$ in $\TT_X$ is $\otimes 1$-nondegenerate, and actually  any $\otimes 1$-non\-degene\-rate ideal $\KK$ in a right tensor $C^*$-precategory $\TT$  may be naturally embedded  as a sub-$C^*$-precategory into $\TT_X$ associated with the $C^*$-correspondence given by relations
 $$
A=:\KK(0,0), \qquad X:=\KK(1,0),\qquad  \langle x,y\rangle_A:=x^*y, \qquad \phi(a)=a\otimes 1,
 $$
where $x,y\in X$ and $a\in A$, and then $\KK_X \subset \KK\subset \TT_X$ (we leave the details to the reader, as we shall not use this fact in the sequel).
  The following statement is a version of Proposition \ref{extensions of representations on Hilbert spaces0} adapted to right tensor $C^*$-precategories.
\begin{proposition} \label{extensions of representations on Hilbert spaces} Suppose that $\KK$ is an ideal in $\TT$ and  $\pi=\{\pi_{n,m}\}_{n,m\in \N}$ is a right tensor representation of $\KK$ in a Hilbert space $H$.  Then  the  projections  $P_m$, $m\in \N$, onto the essential subspaces of $\pi_{m,m}$ are pairwisely commuting 
and if   $\KK$ is  $\otimes 1$-nondegenerate, then we actually have 
\begin{equation}\label{decreasing sequence}
P_1 \geq P_2 \geq P_3 \geq ...
\end{equation}
Moreover, if \eqref{decreasing sequence} holds, then the extension $\overline{\pi}=\{\overline{\pi}_{n,m}\}_{n,m\in \N}$ of $\pi=\{\pi_{n,m}\}_{n,m\in \N}$ described in Proposition \ref{extensions of representations on Hilbert spaces0} is a right tensor representation of $\TT$, and    
\begin{equation}\label{right tensor defining relations}
\overline{\pi}_{n,m}(a)P_{m+k} =\overline{\pi}_{n+k,m+k}(a\otimes 1^k), \qquad P_{n+k}\overline{\pi}_{n,m}(a) =\overline{\pi}_{n+k,m+k}(a\otimes 1^k),
\end{equation}
for   $a\in \TT(n,m)$, $k\in \N$.
\end{proposition}
\begin{proof}
Let us note that the equality $\pi_{m,m}(a)\pi_{m+k,m+k}(b)=\pi_{m+k,m+k}((a\otimes 1^k) b)$ implies $P_{m+k}\pi_{m,m}(a)P_{m+k}=\pi_{m,m}(a)P_{m+k}$ and this together with a similar equality for $a^*$ yields $P_{m+k}\pi_{m,m}(a)=\pi_{m,m}(a)P_{m+k}$. Therefore $P_{m+k}P_{m}=P_{m}P_{m+k}$.
\\
If $\KK$ is $\otimes 1$-nondegenerate, then $(\KK(m,m)\otimes 1)\KK(m+1,m+1)= \KK(m+1,m+1)$ together 
with \eqref{right tensor representation condition} imply that $P_m P_{m+1}=P_{m+1}$, that is $P_m \geq P_{m+1}$. In general, if \eqref{decreasing sequence} holds, then any element $h$ in $P_{m+k}H$ may be written in the form $\pi_{m,m}(b)\pi_{m+k,m+k}(c)h_0$, $b\in \KK(m,m)$, $c\in \KK(m+k,m+k)$ and then for $a\in \TT(n,m)$ we have
\begin{align*}
\overline{\pi}_{n,m}(a)h&=\overline{\pi}_{n,m}(a)\pi_{m,m}(b)\pi_{m+k,m+k}(c)h_0=\pi_{n,m}(ab)\pi_{m+k,m+k}(c)h_0
\\
&=\pi_{n+k,m+k}((ab)\otimes 1^kc)h_0=\pi_{n+k,m+k}((a\otimes 1^k) (b\otimes 1^k)c)h_0
\\
&=\overline{\pi}_{n+k,m+k}(a\otimes 1^k)\pi_{m+k,m+k}( (b\otimes 1^k)c)h_0\\
&=\overline{\pi}_{n+k,m+k}(a\otimes 1^k)\pi_{m,m}(b)\pi_{m+k,m+k}( c)h_0
=\overline{\pi}_{n+k,m+k}(a\otimes 1^k)h.
\end{align*}
Hence $\overline{\pi}_{n,m}(a)P_{m+k} =\overline{\pi}_{n+k,m+k}(a\otimes 1^k)$ and by passing to adjoints one gets $P_{n+k}\overline{\pi}_{n,m}(a) =\overline{\pi}_{n+k,m+k}(a\otimes 1^k)$. In particular, for $a\in \TT(n,m)$ and $b\in \TT(m+k,l)$ we have
\begin{align*}
\overline{\pi}_{n,m}(a)\overline{\pi}_{m+k,l}(b)
&=\overline{\pi}_{n,m}(a)\, P_{m+k}\,\overline{\pi}_{m+k,l}(b)=\overline{\pi}_{n+k,m+k}(a\otimes 1^k)\overline{\pi}_{m+k,l}(b) \\
&=\overline{\pi}_{n+k,l}((a\otimes 1^k)  b),
\end{align*}
that is $\overline{\pi}=\{\overline{\pi}_{n,m}\}_{n,m\in \N}$  is a right tensor representation.
  \end{proof}
\begin{corollary}\label{representations extension core}  
If $\KK$ is an ideal in $\TT$ and $\KK$ is $\otimes 1$-nondegenerate, then we have a one-to-one correspondence between   right tensor representations   of  $\KK$  in the Hilbert space $H$ and   right tensor representations $\{\overline{\pi}_{n,m}\}_{n,m\in \N}$ of $\TT$ in $H$ satisfying
\begin{equation}\label{equlity of essential subspaces}
\overline{\pi}_{m,m}(\TT(m,m))H= \overline{\pi}_{m,m}(\KK(m,m))H,\qquad m\in \N.
\end{equation}
\end{corollary}
\begin{proof} 
Clear by Proposition \ref{extensions of representations on Hilbert spaces} and Lemma \ref{lemma about approximate units}.
\end{proof}
\subsection{Right tensor representations and representations of $C^*$-cor\-res\-pon\-den\-ces}
Let $X$ be a $C^*$-correspondence over $A$. We shall investigate  relationships between   right tensor representations   of the  right tensor $C^*$-precategory $\TT_X$, its ideal 
$
\KK_X=\{\KK( X^{\otimes m},X^{\otimes n})\}_{n,m\in \N}
$, and  representations  of    $X$. 
 \begin{definition} A {\em representation} $(\pi,t)$ of the $C^*$-correspondence 
$X$  in a
$C^*$-algebra
$B$ consists of a linear map $t:X\to B$
and a $^*$-homomorphism $\pi:A\to B$ such that
$$
t(x\cdot a) = t(x)\pi(a),\quad t(x)^*t(y)= \pi(\langle x,y\rangle_A),\quad t(a\cdot x) = \pi(a)t(x),
$$
for $x,y\in X$ and $a\in A$. If $\pi$ is faithful (then automatically $t$ is isometric, cf. \cite{fr}, \cite{mt}) we  say that the  representation  $(\pi,t)$  is \emph{faithful}. If  $B=L(H)$ for  a Hilbert space $H$ we say that   $(\pi,t)$ is a \emph{representation of $X$ in} $H$. 
\end{definition}
\begin{remark}
The above introduced notion  is called a \emph{Toeplitz representation} of $X$ in \cite{fr}, \cite{fmr}, and  an \emph{isometric covariant representation} of $X$ in \cite{ms}.
\end{remark}
The first step is to show that representation $(\pi,t)$ give rise to a right tensor representation of $\KK_X$ and this, in essence,  follows from the results of   \cite{kpw},  \cite{fmr}, \cite{fr}, \cite{p} where it was used in an implicit form. 
\begin{proposition}\label{proposition:rho}
If  $(\pi,t)$ is a representation of   $X$ in a $C^*$-algebra $B$, then there is a unique right tensor representation $\{\pi_{n,m}\}_{m,n\in \N}$ of the ideal $\KK_X$ in the right tensor $C^*$-precategory $\TT_X$, such that 
\begin{equation}\label{takie tam odpowieniosci miedzy reprezentacajmi}
\pi_{0,0}=\pi,\qquad \pi_{1,0}=t.
\end{equation}
We shall denote this representation by $[\pi,t]$. Every right tensor representation of $\KK_X$ is of the form $[\pi,t]$ and 
$$
\ker[\pi,t]=\KK_X(\ker\pi),
$$ 
cf. Definition \ref{notation for ideals definition}.
\end{proposition}
\begin{proof}
Suppose that $\{\pi_{n,m}\}_{m,n\in \N}$ is a right tensor representation of $\KK_X$ such that \eqref{takie tam odpowieniosci miedzy reprezentacajmi} holds. Then  in view of Definition  \ref{tensor representation definition},  for $x_i\in X$, $i=1,...,n$,
we have %$\pi_{n0}( x_1\otimes\dotsm\otimes x_n)=\pi_{10}( x_1)...\, \pi_{10}(x_n)$, that is 
 \begin{equation}\label{definition of pi_n0}
\pi_{n,0}( x_1\otimes\dotsm\otimes x_n)= t(x_1)t(x_2)\dotsm t(x_n).
\end{equation}
Hence $\pi_{n,0}$, $n>0$, is uniquely determined by $t$. 
Furthermore, $
\pi_{n,m}$ is uniquely determined by  $\pi_{n,0}$ and $\pi_{m,0}$, since  for $x\in X^{\otimes n}$ and $y\in X^{\otimes m}$ we have
\begin{equation}\label{definition of pi_nm}
\pi_{n,m}(\Theta_{x,y})=\pi_{n,0}(x)(\pi_{m,0}(y))^*.
\end{equation}
This proves the uniqueness of the representation $[\pi,t]=\{\pi_{n,m}\}_{n,m\in \N}$. For the existence of $[\pi,t]$ note that formula \eqref{definition of pi_n0} give rise to the linear map $\pi_{n,0}: X^{\otimes n}\to B$ such that $(\pi_{n,0},\pi)$ is a  representation of $X^{\otimes n}$, cf. \cite[Prop. 1.8]{fr} or \cite[Lem. 3.6]{fmr}. Arguing as in    \cite[Lem. 2.2]{kpw}   or \cite[Lem. 3.2]{p}, one sees that formula \eqref{definition of pi_nm} defines a contraction   $
\pi_{n,m}:\KK( X^{\otimes m},X^{\otimes n})\to B $ such that the family  $\{\pi_{n,m}\}_{n,m\in \N}$ forms a representation of $\KK_X$ in $B$. To check that $\{\pi_{n,m}\}_{n,m\in \N}$ is a right tensor representation let $x={\otimes}_{i=1}^n x_i\in X^{\otimes n}$, $y={\otimes}_{i=1}^m y_i\in X^{\otimes m}$,  $x'={\otimes}_{i=1}^{m+k} x_i'\in X^{\otimes m+k}$ and $y'={\otimes}_{i=1}^l y_i'\in X^{\otimes l}$  (we adhere to the convention that the indexes indicate  the order of factors). One readily sees that
$$
(\Theta_{x,y}\otimes 1^k )\Theta_{x',y'}=\Theta_{z ,y'}
$$ 
where $z= x\langle\otimes_{i=1}^m y_i,  \otimes_{i=1}^{m} x_i'\rangle \otimes x_{m+1}' \otimes ... \otimes x_{m+k}'\in X^{\otimes m+k}$. On the other hand 
\begin{align*}
	\pi_{n,m}(\Theta_{x,y}) \pi_{m+k,l} (\Theta_{x',y'}) &= \prod_{i=1}^n t(x_i) \Big( \prod_{i=1}^m t(y_i)\Big)^* \prod_{i=1}^{m+k} t(x_i') \Big( \prod_{i=1}^l t(y_i')\Big)^*
	\\
	&=\prod_{i=1}^n t(x_i) \pi(\Big\langle\otimes_{i=1}^m y_i, \otimes_{i=1}^{m} x_i' \Big\rangle_A) \prod_{i=1}^{k} t(x_{m+i}') \Big( \prod_{i=1}^l t(y_i')\Big)^*
	\\
	&=\pi_{n+k,m+k}(\Theta_{z,y'}),
\end{align*}
which means that condition \eqref{right tensor representation condition} is satisfied by "rank one" operators and thereby  by all morphisms of the ideal $\KK_X$. 
\\
Clearly, if $\{\pi_{n,m}\}_{n,m\in\N}$ is a right tensor representation  of $\KK_X$, then the pair $(\pi,t)$ where $\pi_{0,0}:=\pi$, $t:=\pi_{1,0}$,  is a representation of $X$ and we have $[\pi,t]=\{\pi_{n,m}\}_{n,m\in\N}$.  To investigate the form of $\ker\pi_{n,m}$  note that for $x\in X^{\otimes n}$ we have
$$
\|\pi_{n,0}(x)\|^2 =\|\pi_{n,0}(x)^*\pi_{n,0}(x)\| =\|\pi(\langle x,x\rangle_A)\|
\le \|\langle x,x\rangle_A\| = \| x\|^2.
$$
It follows that  $\pi_{n,0}(x)=0 \Longleftrightarrow \langle x,x\rangle_A \in \ker \pi
$. Therefore,   by Hewitt-Cohen Factorization Theorem, $
\ker \pi_{n,0} =X^{\otimes n}\ker\pi
$. In view of formulas \eqref{definition of pi_n0}, \eqref{definition of pi_nm}  $a\in\KK(X^{\otimes m}, X^{\otimes n} )$ belongs to $\ker\pi_{n,m}$  iff the range of $a$ is contained in  $
\ker \pi_{n,0} =X^{\otimes n}\ker\pi
$. 
\end{proof}
\begin{remark}\label{remark to define pi_11}
The mapping $\pi_{1,1}$  in the above theorem %(see \eqref{definition of pi_nm})
 is determined by the formula 
$$
\pi_{1,1}(\Theta_{x,y})=t(x)t(y)^*,\qquad x,y \in X.$$ 
In \cite{ms},  \cite{fmr},\cite{fr},  \cite{mt},  it  is denoted by $\pi^{(1)}$ where  it is used to introduce the notion of coisometricity, see Definition \ref{definition for ideals in correspondences} below. 
\end{remark}
If   $(\pi,t)$ is a  representation of $X$ on a Hilbert space $H$, then within the notation of Proposition \ref{proposition:rho}, the essential subspace of $\pi_{m,m}$, $m=1,2,...$, is  
\begin{equation}\label{remarks concerning projections}
\pi_{m,m}(\KK(X^{\otimes m}))H= \clsp\{ t(x_1)\cdot...\cdot  t(x_m)h:x_1,...,x_m \in X, \,\, h\in H\}.
\end{equation}
Hence the projections $P_m$ onto these spaces %,  associated with  $[\pi,t]$ via  Proposition \ref{extensions of representations on Hilbert spaces0}, 
form a decreasing sequence, cf. \cite[Prop. 1.6]{fr}, \cite[4.3]{fmr}. This could be also derived from  Proposition \ref{extensions of representations on Hilbert spaces}, since $\KK_X$ is $\otimes 1$-nondegenerate ideal in $\TT_X$. As a  consequence of Propositions \ref{proposition:rho},  \ref{extensions of representations on Hilbert spaces} we get
\begin{proposition}\label{representation of TT(X) proposition}
Every  representation $(\pi,t)$  of a $C^*$-correspondence   $X$ in a Hilbert space $H$ give rise to a unique 
 right tensor representation $\{\overline{\pi}_{n,m}\}_{n,m\in \N}$ of  $\TT_X$ 
 such that 
 \begin{equation}\label{takie tam odpowieniosci miedzy reprezentacajmi2}
\overline{\pi}_{0,0}=\pi,\qquad \overline{\pi}_{1,0}=t.
\end{equation}
and
   \begin{equation}\label{equlity of essential subspaces2}
\overline{\pi}_{m,m}(\LL(X^{\otimes m}))H= \overline{\pi}_{m,m}(\KK(X^{\otimes m} ))H,\qquad m=1,2,...\,.
\end{equation}
We shall denote this representation by  $\overline{[\pi,t]}$. An arbitrary  right tensor representation   of $\TT_X$ in a Hilbert space is of the form $\overline{[\pi,t]}$ if and only if it satisfies  \eqref{equlity of essential subspaces2}. Moreover, the kernel of $\overline{[\pi,t]}$ is determined by the kernel of $\pi$:
$$
\ker\overline{[\pi,t]}=\TT_X(\ker\pi),
$$
see Definition \ref{notation for ideals definition}.
\end{proposition}
\begin{proof}
The only thing requiring a comment is the form of $\ker\overline{[\pi,t]}$ but this follows from  \eqref{formula defining extensions of right tensor representations} and the equality $
\ker[\pi,t]=\KK_X(\ker\pi).
$
\end{proof}
%Applying this result to $C^*$-correspondences associated to $C^*$-dynamical systems and  graphs we have the following examples.
\begin{example}\label{endomorphism of unitals}  Let $\al:A\to A $ be an endomorphism of a unital $C^*$-algebra $A$. A triple $(\pi,U,B)$ consisting of a unital $^*$-homomorphism $\pi:A\to B$ and a partial isometry $U\in B$ such that 
$$
\pi(\al(a))=U\pi(a) U^*,\quad\textrm{ for all } a\in A, \quad\textrm{ and }\quad  U^*U\in \pi(A)',
$$ 
is called a \emph{covariant representation of $\al$ in $B$}, see \cite{Lin-Rae}, \cite{kwa-leb1}. It is known \cite{fmr}, \cite{kwa-leb1} that covariant representations of $\al$ are in one-to-one correspondence with representations of the  $C^*$-correspondence $X_\al=\al(1)A$. Thus in view of Proposition \ref{proposition:rho} every right tensor representation of the right tensor $C^*$-category $\TT_\al=\{\al^{n}(1)A\al^m(1)\}_{n,m\in \N}$, cf. Examples \ref{partial morphism category}, \ref{partial morphism category1}, assumes the form  %$\{\pi_{nm}\}_{n,m\in \N}$:  
$$
\pi_{n,m}(a):= U^{*n}\pi(a)U^m,\qquad \quad a\in \al^n(1)\A\al^m(1),\,\, n,m \in \N,
$$
for a covariant representation $(\pi,U,B)$ of $\al$ where  $
\pi=\pi_{0,0}$ and $U=\pi_{0,1}(\al(1))
$. In particular, the $C^*$-algebra generated by the image of the $C^*$-category $\TT_\al$ under  $\{\pi_{n,m}\}_{n,m\in \N}$ is the $C^*$-algebra
$$
C^*(\pi(A),U)=\clsp\{U^{*n}aU^m: a \in A, m,n\in \N\}
 $$ 
 generated by $\pi(A)$ and $U$.\\
In a more general context we may define a \emph{covariant representation of a partial morphism} $\varphi:A\to M(A_0)$  to be  a triple $(\pi,U,H)$ consisting of a Hilbert space $H$, a nondegenerate representation $\pi:A\to L(H)$ and a partial isometry $U\in L(H)$ such that 
$$
U\pi(a) U^*=\overline{\pi}(\varphi(a)),\quad\textrm{ for all } a\in A, \quad\textrm{ and }\quad  U^*U\in \pi(A)',
$$ 
where $\overline{\pi}:M(A_0)\to H$ is a representation given by the conditions 
$$
\overline{\pi}(b)|_{(\pi(A_0)H)^\bot}\equiv 0\quad \textrm{and}\quad  \overline{\pi}(b)\pi(a_0)h=\pi(ba_0)h, \quad b\in M(A_0),\,a_0\in A_0, \,h\in H. 
$$
Then the final subspace of $U$ is $\pi(A_0)H$. One  checks that for the triple  $(\pi,U,H)$ defined above the pair $(\pi,t)$ where
$$
t(x):=U^*\pi(x), \qquad x \in X_\varphi=A_0A,
$$
is a representation of the $C^*$-correspondence $X_\varphi$. Conversely, if $(\pi,t)$ is a representation of $X_\varphi$ in a Hilbert space $H$, then for $a_0\in A_0$ and $h\in H$ we have
$$
\|t(a_0)h\|^2=\langle t(a_0)h, t(a_0)h\rangle=\langle \pi(a_0^*a_0)h, h\rangle=\langle \pi(a_0)h, \pi(a_0)h\rangle=\|\pi(a_0)h\|^2.
$$
Thus relations 
$$
Ut(a_0)h=\pi(a_0)h,\qquad  a_0\in A_0,\,\, h\in H, \qquad \textrm{and }\,\, U|_{(t(A_0)H)^{\bot}}\equiv 0
$$ 
define a partial isometry $U\in L(H)$. For $a_0\in A_0,\, a\in A$ and $h\in H$ we have
\begin{align*}
 U\pi(a)U^*\pi(a_0)h&=U\pi(a)t(a_0)h=Ut(\varphi(a)a_0)h=\pi(\varphi(a)a_0)h
 \\
 &=\overline{\pi}(\varphi(a))\pi(a_0)h,
\end{align*}
that is $U\pi(a) U^*=\overline{\pi}(\varphi(a))$. Furthermore 
\begin{align*}
 \pi(a)U^*Ut(a_0)h&=\pi(a)t(a_0)h=t(\varphi(a)a_0)h=U^*Ut(\varphi(a)a_0)h
 \\
 &=U^*U\pi(a)t(a_0)h. 
\end{align*}
Hence $U^*U\in \pi(A)'$ and the triple  $(\pi, U,H)$ is a covariant representation of $\varphi$. As a consequence,  by Proposition \ref{representation of TT(X) proposition}, it follows  that right tensor  representations $\{\overline{\pi}_{n,m}\}_{n,m\in \N}$ of  $\TT_\varphi$ satisfying  \eqref{equlity of essential subspaces2} are in one-to-one correspondence with covariant representations $(\pi,U,H)$  of $\varphi$.  In particular, the  $C^*$-algebra generated by the image of the $C^*$-precategory $\TT_\varphi$ under  $\{\overline{\pi}_{n,m}\}_{n,m\in \N}$ is  contained the $W^*$-algebra
 $W^*(\pi(A),U)$ generated by $\pi(A)$ and $U$.
\end{example}
\begin{example}\label{E-family example} Let $\TT_E$ be the right tensor $C^*$-precategory  %, cf. Example \ref{precategory of a directed graph example}, 
of a directed graph  $E=(E^0,E^1,r,s)$. By a \emph{Toeplitz-Cuntz-Krieger  $E$-family} it is meant  a collection of mutually orthogonal projections  $\{p_v:v\in E^0\}$ together with a collection of partial isometries with orthogonal ranges $\{s_e:e\in E^1\}$ that satisfy relations
$$
s_e^*s_e=p_{r(e)},\qquad s_{e}s_{e}^*\leq p_{s(e)},\qquad e \in E^1.
$$
Such families are in one-to-one correspondence with representations of $X_E$,  cf. \cite{katsura1}, \cite{mt}, \cite{bprs}. Hence, in view of Proposition \ref{proposition:rho}, every  right tensor representation $\{\pi_{n,m}\}_{n,m\in \N}$ of the ideal   $\KK_{E}=\{\KK( X_E^{\otimes m},X_E^{\otimes n})\}_{n,m\in \N}$ in  $\TT_E$, see Example \ref{precategory of a directed graph example}, is given by the formula 
$$
\pi_{n,m}(\Theta_{\delta_\mu,\delta_\nu})=s_{e_1}...s_{e_m}s_{f_n}^*...s_{f_1}^*$$
where $\mu=(e_1,...,e_m) \in E^m$, $\nu=(f_1,...,f_n)\in E^n 
$,
for a Toeplitz-Cuntz-Krieger $E$-family $\{p_v:v\in E^0\}$, $\{s_e:e\in E^1\}$, where 
$$
p_v=\pi_{0,0}(\delta_v),\quad v\in E^0 \qquad \textrm{ and } \qquad s_e=\pi_{1,0}(\delta_{e}),\quad e\in E^1.
$$
In particular, the algebra generated by the image of the ideal $\KK_{E}$ under representation  $\{\pi_{n,m}\}_{n,m\in \N}$ is the $C^*$-algebra 
$$
C^*(\{p_v:v\in E^0\},\{s_e:e\in E^1\})
$$
 generated by the corresponding Toeplitz-Cuntz-Krieger $E$-family.
\end{example}

\subsection{Ideals of coisometricity for right tensor representations}
We  transfer the concept of coisometricity for  representations of   $C^*$-corollary\-responden\-ces,  introduced by  P. Muhly and B. Solel in \cite{ms}, onto the ground of ideals in right tensor  $C^*$-precategories. %To this end, for
\begin{definition}\label{nie wiem jak to nazwac definition} We put 
 $J(\KK):=(\otimes 1)^{-1}(\KK) \cap \KK$ for any ideal $\KK$  in a right tensor $C^*$-precategory $\TT$ . In other words,   $J(\KK)=\{J(\KK)(n,m)\}_{n,m\in\N}$ is an ideal in $\KK$ where
\begin{equation}\label{nie wiem jak to nazwac}
J(\KK)(n,m):=\big\{a\in \KK(n,m): a\otimes 1\in \KK(n+1,m+1)\big\}.\end{equation}
\end{definition}
The ideal $J(\KK)$ plays the role of the ideal $J(X)=\phi^{-1}(\KK(X))$ introduced by Pimsner  in \cite{p}. Roughly speaking, it consists of the elements  for which  the notion of coisometricity  makes sense.
\begin{proposition}\label{proposition defining  ideals associated with representations}
Let $\KK$ be  an ideal in a  right tensor $C^*$-precategory $\TT$. If $\pi=\{\pi_{n,m}\}_{n,m\in \N}$  is a right tensor representation of $\KK$ in a $C^*$-algebra $B$, then  the spaces 
$$
\JJ(n,m):=\Big\{a\in J(\KK)(n,m):  \pi_{n,m}(a)=\pi_{n+1,m+1}(a\otimes 1)\Big\}
$$   
form an ideal  in $J(\KK)$. Moreover, we have 
$
\JJ\cap (\ker \otimes 1)\subset \ker\pi\textrm{ and }\ker\pi\cap J(\KK)\subset \JJ.
$ In particular,  faithfulness of $\pi$ implies that 
$
\JJ \subset (\ker\otimes 1)^\bot.$
\end{proposition}
\begin{proof}
The first part of proposition  follows directly from the definition of right tensor representation.  To show the second part, note that for $a\in \JJ (n,m)\cap \ker(\otimes 1)(n,m)$ we have $\pi_{n,m}(a)=\pi_{n+1,m+1}(a\otimes 1)=\pi_{n+1,m+1}(0)=0$, that is $a\in  (\ker \pi) (n,m)$. If,  in turn, $a\in \ker\pi (n,m)\cap J(\KK)(n,m)$, then  $$\pi_{n+1,m+1}(a\otimes 1)(\pi_{n+1,m+1}(a\otimes 1))^*=\pi_{n,m}(a)(\pi_{n+1,m+1}(a\otimes 1))^*=0,
$$
that is $\pi_{n+1,m+1}(a\otimes 1)=0$ and therefore $a\in \JJ(n,m)$.
\end{proof}

\begin{definition}\label{definition for ideals in categories}
Let $\KK$ be an ideal in  a right tensor $C^*$-precategory $\TT$ and let $\pi=\{\pi_{n,m}\}_{n,m\in \N}$   be  a right tensor representation of $\KK$. We shall say that $\pi$ is \emph{coisometric} on an ideal $\JJ\subset J(\KK)$ if 
$$
\pi_{n,m}(a)=\pi_{n+1,m+1}(a\otimes 1), \qquad \textrm{for all }\,\,\, a\in \JJ(n,m),\, n,m\in \N. 
$$
The ideal $\JJ$ defined in Proposition \ref{proposition defining  ideals associated with representations} is the biggest ideal on which $\pi$ is coisometric and we shall  call it the \emph{ideal of coisometricity for} $\pi$.
\end{definition}
We devote the rest of this subsection to  discuss and reveal  the relationship between the above definition and the following one.

\begin{definition}[cf. \cite{fmr}, \cite{ms}]\label{definition for ideals in correspondences}
%The role of the ideal \eqref{nie wiem jak to nazwac} in the context of $C^*$-correspondences is played by the ideal $J(X):={\phi}^{-1}(\KK(X))$. To clarify this relationship let us  note that  $$  J(\KK_X)=\KK_X(J(X)),\qquad  J(X)=J(\KK_X)(0,0) $$  and recall the following
A  representation $(\pi,t)$ of a $C^*$-correspondence $X$ is called {\em coisometric} on an ideal $J$ contained in $J(X)=\phi^{-1}(\KK(X))$ if
\begin{equation*}
\pi(a)=\pi_{1,1}(a\otimes 1),\qquad\text{ for all \ } a\in J,
\end{equation*}
where $\pi_{1,1}$ is defined in Remark \ref{remark to define pi_11}. The set $\{a\in J(X):\pi(a)=\pi_{1,1}(a\otimes 1)\}$ is  the biggest ideal on which $(\pi,t)$ is coisometric and we shall call it  an \emph{ideal of coisometricity for} $(\pi,t)$.
\end{definition}
As an immediate  corollary of Propositions \ref{Ideals in TT(X)}, \ref{proposition defining  ideals associated with representations}  we get  
\begin{proposition}\label{coisometricty coincidence}
If  $(\pi,t)$ is  a  representation of $X$ in a $C^*$-algebra $B$ and  
\begin{equation}\label{the biggest ideal of coisometricity}
J=\{a\in J(X): \pi_{1 ,1}(\phi(a))=\pi(a)\},
\end{equation}
then the ideal  of coisometricity for the right tensor representation $[\pi,t]$ of the ideal $\KK_X$, cf. Proposition \ref{proposition:rho},  
is $\KK_X(J)=\{\KK(X^{\otimes m}, X^{\otimes n}J)\}_{n,m\in \N}$. 
\end{proposition}
\begin{corollary}\label{coisometricicity of representations}
Relations \eqref{takie tam odpowieniosci miedzy reprezentacajmi} together with equality $J=\JJ(0,0)$ establish a one-to-one correspondence between representations of a $C^*$-correspondence $X$ coisometric on $J\subset J(X)$ and right tensor representations of   $\KK_X$ coisometric on  $\JJ\subset J(\KK_X)$.
\end{corollary}
\begin{corollary}[Prop. 2.21 \cite{ms}]\label{Proposition by muhly and solel}
If $J\subset J(X)$ is an ideal of coisometricity for a representation $(\pi,t)$ of $X$, then $J\cap \ker\phi\subset \ker \pi$. Hence if $(\pi,t)$ is faithful, then
$
J \subset (\ker\phi)^\bot.
$
\end{corollary}
With the help  of  the following lemma we  get a version of Proposition \ref{coisometricty coincidence} for the extended representation $\overline{[\pi,t]}=\{\overline{\pi}_{n,m}\}_{n,m\in \N}$ introduced in  Proposition \ref{representation of TT(X) proposition}. Here $P_m$ stands for the orthogonal projection onto the subspace \eqref{remarks concerning projections}.
\begin{lemma}\label{lemma almost destroying}
Let   $(\pi,t)$ be  a  representation of $X$  on a Hilbert space $H$,  and let 
$$
J=\{a\in A: \overline{\pi}_{1,1}(\phi(a))=\pi(a)\}.
$$
Then $J_0:=J\cap J(X)$ is an ideal of coisometricity for $(\pi,t)$   and 
\begin{itemize}
\item[i)] an element  $a\in \LL(X^{\otimes m}, X^{\otimes n})$ belongs to $\LL_J(X^{\otimes m}, X^{\otimes n})$ if and only if  $\overline{\pi}_{n,m}(a)$ is supported on $P_{m+1}H$, equivalently,    $\overline{\pi}_{n,m}(a)=\overline{\pi}_{n+1,m+1}(a\otimes 1)$.
\item[ii)] If $a\otimes 1\in \KK(X^{\otimes m+1}, X^{\otimes n+1})$ and $\overline{\pi}_{n,m}(a)$ is supported on $P_{m+1}H$, then $a$ belongs to $\LL_{J_0}(X^{\otimes m}, X^{\otimes n})$ and    in the case  $J_0\subset (\ker\phi)^\bot$ (which is always the case when $(\pi,t)$  is faithful)   $a\in \KK(X^{\otimes m}, X^{\otimes n}J_0)$.
\end{itemize}
\end{lemma}
\begin{proof}
i) Let $a\in \LL_J(X^{\otimes m}, X^{\otimes n})$. Then for any   $x \in X^{\otimes m}$ there exist  $y\in X^{\otimes n}$ and $b\in J$ such that 
$ax =yb$, and we have 
\begin{align*} 
	\overline{\pi}_{n,m}(a)\pi_{m,0}(x) &= \pi_{n,0}(ax)=\pi_{n,0}(yb)=\pi_{n,0}(y) \pi(b)=\pi_{n,0}(y) \overline{\pi}_{1,1}(\phi(b))
	\\
	& =\overline{\pi}_{n+1,1}(y\otimes 1 \cdot \phi(b))=\overline{\pi}_{n+1,1}(yb\otimes 1)= \overline{\pi}_{n+1,1}(ax\otimes 1)
	\\
	&
	= \overline{\pi}_{n+1,m+1}(a\otimes 1)\pi_{m,0}(x).
\end{align*}
Thus  $\overline{\pi}_{n,m}(a)=\overline{\pi}_{n+1,m+1}(a\otimes 1)$ which by  \eqref{right tensor defining relations} is equivalent to   $\overline{\pi}_{n,m}(a)$
being supported on $P_{m+1}H$. 
\\
 Conversely, assume  $a$ is such that $\overline{\pi}_{n,m}(a)$
 is supported on $P_{m+1}H$ (equivalently   $
 \overline{\pi}_{n,m}(a)=\overline{\pi}_{n+1,m+1}(a\otimes 1) $).
Multiplying  $\overline{\pi}_{n,m}(a)$   by $\pi_{0,n}(\flat(x))$, $x\in X^{\otimes n}$, from the left  and by  $\pi_{m,0}(y)$, $y\in X^{\otimes m}$, from the right   we get
 $$
  \pi_{0,n}(\flat(x))\overline{\pi}_{n,m}(a) \pi_{m,0}(y)=\pi(\flat(x)\cdot a \cdot y)=\pi(\langle x,a y\rangle_A).
 $$
Analogously for $\overline{\pi}_{n+1,m+1}(a\otimes 1)$ we have
  \begin{align*} 
  \pi_{0,n}(\flat(x))\overline{\pi}_{n+1,m+1}(a\otimes 1) \pi_{m,0}(y)&=\overline{\pi}_{1,1}(\flat(x)\otimes 1 \cdot (a\otimes1) \cdot y\otimes 1) 
  \\
  &=\overline{\pi}_{1,1}(\langle x ,ay\rangle \otimes1)= \overline{\pi}_{1,1}(\phi(\langle x ,ay \rangle_A)).
 \end{align*}
This implies that $\pi(\langle x,a y\rangle_A)= \overline{\pi}_{1,1}(\phi(\langle x ,ay\rangle_A))$ and hence $\langle x,a y\rangle_A\in J$. By arbitrariness of $x$ and $y$ we conclude  that $a\in \LL_J(X^{\otimes m}, X^{\otimes n})$.
 \\
ii) If $a\otimes 1\in \KK(X^{\otimes m+1}, X^{\otimes n+1})$ and $\overline{\pi}_{n,m}(a)$ is supported on $P_{m+1}H$, then the argument form  the proof of i) shows that  $\pi(\langle x,a y\rangle_A)=\overline{\pi}_{1,1}(\phi(\langle x ,ay \rangle_A))$, and similarly as in the proof of \cite[Lem. 4.2 ii)]{fmr},  one sees that $\phi(\langle x ,ay \rangle_A) \in \KK(X)$. Hence we deduce that $\langle x,a y\rangle_A \in J\cap J(X)$ and as a consequence  $a\in \LL_{J_0}(X^{\otimes m}, X^{\otimes n})$. If additionally $J_0\subset (\ker\phi)^\bot$, then $a\in \KK(X^{\otimes m}, X^{\otimes n}J_0)$ by Lemma \ref{lemma embedding of tensors} iii). 
 \end{proof}
\begin{proposition}\label{theorem with no label}
If  $(\pi,t)$ is  a   representation of $X$ on a Hilbert space $H$ and 
$$
J=\{a\in A: \overline{\pi}_{1,1}(\phi(a))=\pi(a)\},
$$
 then the ideal  of coisometricity for  the right tensor representation $\overline{[\pi,t]}$ of   $\TT_X$ is $\TT_X(J)$.
In particular, if $(\pi,t)$ is   faithful, then   $J$ contained in $(\ker\phi)^\bot$.  
\end{proposition}
    \begin{proof} Apply Proposition  \ref{extensions of representations on Hilbert spaces} iii) and  Lemma \ref{lemma almost destroying} i).
\end{proof}
    \begin{example}\label{endomorphism of unitals2} Let $\TT_\al=\{\al^{n}(1)A\al^m(1)\}_{n,m\in \N}$ be  a right tensor $C^*$-category associated with  an endomorphism 
$\al:A\to A $  of a unital $C^*$-algebra $A$, cf. Examples \ref{partial morphism category}, \ref{partial morphism category1}. If $(\pi,U,B)$ is a  covariant representations   of $\al$ and $\{\pi_{n,m}\}_{n,m\in \N}$ is the associated right tensor representation of $\TT_\al$, see Example \ref{endomorphism of unitals}, then denoting by $\JJ$  the ideal of coisometricity for $\{\pi_{n,m}\}_{n,m\in \N}$ we have  
$$
\JJ(n,m)=\al^{n}(1)J \al^{m}(1) \quad \textrm{where}\quad J=\{a\in A: U^*U\pi(a)=\pi(a)\},
$$
cf.  \cite{kwa-leb1}, \cite[Ex. 1.6]{fmr}. Thus in  terminology of \cite{kwa-leb1} the triple $(\pi,U,B)$ is called a  \emph{covariant representation $(\pi,U,B)$ associated with the ideal} $J$, and we have a one-to-one correspondence between  right tensor representations  of $\TT_\al$ with the ideal of coisometricity $\JJ=\{\al^{n}(1)J \al^{m}(1)\}_{n,m\in \N}$ and covariant representations  of $\al$ associated with  $J=\JJ(0,0)$. Furthermore, if   there exists  a complete transfer operator $\LL$ for  $\al$, cf. Example  \ref{partial morphism category1}, and $J=(\ker\al)^\bot$ (equivalently $\JJ=(\ker\otimes 1)^\bot$),  then   $$ U\pi(a)U^*=\pi(\al(a)),\qquad U^*\pi(a)U=\pi(\LL(a)), $$ that is $(\pi,U,B)$ is a covariant representation in the sense of  \cite{Bakht-Leb}, \cite{kwa-leb3}, \cite{Ant-Bakht-Leb}, cf. \cite{kwa-leb1}. 
In Example \ref{endomorphism of unitals} we have defined the notion of a covariant representation  $(\pi,U,H)$   of a general partial morphism $\varphi:A\to M(A_0)$ in a Hilbert space $H$. Every such  representation give rise to the  right tensor representation  $\{\pi_{n,m}\}_{n,m\in \N}$ of the ideal $\KK_{\varphi}=\{ \KK(X_\varphi^{\otimes m}, X_\varphi^{\otimes n})\}_{n,m\in\N}$ in  $\TT_\varphi$, which in turn extends to the  right tensor representation $\{\overline{\pi}_{n,m}\}_{n,m\in \N}$  of  $\TT_\varphi$. The ideals of coisometricity $\JJ$ and  $\JJ_0$ for $\{\overline{\pi}_{n,m}\}_{n,m\in \N}$ and   $\{\pi_{n,m}\}_{n,m\in \N}$ are established   by the ideals  
  $$
  J=\{a\in A: U^*U\pi(a)=\pi(a)\} \quad\textrm{ and }\quad J_0=J\cap \varphi^{-1}(A_0),
  $$
 respectively.     One may see that,  if   $\varphi$ arises from partial automorphism, cf. Example \ref{partial morphism category1.0}, then $(\pi,U,H)$ is a covariant representation  in the sense of \cite{exel1} iff $J_0=(\ker\varphi)^\bot \cap \varphi^{-1}(A_0)$ (equivalently $\JJ_0=(\ker\otimes 1)^\bot\cap J(\KK_{\varphi})$).
\end{example}

\begin{example}\label{example conerning (E,V)-families} As in Example \ref{E-family example}, let $\TT_E$ be right tensor $C^*$-precategory of a directed graph  $E=(E^0,E^1,r,s)$, $\{\pi_{n,m}\}_{n,m\in \N}$  a right tensor representation  of $\KK_{E}=\{ \KK(X_E^{\otimes m}, X_E^{\otimes n})\}_{n,m\in\N}$ and   $\{p_v:v\in E^0\}$, $\{s_e:e\in E^1\}$ the corresponding Toeplitz-Cuntz-Krieger $E$-family. Then the ideal  of coisometricity $\JJ$  for $\{\pi_{n,m}\}_{n,m\in \N}$ is established via \eqref{ideals of coimoterictiy vs vertices} by a set of vertices $V\subset E^0$ where %.  Actually, $V$  is given by the following relation
$$
p_v=\sum_{s(e)=v}s_e s_e^*,\qquad \textrm{ if and only if } v\in V.
$$
Authors of \cite{mt} called an $E$-family satisfying $p_v=\sum_{s(e)=v}s_e s_e^*$, for  $v\in V$, a
  \emph{Cuntz-Krieger $(E,V)$-family}. Thus we have a one-to-one correspondence between Cuntz-Krieger $(E,V)$-families and right  tensor representations of $\KK_{E}$ coisometric on the ideal $\JJ$ corresponding to $V$. In particular,   
  $$
  \JJ=J(\KK_{E})\cap (\ker\otimes 1)^\bot\quad  \textrm{ iff } \quad V=\{v\in E^0:0<|s^{-1}(v)|<\infty\},
  $$ 
 and if these equivalent conditions hold,  the corresponding  $(E,V)$-family is  called a \emph{Cuntz-Krieger $E$-family} \cite{mt},  \cite{bhrs}, \cite{bprs}.
\end{example}
  \section{$C^*$-algebra $\OO_\TT(\KK,\JJ)$} % of an  ideal $\KK$   relative to an ideal $\JJ$ in  a right tensor $C^*$-precategory $\TT$}
  \label{section four}
  In this section we present a construction of the title object of the  paper - the $C^*$-algebra $\OO_\TT(\KK,\JJ)$. Before  approaching this task, we formulate a universal definition of $\OO_\TT(\KK,\JJ)$ and briefly discuss its relations with relative Cuntz-Pimsner algebras, various crossed products and graph $C^*$-algebras.
  \subsection{Definition  and examples of  $\OO_\TT(\KK,\JJ)$}

The main goal of this section is to construct an object, existence of which is ensured by the  following
\begin{theorem}[\textbf{Characterization of $\OO_\TT(\KK,\JJ)$ as a universal algebra}]
\label{????} For every ideal $\KK$ in a right tensor $C^*$-precategory $\TT$  and every ideal $\JJ$ in  $J(\KK)$ there exists a pair  $\left(\OO_\TT(\KK,\JJ), \iota \right)$ consisting of a $C^*$-algebra $\OO_\TT(\KK,\JJ)$ and a right tensor representation $\iota=\{\iota_{n,m}\}_{n,m\in \N}$  of $\KK$ in $\OO_\TT(\KK,\JJ)$  coisometric on $\JJ$, such that 
 \begin{itemize}
 \item[1)]  the $C^*$-algebra $\OO_\TT(\KK,\JJ)$ is generated by the image of  $\iota=\{\iota_{n,m}\}_{n,m\in \N}$, i.e. 
 $$
 \OO_\TT(\KK,\JJ)=C^*(\{\iota_{n,m}\left(\KK(n,m))\}_{n,m\in \N}\right).
 $$
\item[2)] Every  right tensor representation $\pi=\{\pi_{n,m}\}_{n,m\in \N}$ of  $\KK$ coisometric on $\JJ$ integrates to a representation $\Psi_\pi$ of $\OO_\TT(\KK,\JJ)$ given  by
$$
\Psi_{\pi}(\iota_{n,m}(a))=\pi_{n,m}(a),\qquad a \in \KK(n,m),\,\, n,m \in \N. 
$$
\end{itemize}
Moreover, the pair $\left(\OO_\TT(\KK,\JJ), \iota\right)$ is uniquely determined  in the sense that if $\left(C, \kappa\right)$ is any other pair consisting of a $C^*$-algebra $C$ and a right tensor representation $\kappa=\{\kappa_{n,m}\}_{n,m\in \N}$ of $\KK$  in  $C$ coisometric on $\JJ$, that  satisfies   conditions 1) and  2),  then the mappings
$
\iota_{n,m}(a) \longmapsto  \kappa_{n,m}(a)$, $a \in \KK(n,m)$,
give rise to the canonical isomorphism $\OO_\TT(\KK,\JJ)\cong C$.
\end{theorem}
\begin{proof} The first part of assertion follows from the well known fact on universal $C^*$-algebras, see \cite{blackadar}.  The  second part is straightforward --  by the  universality of    $\left(\OO_\TT(\KK,\JJ), \iota\right)$ and $\left(C, \kappa\right)$ the mappings $\iota_{n,m}(a) \longmapsto  \kappa_{n,m}(a)$ and $\kappa_{n,m}(a) \longmapsto  \iota_{n,m}(a)$,  $a \in \KK(n,m)$, extend to mutually inverse homomorphisms. 
 \end{proof}

%\subsubsection{Representation theorems for $O_\TT(\KK,\JJ)$}
  
\begin{definition}\label{definition of main object}
For any ideals $\KK$, $\JJ$ in a right tensor $C^*$-precategory $\TT$ such that  $\JJ \subset J(\KK)$, the object $
\OO_{\TT}(\KK,\JJ)$ described in  the statement of Theorem \ref{????} will be called  a  \emph{$C^*$-algebra  of the ideal $\KK$ in the right tensor $C^*$-category  $\TT$   relative to the ideal} $\JJ$. 
\end{definition}
 \begin{remark} Uniqueness of $\left(\OO_\TT(\KK,\JJ), \iota\right)$ implies that  there is a canonical circle action $\gamma:S^1\to Aut(\OO_\TT(\KK,\JJ))$ by automorphisms of the $C^*$-algebra $\OO_\TT(\KK,\JJ)$ given by
$$
\gamma_z(\iota_{n,m}(a))=z^{n-m} \iota_{n,m}(a),\qquad a \in \KK(n,m),\,\, z\in S^1.
$$
We shall refer to $\gamma$  as a \emph{gauge action} on $\OO_\TT(\KK,\JJ)$.  
  \end{remark}
\begin{remark}\label{remark on definition of main object}
With analogy to $C^*$-algebras arising from $C^*$-correspondences, see Example \ref{Relative Cuntz-Pimsner algebras par} below, it is reasonable  to call the algebras
$$
\OO_{\TT}(\KK, \{0\}),\qquad  \OO_{\TT}(\KK,\JJ(\KK)\cap (\ker\otimes 1)^\bot),\qquad \OO_{\TT}(\KK,\JJ(\KK))
$$
a \emph{Toeplitz algebra} of $\KK$,   \emph{(Katsura's) $C^*$-algebra of the ideal} $\KK$, and a \emph{Cuntz-Pimsner algebra}  of $\KK$, respectively.
\end{remark}
In the event $\TT$ is a right tensor $C^*$-category, the Doplicher-Roberts algebra $\DR(\TT)$ together with the family of natural maps $\iota_{n+k,n}:\TT(n+k,n)\to \underrightarrow{\,\lim\,\,}\TT(r+k,r)$ (cf. the definition of $\DR(\TT)$ on  page \pageref{Doplichery i cirlce}) is naturally identified with the algebra $\OO_\TT(\TT,\TT)$. We postpone  a relevant discussion until section  \ref{Doplichery i cirlce}, and now present a survey of  examples of $
\OO_{\TT}(\KK,\JJ)$ associated with $C^*$-correspondences.
\begin{example}[Relative Cuntz-Pimsner algebras]\label{Relative Cuntz-Pimsner algebras par}  Let $\TT=\TT_X$ be the right tensor $C^*$-precate\-gory  of a $C^*$-correspondence $X$ and let   $\JJ$  be an ideal in $J(\KK_X)$. %by Theorem \ref{coisometricty coincidence}
We know that $
 \JJ=\KK_{X}(J)=\{\KK( X^{\otimes m},X^{\otimes n}J)\}_{n,m\in \N}
$
where $J=\JJ(0,0)\subset J(X)$, cf. Proposition \ref{Ideals in TT(X)}. In view of Corollary  \ref{coisometricicity of representations} and \cite[Prop. 1.3]{fmr} the algebra $\OO_{\TT}(\KK_X,\JJ)$ coincides with  the \emph{relative Cuntz-Pimsner algebra} of Muhly and Solel \cite{ms}: 
$$
 \OO_{\TT}(\KK_X,\JJ)= \OO(J,X). 
$$
 In particular,  $\OO_{\TT}(\KK_X,\{0\})=\OO(\{0\},X)$ is called a \emph{Toeplitz algebra }of $X$. In the case $\phi$ is injective $\OO_{\TT}(\KK_X,J(\KK_X))=\OO(J(X),X)$ is the algebra originally introduced by Pimsner, and  in the general case the
 algebra 
$$
\OO_{\TT}(\KK_X, J(\KK_X) \cap (\ker \otimes 1)^{\bot})=\OO(J(X)\cap (\ker\phi)^{\bot},X),
$$
is the \emph{$C^*$-algebra $\OO_X$ of the correspondence} $X$ investigated and popularized by  T. Katsura \cite{katsura1},  \cite{katsura}, \cite{katsura2}.
  \end{example}  
 As subclasses of the above algebras we get
\begin{example}[Relative graph algebras] \label{paragraf 3} 
Let $\TT=\TT_E$ be a right tensor $C^*$-precategory  of a directed graph $E$ and let  $\KK$ be the ideal in $\TT$ spanned by the matrix units \eqref{matrix units of graphs}.    Then any  ideal  $\JJ$ in $\KK$ is determined by a set of vertices  $V\subset E^0$, see \eqref{ideals of coimoterictiy vs vertices}. In view of Example \ref{example conerning (E,V)-families}, $V \subset \{v \in E^0: 0<|s^{-1}(v)|<\infty\}$ iff $
\JJ\subset J(\KK)\cap (\ker\otimes 1)^\bot$ and if  this is the case   we have 
$$
\OO_\TT(\KK,\JJ)=C^*(E,V)
$$
where $C^*(E,V)$ is the \emph{relative graph algebra} introduced in \cite[Def. 3.5]{mt}. In particular, $\OO_\TT(\KK,\{0\})=C^*(E,\emptyset)$ is the \emph{Toeplitz algebra of} $E$ as defined in \cite{fr}, and  
$$
\OO_\TT(\KK,J(\KK)\cap (\ker\otimes 1)^\bot)=C^*(E, \{v \in E^0: 0<|s^{-1}(v)|<\infty\})=C^*(E)
$$
where  $C^*(E)$ is the \emph{graph algebra}, cf.  \cite{mt}, \cite{bhrs}, \cite{bprs}. We extend \cite[Def. 3.5]{mt} and define  $C^*(E,V):=\OO_\TT(\KK,\JJ)$ for any $V \subset \{v \in E^0: |s^{-1}(v)|<\infty\}$.   \end{example}

\begin{example}[Crossed products by  partial morphisms]\label{paragraf 1}
Let $\TT=\TT_\varphi$ be the right tensor $C^*$-precategory   associated with a partial morphism $\varphi$. For the ideal $\KK$ described by \eqref{ideal for partial morhpsims} we have
$$
\OO_{\TT}(\KK, J(\KK) \cap (\ker \otimes 1)^{\bot})=A\rtimes_{\varphi} \N
$$
where $A\rtimes_{\varphi} \N$ is the \emph{crossed product by partial morphism} defined in \cite{katsura1}. In particular, if $\varphi$  arises from partial automorphism, this algebra coincides with the Exel's crossed product by partial automorphism \cite{exel1}.
   \end{example}

\begin{example}[Partial isometric crossed products by  endomorphisms] \label{paragraf 2} 
Let $\TT$ be the right tensor $C^*$-category  $\TT_\al=\{\al^{n}(1)A \al^m(1)\}_{n,m\in \N}$ associated with  an endomorphism $\al:A\to A$   of a unital  $C^*$-algebra $A$, cf. Example \ref{partial morphism category1}.  Then every ideal $\JJ$ in $\TT$ is of the form $\JJ=\{\al^{n}(1)J \al^m(1) \}_{n,m\in \N}$ for  an ideal $J$  in $A$, and
$$
\OO_\TT(\TT,\JJ)=C^*(A,\al;J)
$$
 where $C^*(A,\al;J)$ is the \emph{relative crossed product} considered in  \cite{kwa-leb1}, cf. Example \ref{endomorphism of unitals2}. Thereby, see \cite{kwa-leb1},  for   $J=A$, $J=(\ker\al)^\bot$ and  $J=\{0\}$ one arrives at partial-isometric crossed product from  \cite{Adji_Laca_Nilsen_Raeburn}, \cite{kwa} and \cite{Lin-Rae}, respectively.
 Furthermore, if there exists a complete transfer operator $\LL$ for $\al$, then
  $$
C^*(A,\al;(\ker\al)^\bot)= \OO_\TT(\TT, (\ker\otimes 1)^\bot)
$$
is the crossed product investigated in \cite{kwa-leb3}, \cite{Ant-Bakht-Leb}.
A distinctive property of $C^*$-algebras of this type  is that they are   generated by a homomorphic image of $A$ and a   partial isometry. 
  \end{example}
The last two examples  indicate   that it is reasonable to adopt the following   definition  which embraces  all the   crossed products  mentioned above. %(we shall use it in Example \ref{structure theorem for partial morphisms}).
\begin{definition}\label{definition of crossed product by partial morphisms}
Let  $\varphi:A\to M(A_0)$ be a partial morphism and let $J$ be an ideal   in $\varphi^{-1}(A_0)$.  We define the \emph{relative crossed product} $C^*(\varphi;J)$ of  $\varphi$ relative to $J$ to be  the $C^*$-algebra $\OO(J,X_\varphi)$.
\end{definition}
\subsection{Construction of $\OO_\TT(\KK,\JJ)$}\label{section construction}
Let $\KK$ and $\JJ\subset \KK$ be ideals in a right tensor $C^*$-precategory $\TT$. A construction of the $C^*$-algebra  $\OO_\TT(\KK,\JJ)$ will consist of two steps.
Firstly, we describe a  matrix calculus that yields   a purely algebraic structure - a $^*$-algebra $\M_\TT(\KK)$. Secondly, we  use $\JJ$ to define a seminorm $\|\cdot \|_\JJ$ on  $\M_\TT(\KK)$  such that completing the quotient space  $\M_\TT(\KK)/\|\cdot \|_\JJ$ yields $\OO_\TT(\KK,\JJ)$.
\subsubsection{An algebraic framework  for $\OO_\TT(\KK,\JJ)$}
\label{2}
 We let $\M_\TT$ be the set of all infinite matrices $\{a_{n,m}\}_{n,m\in\N} $  where  $a_{n,m}\in \TT(n,m)$, and denote by  $\M_\TT(\KK)$ the subset  of $\M_\TT$ consisting of matrices
$\{a_{n,m}\}_{n,m\in\N} $    such that 
$$
a_{n,m}\in \KK(n,m), \qquad  n,\,m\in \N,
$$
and there is at most finite number of elements $a_{n,m}$ which are non-zero. 
  We define the
addition, multiplication by scalars, and
involution on $\M_\TT(\KK)$  in a quite natural manner: for $a=\{a_{n,m}\}_{n,m\in\N}$ and $b=\{b_{n,m}\}_{n,m\in\N}$ we put
\begin{gather}\label{add1}
 (a+b)_{n,m}:=a_{n,m}+b_{n,m},\\[6pt]
\label{mulscal1}
 (\lambda a)_{n,m}:=\lambda a_{n,m}\\[6pt]
\label{invol1}
 (a^*)_{n,m}:=a_{m,n}^*.
\end{gather}
A  multiplication "$\star$" on
$\M_\TT(\KK)$  is more involved. We set
\begin{equation}\label{star1}
a\star b:= a\cdot\sum_{k=0}^\infty  \Lambda^k(b)+ \sum_{k=1}^\infty \Lambda^k(a)\cdot b
\end{equation}
where "$\cdot$" stands for the standard multiplication of matrices and   $\Lambda:
\M_\TT(\KK)
\rightarrow \M_\TT$ is defined
 to act as follows
 \begin{equation}\label{Lambda}
\Lambda(a)= \left(
\begin{array}{cc c c c }
 0      &         0           &           0         &            0       & \cdots  \\
 0      &  a_{0,0}\otimes 1  &  a_{0,1}\otimes 1  &  a_{0,2}\otimes 1 & \cdots  \\
 0      &  a_{1,0}\otimes 1 & a_{1,1}\otimes 1  &  a_{1,2}\otimes 1 & \cdots  \\
 0      &  a_{2,0}\otimes 1 &  a_{2,1}\otimes 1 &  a_{2,2}\otimes 1 & \cdots  \\
\vdots  &       \vdots        &       \vdots        &       \vdots       & \ddots
\end{array}\right),
\end{equation}
that is  $\Lambda(a)_{n,m}=a_{n-1,m-1}\otimes 1$, for $n,m> 1$, and $\Lambda(a)_{n,m}=0$ otherwise. 
Note that even though the entries of $\Lambda(a)$ need not belong to $\KK$, entries of $a\star b$ are in $\KK$  (because $\KK$ is an ideal in $\TT$)  and hence $a\star b$ is well defined.
\par
We denote by $\iota_{n,m}^\M:\TT(n,m) \to \M_\TT$, $n,m\in\N$, the natural embeddings, that is   $\iota_{n,m}^\M(a)$ is  a matrix $\{a_{i,j}\}_{i,j\in \N}$  satisfying $a_{i,j}= \delta_{i,n} \delta_{j,m} a$  where $\delta_{k,l}$ is the Kronecker symbol.
\begin{proposition}\label{representations  extensions  proposition}\label{*-algebra proposition}
The set $\M_\TT(\KK)$ with operations \eqref{add1},
\eqref{mulscal1}, \eqref{invol1}, \eqref{star1} becomes an algebra with involution. Moreover,
\begin{itemize}
\item[i)]  the family  $\{\iota_{n,m}^\M\}_{n,m\in \N}$ forms a right tensor representation of  $\KK$ in  $\M_\TT(\KK)$.
\item[ii)] we have a one-to-one correspondence, established by relations 
$$
\Psi( \iota_{n,m}^\M (a_{n,m}))=\pi_{n,m}(a_{n,m}),
$$
between right tensor representation $\{\pi_{n,m}\}_{n,m\in \N}$ of  $\KK$  and representations $\Psi$ of $\M_{\TT}(\KK)$. We shall denote by $\Psi_{\pi}$ the representation
$\Psi_{\pi}( \{a_{n,m}\}_{n,m\in\N})=\sum_{n,m=0}^\infty \pi_{n,m}(a_{n,m})$
corresponding to $\pi=\{\pi_{n,m}\}_{n,m\in \N}$.
\end{itemize}
\end{proposition}
\begin{proof}  The mapping $\sum_{k=0}^\infty  \Lambda^k:\M_\TT(\KK)\to \M_\TT$ embeds $\M_\TT(\KK)$ equipped with operations \eqref{add1},
\eqref{mulscal1}, \eqref{invol1}, \eqref{star1} into the $^*$-algebra $\M_\TT$  equipped with standard matrix operations. Indeed,  a moment of thought shows  that  $\sum_{k=0}^\infty  \Lambda^k:\M_\TT(\KK)\to \M_\TT$ is an injective linear map that preserves involution, and
\begin{align*}
\sum_{k=0}^\infty  \Lambda^k(a\star b) &= \sum_{k=0}^\infty  \Lambda^k \left(a\cdot\sum_{l=0}^\infty  \Lambda^l(b)+ \sum_{l=1}^\infty \Lambda^k(a)\cdot b\right)
\\
 &=\sum_{k=0}^\infty \left( \Lambda^k (a)\sum_{l=0}^\infty \Lambda^{k+l}(b)+  \sum_{l=1}^\infty\Lambda^{k+l}(a)\cdot \Lambda^{k}(b)\right)
 \\
 &=\sum_{k,l=0}^\infty  \Lambda^k(a)\Lambda^{l}(b)=\left(\sum_{k=0}^\infty  \Lambda^k(a)\right)\left(\sum_{l=0}^\infty  \Lambda^l(b) \right).
\end{align*}
Hence $\sum_{k=0}^\infty  \Lambda^k:\M_\TT(\KK)\to \M_\TT$ is an injective $^*$-homomorphism and therefore $\M_\TT(\KK)$ is a $^*$-algebra. 
\\
i) Checking that $\{\iota_{n,m}^\M\}_{n,m\in \N}$ is a right tensor representation is straightforward. In particular, the relation $
\iota_{n,m}^\M (a_{n,m})\star  \iota_{m+k,l}^\M (b_{m+k,l}))=\iota_{n+k,l}^\M ((a_{n,m}\otimes 1^k)b_{m+k,l})
$ follow directly from \eqref{star1}.
\\
ii) If $\{\pi_{n,m}\}_{n,m\in \N}$ is a right tensor representation  of  $\KK$, then $\Psi( \iota_{n,m}^\M (a_{n,m}))=\pi_{n,m}(a_{n,m})$ determines uniquely a linear map $\Phi$. Since
$$
\Psi( \iota_{n,m}^\M (a_{n,m})^*)=\Psi( \iota_{m,n}^\M (a_{n,m}^*))=\pi_{m,n}(a_{n,m}^*)=\pi_{n,m}(a_{n,m})^*
$$
and 
\begin{align*}
\Psi( \iota_{n,m}^\M (a_{n,m})\star  \iota_{m+k,l}^\M (b_{m+k,l}))&=\Psi( \iota_{n+k,l}^\M ((a_{n,m}\otimes 1^k)b_{m+k,l})
 \\
 &=\pi_{n+k,l}((a_{n,m}\otimes 1^k) b_{m+k,l})
 \\
 &=\pi_{n,m} (a_{n,m})  \pi_{m+k,l}(b_{m+k,l})
\\
&=\Psi( \iota_{n,m}^\M (a_{n,m}))\Psi(  \iota_{m+k,l}^\M (b_{m+k,l})),
\end{align*}
it follows that  $\Psi$ is a representation of $\M_{\TT}(\KK)$. Conversely, if $\Psi$ is a representation of $\M_{\TT}(\KK)$, then $\pi=\{\pi_{n,m}\}_{n,m\in \N}$ as a ``composition'' of $\Psi$ and the right tensor representation $\{\iota_{n,m}^\M\}_{n,m\in \N}$ is a   right tensor representation.
\end{proof}

\subsubsection{Grading and seminorms $\|\cdot\|_{\JJ}$ in $\M_\TT(\KK)$}
We define  a seminorm  $\|\cdot\|_{\JJ}$ in $\M_\TT(\KK)$ using a natural grading of  $\M_\TT(\KK)$.  Also we take an opportunity to show that this  seminorm   could be   defined within the inductive limit frames similar to   Doplicher-Roberts approach.
For that purpose, for $r\in \N$, $k\in \Z$, $r+k\geq 0$,  we put
$$
\M(r+k,r):=\{\{a_{n,m}\}_{m,n\in\N}\in \M_\TT(\KK): \,\,a_{n,m}\neq 0  \Longrightarrow\,\, n-m=k,\,\, m\leq r\}.
$$
Hence an element  $a\in \M_\TT(\KK)$ is in  $\M(r+k,r)$ iff it is of the form
\begin{center}\setlength{\unitlength}{1mm}
 \begin{picture}(100,30)(62,-13)
    
  \put(92,0){, $\,\,$ if $\,\,k \geq 0$, or }
 \put(61,0){$
 \left(\begin{array}{c}\begin{xy}
\xymatrix@C=-1pt@R=4pt{
    &      \, \,   &    \qquad 0   \\
     &      \,     &     \\
                &     \,      &   \qquad \quad\,\,\,    \\
                      &     \,      &       \\
                  &     \,      &       \\
                   0 \,\, &       &    \\
                    &     \,      &       
                   }
  \end{xy}
  \end{array}\right)$}
  \scriptsize
  \put(64,5){$a_{k,0}$}
   \put(75,-6){$a_{r+k,r}$}
\put(58,10.5){$\begin{cases}  & \\ 
 &  \end{cases}$}
 \put(56,10.5){$k$}
 \put(58,-1){$\begin{cases}  & \\ 
 & \\
 & \\
 & \end{cases}$}
 \put(51.5,-1){$r+1$}
  \qbezier[8](69,3)(72.5,0)(76,-3)
       \end{picture}

\begin{picture}(50,0)(-38,-15)
\put(-7,0){$
 \left(\begin{array}{c}\begin{xy}
\xymatrix@C=-1pt@R=3pt{
   &      \, \,   &    \qquad \qquad\,\,\,\,\,\,\,\,\,0 \,   \\
     &      \,     &       \\
          &      \,     &       \\
               &      \,     &       \\
      0 &        &      \\
           &       &  \qquad 
        }
  \end{xy}
  \end{array}\right)
 $  }
 \scriptsize
  \put(3,10){ $a_{0,-k}$}
    \put(14,-1){ $a_{r+k,r}$}
  \put(-4,13){$\overbrace{\quad\quad\,\,\,}^{-k}$} \put(4,13){$\overbrace{\qquad\qquad\qquad\,\,}^{r+k+1}$}

 \qbezier[8](10, 8)(13.5,5)(17,2)
  
\normalsize
   \put(32.5,2){, if $\,\,k< 0$.}
    \end{picture}
      
     \end{center}
For every $k\in \Z$  we get an increasing family $\{\M(r+k,r)\}_{r\in \N,\,\, r+k \geq 0}$ of  subspaces of  $\M_\TT(\KK)$, and the spaces 
$$
\M^{(k)}_\TT(\KK):=\bigcup_{r\in \N, \atop r+k \geq 0}\M(r+k,r),\qquad k\in \Z,
$$
define a $\Z$-grading on  $\M_\TT(\KK)$. In particular 
$$
\M_\TT(\KK)=\bigoplus_{k\in \Z} \M^{(k)}_\TT(\KK).
$$
The first step in defining $\|\cdot\|_{\JJ}$ is to  equip  the family $\{\M(n,m)\}_{n,m\in \N}$ with the structure of  a right tensor $C^*$-precategory. We recall that $q_\JJ$ denotes  the quotient map from $\TT$ onto $\TT/\JJ$, cf. Proposition \ref{quotient C*-category}. 
\begin{theorem}[Construction of the right tensor $C^*$-precategory $\KK_\JJ$]\label{norm formulas on 0-spectral subspace}
For every $r\in \N$, $k\in \Z$, $r+k\geq 0$, the   formula 
$$
\|a\|_{r+k,r}^\JJ:=\max \left\{\max_{s=0,...,r-1}\big\{ \|q_\JJ\Big(\sum_{i=0, \atop i+k\geq 0}^{s}
   a_{i+k,i}\otimes 1^{s-i}\Big)\|\big\}, \big\|\sum_{i=0, \atop i+k\geq 0}^{r}
   (a_{i+k,i}\otimes 1^{r-i})\big\| \right\}
$$
defines a  seminorm on $\M(r+k,r)$,   such that the family of quotients spaces $$\KK_\JJ:=\{\M(n,m)/\|\cdot\|_{n,m}^\JJ\}_{n,m\in \N}$$
form a right tensor $C^*$-precategory  with a right tensoring $\otimes_\JJ 1$ induced by the inclusions $\M(n,m)\subset \M(n+1,m+1)$,    $m,n\in \N$. Moreover  the right tensoring $\otimes_\JJ 1$ is faithful iff  $\JJ\subset (\ker\otimes 1)^\bot$.
\end{theorem}
\begin{proof} By \eqref{invol1} we have  $\M(n,m)^*=\M(m,n)$, $n,m\in \N$. To  show
$$\M(r+k,r) \star \M(r,r+l)\subset  \M(r+k,r+l), \qquad  r\in \N, k,l\in \Z,$$ 
 note that $\M(r+k,r)=\spane\{\iota_{i+k,i}^\M(\KK(i+k,i)): i=0,...,r,\,\, i+k\geq 0\}$ and  for $a\in\KK(i+k,i)$ and $b\in \KK(j,j+l)$, $i,j=0,...,r$,  the element
\begin{equation}\label{product of monomials}
 \iota_{i+k,i}^\M(a)\star \iota_{j,j+l}^\M(b)=\begin{cases}
  \iota_{j+k,j+l}^\M((a\otimes 1^{j-i})b), & \textrm{ if } i\leq j,\\
  \iota_{i+k,i+l}^\M(a(b\otimes 1^{i-j})), & \textrm{ if }j< i,
 \end{cases}\quad 
\end{equation}
 is in $\M(r+k,r+l)$. 
 Thus   $\{\M(n,m)\}_{n,m\in \N}$ is a $^*$-precategory. 
 \\
It is clear that functions $\|\cdot\|_{r+k,r}^\JJ$ are seminorms on $\M(r+k,r)$. To see that they are  submultiplicative   let $a=\{a_{n,m}\}_{n,m\in \N}\in \M(r+k,r)$, $b=\{b_{n,m}\}_{n,m\in \N} \in\M(r,r+l)$  and  $s=0,1,...,r+l$. We adopt to the convention that $a_{n,m}\otimes 1^p=0$ whenever $m$, $n$, or $p$ is less than $0$. Using \eqref{product of monomials} one gets
 $$
 (a\star b)_{p+k-l,p}=\sum_{i=0}^{p-l} (a_{i+k,i}\otimes 1^{p-l-i})b_{p-l,p} + \sum_{j=0}^{p-l-1} (a_{p+k-l,p-l}(b_{j,j+l})\otimes 1^{p-l-j}),
 $$
for $p=0,...,s$, and  after rearranging the addends one obtains 
\begin{equation}\label{addends rearranged}
\sum_{p=0}^{s}(a\star b)_{p+k-l,p}\otimes 1^{s-p}  =\left(\sum_{p=0}^{s}  a_{p+k-l,p-l}\otimes 1^{s-p}\right) \left(\sum_{p=0}^{s} b_{p-l,p}\otimes 1^{s-p}\right).
\end{equation}
%\begin{align*} &\sum_{p=0}^{s}(a\star b)_{p+k-l,p}\otimes 1^{s-p}  \\  &=\sum_{p=0}^{s} \sum_{i=0}^{p-l} (a_{i+k,i}\otimes 1^{s-l-i})b_{p-l,p}\otimes 1^{s-p} + \sum_{i=0}^{p-l-1} (a_{p+k-l,p-l}\otimes 1^{s-p} )b_{i,i+l}\otimes 1^{s-l-i}   \\  &=\left(\sum_{p=0}^{s}  a_{p+k-l,p-l}\otimes 1^{s-p}\right) \left(\sum_{p=0}^{s} b_{p-l,p}\otimes 1^{s-p}\right). \end{align*}
Accordingly, since  $q_\JJ$ preserves the operations in $\TT$ and  by submultiplicativity of  norm in $\TT$ we see that $\|q_\JJ \big(\sum_{p=0}^{s}(a\star b)_{p+k-l,p}\otimes 1^{s-p}\big)\|$ is not larger than
$$
 \|q_\JJ\Big(\sum_{p=0}^{s}  a_{p+k-l,p-l}\otimes 1^{s-p}\Big)\|\cdot  \|q_\JJ\Big(\sum_{p=0}^{s} b_{p-l,p}\otimes 1^{s-p}\Big)\| 
 $$
 and  
$$
\|\sum_{p=0}^{r+l}(a\star b)_{p+k-l,p}\otimes 1^{s-p}\|\leq \|\sum_{p=-l}^{r}  a_{p+k,p}\otimes 1^{r-p}\|\cdot  \|\sum_{p=0}^{r+l} b_{p-l,p}\otimes 1^{r+l-p}\|.   $$ 
Therefore 
 $$
\|a\star b\|_{r+k,r+l}^\JJ
 \leq \|a\|_{r+k,r}^\JJ \|b\|_{r,r+l}^\JJ.
$$
Putting in  \eqref{addends rearranged}, $k=l$,  $b= a^*$ and using the $C^*$-equality in $\TT$ one sees that the above inequalities become equalities,  and therefore 
    $$
    \|a\star a^{*}\|_{r+k,r+k}^\JJ=\big(\|a\|_{r+k,r}^\JJ\big)^2.
    $$
  Now let us consider the quotient  space  $\M(r+k,r)/\|\cdot\|_{r+k,r}^\JJ$ with its quotient norm. Since  $q_\JJ$ and $\otimes 1$ are contractive  one gets 
$$
    \|a\|_{r+k,r}^{\JJ} \leq r \max_{i=0,...,r \atop i+k \geq 0} \|a_{i+k,i}\|, \qquad a\in\M(r+k,r), 
  $$
  and it follows that   we have a bounded linear epimorphism
    $$
    \bigoplus_{i=0\atop i+k\geq 0}^{r}\KK(i+k,i) \ni \{a_{i+k,i} \}_{i=0\atop i+k\geq 0}^r \longmapsto \sum_{i=0, \atop i+k\geq 0}^{s}
   \iota^\M_{i+k,i}(a_{i+k,i}) \in \M(r+k,r)/\|\cdot\|_{r+k,r}^\JJ
    $$
   from  the direct sum  of Banach spaces $\bigoplus_{i=0\atop i+k\geq 0}^r\KK(i+k,i)$ onto $\M(r+k,r)/\|\cdot\|_{r+k,r}^\JJ$ (we abuse the notation concerning $\iota^\M_{i+k,i}$). Hence by the closed graph  theorem  the space  $\M(r+k,r)/\|\cdot\|_{r+k,r}^\JJ$ is a Banach space  (isomorphic to the quotient of $\bigoplus_{i=0\atop i+k\geq 0}^r\KK(i+k,i)$ by the kernel of the introduced epimorphism).  
Accordingly, in view of the above  $\KK_\JJ=\{\M(n,m)/\|\cdot\|_{n,m}^\JJ\}_{n,m\in \N}$ is a right tensor $C^*$-precategory.
\par
To see  that the inclusion  $\M(r+k,r)\subset \M(r+k+1,r+1)$ factors through to the mapping
$$
\M(r+k,r)/\|\cdot\|_{r+k,r}^\JJ \stackrel{\,\,\otimes _\JJ 1\,\,}{\longrightarrow}  \M(r+k+1,r+1)/\|\cdot\|_{r+k+1,r+1}^\JJ
$$ 
let  $a\in \M(r+k,r)$ be such that   $\|a\|_{r+k,r}^\JJ=0$. Then 
   $\sum\limits_{i=0}^{s}
   a_{i+k,i}\otimes 1^{s-i}\in \JJ(s+k,s)$, for all $s=0,1,....,r$, where  for $s=r$ we actually have $\sum\limits_{i=0}^{r}
   a_{i+k,i}\otimes 1^{r-i}=0$. Moreover, 
$$
\sum_{i=0}^{r+1}
  a_{i+k,i}\otimes 1^{r+1-i} =\sum_{i=0}^{r}
  a_{i+k,i}\otimes 1^{r+1-i}=\Big(\sum_{i=0}^{r}
  a_{i+k,i}\otimes 1^{r-i}\Big)\otimes 1= 0.
$$
Thus $\|a\|_{r+1,r+1}^\JJ=0$ and $\otimes_\JJ 1$ is a well defined   right tensoring on $\KK_\JJ$.
\\
Assume  now  that  $\JJ\subset (\ker\otimes1)^\bot$ and let  $a\in \M(r+k,r)$ be such that  $\|a\|_{r+k+1,r+1}^\JJ=0$. Then  $\sum\limits_{i=0}^{s}
   a_{i,i+k}\otimes 1^{s-i}\in \JJ(s+k,s)$,  for each  $s=0,1,....,r$, and
$$
\sum_{i=0}^{r+1}
   a_{i+k,i}\otimes 1^{r+1-i} =\Big(\sum_{i=0}^{r}
   a_{i+k,i}\otimes 1^{r-i}\Big)\otimes 1=0.
$$
In particular,  $\sum\limits_{i=0}^{r}
   a_{i+k,i}\otimes 1^{r-i}\in \ker (\otimes 1) (r+k,r)$ and $\sum\limits_{i=0}^{r}
   a_{i+k,i}\otimes 1^{r-i}\in \JJ(r+k,r)$. Since $ \JJ \cap \ker (\otimes 1)= \{0\}$, we get  $\sum\limits_{i=0, \atop i+k\geq 0}^{r}
   a_{i+k,i}\otimes 1^{r-i}=0$ and thereby $\|a\|_{r+k,r}^\JJ=0$.  This proves the injectivity of $\otimes _\JJ 1$. 
 \\
 Conversely, if $\JJ\nsubseteq (\ker\otimes1)^\bot$, then  there exists  $a\in  (\JJ\cap \ker\otimes1) (n,m)$ such that $a\neq 0$. Then $\|\iota_{n,m}(a)\|_{n,m}^\JJ\neq 0$ and  $\|\iota_{n,m}(a)\|_{n+1,m+1}^\JJ= 0$, that is $\otimes_\JJ 1$ is not injective.
\end{proof}

The foregoing statement allow us to   apply the "Doplicher-Roberts method" of constructing $C^*$-algebras from right tensor $C^*$-categories presented in \cite{dr}. Namely, for  $k\in \Z$, we let 
$ \OO^{(k)}_\TT(\JJ,\KK):=\underrightarrow{\,\lim\,\,}\KK_\JJ(r+k,r)$  be the Banach space inductive limit of the inductive sequence   $$  \label{dr-definition of our algebras}
\KK_{\JJ}(r+k,r) \stackrel{\otimes_\JJ 1}{\longrightarrow}\KK_\JJ(r+k+1,r+1) \stackrel{\otimes_\JJ 1}{\longrightarrow} \KK_\JJ(r+k+2,r+2) \stackrel{\otimes_\JJ 1}{\longrightarrow} ...
$$
defined for $r\in \N$, $k+r \geq 0$.  The algebraic direct sum $\bigoplus_{k\in \Z} \OO^{(k)}_\TT(\KK,\JJ)$ has a natural structure of $\Z$-graded $^*$-algebra  and we define  $\OO_{\TT}(\KK,\JJ)$ to be   the $C^*$-algebra obtained by completing $\bigoplus_{k\in \Z} \OO^{(k)}_\TT(\KK,\JJ)$ in the unique $C^*$-norm for which the automorphic action defined by the grading is isometric, see \cite[Thm. 4.2]{dr}.
\\
We may also consider  $\OO^{(k)}_\TT(\KK,\JJ)$ as a completion of the quotient space  $\M^{(k)}_\TT(\KK)/ \|\cdot\|_\JJ$ where $\|\cdot\|_\JJ$ is a seminorm described below, and thus obtain another  construction of  $\OO_\TT(\KK,\JJ)$.
\begin{proposition}[Construction of $\OO_\TT(\KK,\JJ)$]\label{propostion nie udowodnione}
The  formula 
\begin{align*}
\|a\|_{\JJ}=\sum_{k\in\Z} \lim_{r \to \infty} \max \left\{ \max_{ s=0,...,r-1}\big\{ \|q_\JJ\Big(\sum_{i=0, \atop i+k\geq 0}^{s}
    \right.& a_{i+k,i}\otimes 1^{s-i}\Big)\|\big\},\, 
   \\
   & \left. \big\|\sum_{i=0, \atop i+k\geq 0}^{r}
   (a_{i+k,i}\otimes 1^{r-i})\big\| \right\}
\end{align*}defines a submultiplicative $^*$-seminorm on $\M_\TT(\KK)$ such that the  enveloping $C^*$-algebra of the quotient $^*$-algebra $\M_\TT(\KK)/\|\cdot\|_\JJ$ is naturally isomorphic to $\OO_\TT(\KK,\JJ)$.
\end{proposition}
\begin{proof}
It suffices to notice  that  the norm of the image of an element $a\in \M^{(k)}_\TT(\KK)$ in the inductive limit space $\OO^{(k)}_\TT(\KK,\JJ)$ coincides with the value $\|a\|_{\JJ}$, that is 
$\|a\|_{\JJ}=\lim_{r\to \infty}\|a\|_{r+k,r}^\JJ$. In particular, $\|\cdot \|_{\JJ}$ is submultiplicative $^*$-seminorm on $\M_\TT(\KK)$ and  enveloping $C^*$-norm on $\M_\TT(\KK)/\|\cdot\|_\JJ$ satisfy the conditions of  \cite[Thm. 4.2]{dr}.
\end{proof}
The only thing  left for us  to prove is  that $\OO_\TT(\KK,\JJ)$  satisfies the universal conditions presented in Theorem  \ref{????}.

\begin{theorem}[Universality of $\OO_\TT(\KK,\JJ)$]\label{universality checking}
Let $\JJ$ be an ideal in $\TT$ such that $\JJ\subset J(\KK)$. 
The  right tensor representation   described in Proposition \ref{*-algebra proposition} i)  composed with the quotient map that arise  in  Proposition \ref{propostion nie udowodnione} yields the  right tensor representation $\iota=\{\iota_{n,m}\}_{n,m\in \N}$ of $\KK$  in  $\OO_\TT(\KK,\JJ)$ which is coisometric on  $\JJ$. The pair $\left(\OO_\TT(\KK,\JJ), \iota\right)$ satisfies the following conditions: 
 \begin{itemize}
 \item[1)]  the $C^*$-algebra $\OO_\TT(\KK,\JJ)$ is generated by the image of the representation $\{\iota_{n,m}\}_{n,m\in \N}$, i.e. 
 $$
 \OO_\TT(\KK,\JJ)=C^*(\{\iota_{n,m}\left(\KK(n,m))\}_{n,m\in \N}\right).
 $$
\item[2)] for every  right tensor representation $\pi=\{\pi_{n,m}\}_{n,m\in \N}$ of  $\KK$ which is coisometric on $\JJ$ the  representation $\Psi_\pi$ described in Proposition \ref{representations  extensions  proposition} factors through to the representation of $\OO_\TT(\KK,\JJ)$ such that
$$
\Psi_{\pi}(\iota_{n,m}(a))=\pi_{n,m}(a),\qquad a \in \KK(n,m). 
$$
\end{itemize}
\end{theorem}
\begin{proof}
Since the quotient map in  Prop. \ref{propostion nie udowodnione} is a $^*$-homomorphism and  $\{\iota_{n,m}^\M\}_{n,m\in \N}$ is a right tensor representation, their composition is a    right tensor representation $\iota=\{\iota_{n,m}\}_{n,m\in \N}$. To see that $\iota$ is coisometric on $\JJ$ take $a \in J(\KK)(n,m)$, $n,m\in \N$, and notice that
$$
\begin{array}{rcl}
\iota_{n,m}(a)=\iota_{n+1,m+1}(a\otimes 1) & \,\,\,\Longleftrightarrow\, \,\,& \|\iota_{n,m}^\M(a) -\iota_{n+1,m+1}^\M(a\otimes 1)\|_\JJ=0 \\[4pt]
                                        & \Longleftrightarrow & a \in \JJ(n,m)
	\end{array} 
$$
Item 1) is clear. Let us prove 2). By the form of  $\|\cdot\|_{\JJ}$ (it is the sum of seminorms on spectral subspaces) it is enough to show that $\Psi_{\pi}$ factors through to a representation of the subspace $\M_\TT^{(k)}(\KK)/\|\cdot\|_\JJ$, $k\in\Z$. %, cf. the proof of Proposition \ref{norm formulas on 0-spectral subspace}. 
For that purpose let $a\in \M_\TT^{(k)}(\KK)$ be such that  $\|a\|_{\JJ}=0$. Then
$$
\sum_{i=0, \atop i+k\geq 0}^{s}
   a_{i+k,i}\otimes 1^{s-i}\in \JJ(s,s),\,\,\,\,s\in \N \qquad \lim_{r\to\infty} \|\sum_{i=0, \atop i+k\geq 0}^{r}
   a_{i+k,i}\otimes 1^{r-i}\|=0.
$$ 
Thus,  for sufficiently large $N$,  we have 
\begin{align*}
\|\Psi_{\pi}( a)\|&=\|\sum_{i=0, \atop i+k\geq 0}^N \pi_{i+k,i}(a_{i+k,i})\|=\|\pi_{N+k,N}(\sum_{i=0, \atop i+k\geq 0}^N a_{i+k,i}\otimes 1^{N-i})\| 
\\
 &\leq \|\sum_{i=0, \atop i+k\geq 0}^N a_{i+k,i}\otimes 1^{N-i}\| \longrightarrow 0, \qquad \textrm{as}\qquad N\longrightarrow \infty.
\end{align*}
Hence item 2) follows.
%\\ The form of the kernel of the representation $\{\iota_{n,m}\}_{n,m\in \N}$ is  clear. In particular it contains the ideal $\JJ \cap \ker\otimes 1$. In case $\JJ \cap \ker\otimes 1=\{0\}$ the faithfulness of $\{\iota_{n,m}\}_{n,m\in \N}$ follows from Proposition \ref{norm formulas on 0-spectral subspace}.
 \end{proof}
\begin{remark}
In the process of construction of  $\OO_\TT(\KK,\JJ)$, until the above statement, we did not require that $\JJ\subset J(\KK)$. Actually  for any ideal $\JJ$ in $\TT$, Proposition \ref{propostion nie udowodnione} defines the $C^*$-algebra $\OO_\TT(\KK,\JJ)$. However, without the  condition  $\JJ\subset J(\KK)$ the relationship  between    representations  of $\KK$  and $\OO_\TT(\KK,\JJ)$  breaks down.
\end{remark}

\subsection{Immediate corollaries of  the construction of $\OO_\TT(\KK,\JJ)$}
%The construction   presented in the previous subsections  to reveals a  piece of the structure of $\OO_\TT(\KK,\JJ)$. In particlar, we have
\begin{theorem}[Norm of elements in spectral subspaces] \label{proposition of norm formulas} If we let 
$\left(\OO_\TT(\KK,\JJ), \iota \right)$ be the universal pair as in Theorem \ref{????}, then for every $k\in \Z$, the  elements of the form
\begin{equation}\label{element a form}
a= \sum_{j=0, \atop j+k\geq 0}^{r_0}
   \iota_{j+k,j}(a_{j+k,j}), \qquad  r_0\in \N,
\end{equation}
constitute a dense subspace of the $k$-th spectral subspace of $\OO_\TT(\KK,\JJ)$, and the norm  $\|a\|$ of $a$ is given by the following  limit
$$
\lim_{r\to \infty} \max \left\{\max_{s=0,...,r-1}\big\{ \|q_\JJ\Big(\sum_{i=0, \atop i+k\geq 0}^{s}
   a_{i+k,i}\otimes 1^{s-i}\Big)\|\big\},\, \big\|\sum_{i=0, \atop i+k\geq 0}^{r}
   (a_{i+k,i}\otimes 1^{r-i})\big\| \right\}.
$$
If $\JJ\subset (\ker\otimes 1)^\bot$ the above formula for $\|a\|$ reduces to the following one
$$
\max \left\{\max_{s=0,...,r_0-1}\big\{ \|q_\JJ\Big(\sum_{i=0, \atop i+k\geq 0}^{s}
   a_{i+k,i}\otimes 1^{s-i}\Big)\|\big\},\, \big\|\sum_{i=0, \atop i+k\geq 0}^{r_0}
   (a_{i+k,i}\otimes 1^{r_0-i})\big\| \right\}.
$$
\end{theorem}
\begin{proof}
In view of Theorem \ref{norm formulas on 0-spectral subspace}, the function $\|\cdot \|_\JJ$ defined in Proposition \ref{propostion nie udowodnione}  satisfies the equality: $\|a^*a\|_\JJ=\|a\|_\JJ^2$,  for all $a\in \M^{(k)}_\TT(\KK)$, $k\in Z$. In particular, it is a $C^*$-seminorm on $\M^{(0)}_\TT(\KK)$ and hence (by the uniqueness of the $C^*$-norm) the norm of the element \eqref{element a form} is given by the same formula as the $\|\cdot\|_\JJ$-norm of the corresponding element of  $\M^{(k)}_\TT(\KK)$. The last part of the statement follows from the second part of Theorem \ref{norm formulas on 0-spectral subspace}.
\end{proof}
\begin{corollary}[Kernel of universal representation]\label{kernel of universal representation corollary}
The  universal representation  $\iota=\{\iota_{n,m}\}_{n,m\in \N}$ of $\KK$  in $\OO_\TT(\KK,\JJ)$ is faithful if and only if  $\JJ\subset (\ker\otimes 1)^{\bot}$, and in general we  have 
$$
a \in \ker \iota_{n,m} \,\,\, \Longleftrightarrow  \,\,\, \lim_{k\to \infty}\| a\otimes 1^k\|=0 \quad\textrm{ and }\quad a\otimes 1^k \in \JJ \,\,\,\textrm{ for all } k\in \N. 
$$
\end{corollary}
\begin{corollary}\label{Kernel of universal representation correspondence}
The universal  representation  of a $C^*$-correspondence $X$ in the $C^*$-algebra $\OO(J,X)$ is faithful if and only if
  $J\subset (\ker \phi )^{\bot}$.
\end{corollary}
The ideal described in Corollary \ref{kernel of universal representation corollary}  - the kernel of $\iota=\{\iota_{n,m}\}_{n,m\in \N}$, together with    condition \eqref{necessary condition for faithfulness7} presented below  play  essential role in investigation of  faithful  representations arising form right tensor representations, cf. Theorem \ref{Gauge invariance theorem for O_T(I,J)1}. In the context of $C^*$-correspondences condition \eqref{necessary condition for faithfulness7}   leads to the notion of a $T$-pair introduced in \cite{katsura2}, see Subsection \ref{Ideal structure  for relative Cuntz-Pimsner algebras subsection}.

\begin{proposition} \label{another one bites the proposition}
The ideal $\JJ$ is an ideal of coisometricity for the universal representation   of $\KK$ in $\OO_\TT(\KK,\JJ)$.
Moreover, we have
\begin{equation}\label{necessary condition for faithfulness7}
 \JJ(n,n)=\iota^{-1}_{n,n}\left( \sum_{j=1}^{k}\iota_{n+j,n+j}(\KK(n+j,n+j))\right),\,\,\, \textrm{ for all }n\in \N \textrm{ and }k>0.
\end{equation}
\end{proposition}
\begin{proof}
The form of the  ideal of coisometricity for $\{\iota_{n,m}\}_{n,m\in \N}$  follows from the first part of the proof of Theorem \ref{universality checking}. In particular, if $k>0$, then  $\JJ(n,n)\subset \iota^{-1}_{n,n}\left( \sum_{j=1}^{k}\iota_{n+j,n+j}(\KK(n+j,n+j))\right)$ and for the converse inclusion note that 
$$
\iota_{n,n}(a)=\sum_{j=1}^k \iota_{n+j,n+j}(b_j) \Longrightarrow \|\iota_{n,n}^\M(a)-\sum_{j=1}^k \iota_{n+j,n+j}^\M(b_j)\|_\JJ=0 \Longrightarrow a\in \JJ(n,n).
$$ 
\end{proof}
\begin{corollary} Every ideal $\JJ$ in $J(\KK)$ is an ideal of coisometricity for a certain right tensor representation of $\KK$.
\end{corollary}
\begin{corollary}\label{coisometricity corollary for correspondences}  Let $X$ be a  $C^*$-correspondence. Every ideal $J$ in $J(X)$ is an ideal of coisometricity  for a certain  representation $(\pi,t)$ of $X$.
\end{corollary}
 \begin{example} Let   $\varphi:A\to M(A_0)$ be a partial morphism. In view of Example \ref{endomorphism of unitals2}, for every ideal $J$ in  $\varphi^{-1}(A_0)$ there exists a covariant representation $(\pi,U,H)$ such that 
 $$
 J=\{a\in A: U^*U\pi(a)=\pi(a)\},
 $$
consult  with analogous results of \cite{kwa-leb1}.    
\end{example}
 \begin{example} For any graph  $E=(E^0,E^1,r,s)$ and every set of vertices $V\subset \{v\in E^0: |s^{-1}(v)|<\infty\}$  there exists a Cuntz-Krieger  $(E,V)$-family, which is not an $(E,V')$-family for any $V'$ bigger than $V$, cf. Example \ref{example conerning (E,V)-families}.
\end{example}

\section{$\OO_\TT(\KK,\JJ)$,  $\DR(\TT)$ and algebras with circle actions}
\label{Doplichery i cirlce}
The  \emph{Doplicher-Roberts algebra  $\DR(\TT)$ of a right tensor $C^*$-category} $\TT$    is defined  to be the completion of the algebraic direct sum 
$\bigoplus_{k\in\Z}\DR^{(k)}(\TT)$, where
$$
\DR^{(k)}(\TT):=\underrightarrow{\,\lim\,\,}\TT(r+k,r),
$$
in the unique $C^*$-norm for which the automorphic action defined by the grading is isometric \cite[p. 179] {dr}. It is evident that  this construction could be successfully applied to $C^*$-precategories (or even  to ideals in $C^*$-precategories). Thus we slightly extend existing nomenclature  and, for any right tensor $C^*$-precategory $\TT$, call the $C^*$-algebra $\DR(\TT)$ defined above  a \emph{Doplicher-Roberts algebra of the right tensor $C^*$-precategory} $\TT$.
\begin{proposition}[Relationship between $\OO_\TT(\KK,\JJ)$ and $\DR(\TT)$]\label{roberts and doplicher wear a T-shirt} 
Every algebra  $\DR(\TT)$ is an algebra of the type  $\OO_\TT(\KK,\JJ)$. Namely we have a natural isomorphism 
$$ 
\DR(\TT)\cong \OO_\TT(\TT,\TT).
$$ 
Conversely, for every  ideals $\KK$, $\JJ$ in $\TT$ such that $\JJ\subset J(\KK)$ the mappings
$$
\KK(n,m)\ni a\longmapsto \iota_{n,m}^\M(a)\in \KK_\JJ
$$
where $\{\iota_{n,m}^\M\}_{n,m\in\N}$ is the representation from Proposition \ref{*-algebra proposition} i), induce  an  isomorphism 
$$
 \OO_\TT(\KK,\JJ)\cong \DR(\KK_\JJ).
 $$ 
\end{proposition}
\begin{proof}
It is an immediate consequence of Theorem \ref{norm formulas on 0-spectral subspace} and the  definition of $\OO_\TT(\KK,\JJ)$ presented  below Theorem \ref{norm formulas on 0-spectral subspace}.
\end{proof}
\begin{corollary}
A $C^*$-algebra is a Doplicher-Roberts algebra associated with a $C^*$-category (as defined in \cite{dr}) if and only if it is an algebra of type   $\OO_\TT(\TT,\JJ)$ for a right tensor $C^*$-category $\TT$.
\end{corollary}
\begin{proof}
 If $\TT$ is a $C^*$-category, then the right tensor $C^*$-precategory $\TT_\JJ$ defined in Theorem \ref{norm formulas on 0-spectral subspace} (where $\KK=\TT$) is a $C^*$-category.
\end{proof}
%The foregoing statements indicate that  it is perhaps  reasonable to treat algebras of  the type  $  \OO(\TT,\JJ)\cong \DR(\TT_\JJ)$ as a \emph{relative Doplicher-Roberts algebras associated to $\TT$, relative to $\JJ$}.
 
It follows that the classes of $C^*$-algebras of type $\OO_\TT(\KK,\JJ)$ and type $\DR(\TT)$ coincide (in fact they are both equal to the class of all $C^*$-algebras). Obviously what distinguishes these algebras is an additional structure, particularly the associated gauge action. We  show that an arbitrary circle action may be viewed as a gauge action induced by our construction.
  \\
 Let  $\gamma:S^1 \to \Aut B$ be a circle action on a $C^*$-algebra $B$ and let $\{B_{n}\}_{n\in \Z}$ be the family of \emph{spectral subspaces}:
 $$
 B_{n}=\{b\in B: \gamma_z(b)=z^nb\,\, \textrm{ for }\,z\in S^1 \}. 
 $$ 
We set   $\TT=\{\TT(n,m)\}_{n,m\in\N}$ where  $\TT(r+k,r)=B_k$ and $\TT(r,r+k)=B_{-k}$ for all $r,k\in \N$.  Then  $\TT$ with  operations inherited from $B$ is a $C^*$-precategory, matricially presented in the following way 
\begin{equation}\label{trivial C - category}
\TT= \left( \begin{array}{c c c c  }  B_0      &         B_{-1}           &           B_{-2}         &           \cdots  
\\ B_{1}      &  B_0  &  B_{-1}  &   \cdots  \\  B_{2}      &  B_{1} & B_0  &   \cdots  \\  \vdots  &       \vdots        &       \vdots        & \ddots \end{array}\right).
\end{equation}
We  equip   $\TT$ with a right tensoring $\otimes 1$ which simply   slides the elements along the diagonals  in \eqref{trivial C - category}, that is $\TT(n,m)\ni a \stackrel{\otimes 1}{\longrightarrow} a \in \TT(n+1,m+1)$. In this way  $\TT$ becomes  a   right tensor $C^*$-precategory  and we have a natural gauge invariant isomorphism 
$$B\cong  \OO_\TT(\TT,\TT)=\DR(\TT).$$
 Therefore we get 
\begin{theorem}%\label{theorem about gauge actions}
Every $C^*$-algebra with a circle action is gauge invariantly isomorphic to an algebra of type $\DR(\TT)$ (and all the   more to an algebra of type $\OO_\TT(\KK,\JJ)$).
 \end{theorem}
 
The above statement  shows that taking into account gauge actions the class of relative Cuntz-Pimsner algebras  $\OO(J,X)$ is strictly smaller than the class of algebras of  type $\OO_\TT(\KK,\JJ)$. Indeed,  we have   
\begin{theorem}[Thm. 3.1 \cite{aee}]\label{theorem about gauge actions} A $C^*$-algebra with a circle action $\gamma$ is gauge invariantly isomorphic to a relative Cuntz-Pimsner algebra $\OO(J,X)$ if and only if the action $\gamma$ is semi-saturated.  \end{theorem}

 The right tensor $C^*$-precategory \eqref{trivial C - category} has an advantage that the universal representation  embeds the  spaces $\TT(k,0)$, $k\in \N$,  into $\OO_\TT(\TT,\TT)$ as spectral subspaces. We now present  necessary and sufficient conditions  assuring this in the  general case. We start with a statement that illustrates the forthcoming definition.
\begin{proposition}\label{complete}
Let $\KK$ be an ideal in a right tensor $C^*$-precategory $\TT$. The following
conditions are equivalent.
\begin{itemize}
\item[i)] There is a system of mappings   $\LL:\{\KK(n,m)\}_{n,m=1}^\infty \to \KK$ where  $\LL:\KK(n+1,m+1)\to \KK(n,m)$    satisfies the properties
$$
\LL(a)\otimes 1 =a \qquad (\textrm{right-inversion of }\otimes 1),
$$
$$
\LL(a (b\otimes 1))=\LL(a) b\qquad  (\textrm{transfer operator property}),
$$
for all $a\in \KK(n+1,m+1)$, and $b\in \KK(m,l)$  $n ,m, l\in \N$ (we do not explicitly impose any other  algebraic properties on $\LL$).
\item[ii)] The morphism $
\otimes 1: J(\KK)\cap (\ker\otimes 1)^\bot\longrightarrow \{\KK(n,m)\}_{n,m=1}^\infty$ 
is an epimorphism (then it is   automatically  an isomorphism).
\end{itemize}
If these conditions are satisfied, then  $\LL$ in item $i)$ is  determined uniquely: it is the inverse of the isomorphism $\otimes 1: J(\KK)\cap (\ker\otimes 1)^\bot\longrightarrow \{\KK(n,m)\}_{n,m=1}^\infty$.
In particular, $\LL:\{\KK(n,m)\}_{n,m=1}^\infty \to \KK$ is a homomorphism of $C^*$-precategories.
\end{proposition}
\begin{proof}
i) $\Longrightarrow$   ii). Let $a\in \KK(n+1,m+1)$. Since $\LL(a)\otimes 1 =a$ we get  $\LL(a)\in J(\KK)(n,m)$. For every  $b\in (\KK\cap \ker\otimes 1)(m,l)$ we have
$$
\LL(a) b=\LL(a (b\otimes 1))=\LL(0)=\LL(0 (0\otimes 1) )=\LL(0) \cdot 0=0.
$$
Hence $\LL(a)\in (\ker\otimes 1)^{\bot}$ and  it follows that the image of $\LL$ is contained in $ J(\KK)\cap (\ker\otimes 1)^\bot$. This together with the equality  $\LL(a)\otimes 1 =a$ implies that $\otimes 1: J(\KK)\cap (\ker\otimes 1)^\bot\longrightarrow \{\KK(n,m)\}_{n,m=1}^\infty$ is a surjection.
Furthermore,  since $\otimes 1$ is isometric on $(\ker\otimes 1)^\bot$, the image of  $\LL$  coincides with  $J(\KK)\cap (\ker\otimes 1)^\bot$ and $\LL$  coincides with the inverse of the isomorphism $\otimes 1: J(\KK)\cap (\ker\otimes 1)^\bot\longrightarrow \{\KK(n,m)\}_{n,m=1}^\infty$. 
\\
ii) $\Longrightarrow$   i). Define $\LL$ as  the inverse of  $\otimes 1: J(\KK)\cap (\ker\otimes 1)^\bot\rightarrow \{\KK(n,m)\}_{n,m=1}^\infty$. Relation  $\LL(a)\otimes 1 =a$ is trivially satisfied.  To show the
"transfer operator property" let $a\in \KK(n+1,m+1)$ and $b\in \KK(m,l)$, $n ,m, l\in \N$.  Then 
$$
\LL(a (b\otimes 1))\otimes 1= a (b\otimes 1)= (\LL(a)\otimes 1)( b\otimes 1) =(\LL(a) b)\otimes 1
$$
and since both  $\LL(a (b\otimes 1))$ and $\LL(a) b$ belong to $J(\KK)\cap (\ker\otimes 1)^\bot(n,l)$,
we get $\LL(a (b\otimes 1))=\LL(a) b$.
\end{proof}
 \begin{definition}
We  say that an ideal  $\KK$ in a right tensor $C^*$-precategory \emph{admits a transfer homomorphism} if it satisfies the equivalent conditions of Proposition \ref{complete}. 
\end{definition}
The role of the above introduced notion is explained by  the following  
\begin{theorem}\label{complete ideal in the $C^*$-category}
Let $\KK$ and $\JJ$ be ideals in a right tensor $C^*$-precategory $\TT$ such that $\JJ\subset J(\KK)$. The universal representation embeds each  space $\KK(k,0)$, $k\in\N$,    into  $\OO_\TT(\KK,\JJ)$ as the $k$-th spectral subspace for the associated gauge action $\gamma$ if and only if 
$$
\JJ=J(\KK)\cap (\ker\otimes 1)^\bot   
$$
and the ideal $\KK$ admits a transfer homomorphism.
\end{theorem} 
\begin{proof}
Assume that $
\JJ=J(\KK)\cap (\ker\otimes 1)^\bot
$ and $\KK$ admits a transfer homomorphism $\LL$. By Corollary \ref{kernel of universal representation corollary} the universal representation $\{\iota_{n,m}\}_{n,m\in\N}$ is faithful. To see that the image of the space  $\KK(k,0)$ is the  $k$-spectral subspace of $\OO_\TT(\KK,\JJ)$ let  $a\in \KK(r+k,r)$, $r>0$. Since $a=\LL(a)\otimes 1$  and $\LL(a)\in \JJ$, using Theorem \ref{proposition of norm formulas} we have $\iota_{r+k,r}(a)=\iota_{r-1+k,r-1}(\LL(a))$ 
 and following  in this way by induction one gets $\iota_{r+k,r}(a)=\iota_{k,0}(\LL^r(a))$.
\\
Now suppose that the  representation $\{\iota_{n,m}\}_{n,m\in\N}$ embeds  $\KK(k,0)$, $k\in\N$, as the  $k$-spectral subspace of $\OO_\TT(\KK,\JJ)$. Then $\JJ\subset J(\KK)\cap (\ker\otimes 1)^\bot$ by Corollary \ref{kernel of universal representation corollary}. To show the converse inclusion let $ a\in J(\KK)\cap (\ker\otimes 1)^\bot(n,n)$. Since $\iota_{n,n}(a)$ and $\iota_{n+1,n+1}(a\otimes 1)$ lie  in the spectral subspace $\iota_{0,0}(\KK(0,0))$  there exists $b\in \KK(0,0)$ such that $\iota_{0,0}(b)=\iota_{n,n}(a)-\iota_{n+1,n+1}(a\otimes 1)$.  Using Theorem \ref{proposition of norm formulas} we get 
$$
 b\otimes 1^n - a \in \JJ(n,n),\qquad b\otimes 1^{n+1}=0.
$$
It follows that   $b\otimes 1^n$ belongs both to $(\ker\otimes 1)^\bot(n,n)$ and $\ker\otimes 1(n,n)$. Therefore
$b\otimes 1^n=0$ and consequently $a\in \JJ(n,n)$. This  shows that $\JJ=J(\KK)\cap (\ker\otimes 1)^\bot$ and to prove that $
\otimes 1: J(\KK)\cap (\ker\otimes 1)^\bot\longrightarrow \{\KK(n,m)\}_{n,m=1}^\infty$ 
is  an epimorphism let $a\in \KK(n+1+k,n+1)$, $n,k\in \N$. Then there exists $b\in \KK(k,0)$ such that
$
\iota_{0,0}(b)=\iota_{n+1+k,n+1}(a)$  or equivalently 
$$
b\otimes 1^r\in \JJ(r+k,r),\quad  r=0,...,n,\qquad b\otimes 1^{n+1}= a. 
$$  
Hence
$
 a=c\otimes 1$  where $c=b\otimes 1^n \in \JJ(n+k,n)=J(\KK)\cap (\ker\otimes 1)^\bot(n+k,n).
$
\end{proof}
We  interpret this result  on the level of $C^*$-algebras associated with   $C^*$-cor\-respon\-den\-ces.
 \begin{proposition}\label{complete bimodule spectral subspace algebra8}
An ideal $\KK_X=\{\KK( X^{\otimes m},X^{\otimes n})\}_{n,m\in \N}$ in a right tensor $C^*$-precategory $\TT=\TT_X$ associated with a $C^*$-correspondence $X$ over a $C^*$-algebra $A$   admits a transfer homomorphism if and only if $X$ is a Hilbert $A$-bimodule.
\end{proposition} 
\begin{proof}
Apply Proposition \ref{equivalent characteristics of completeness} iii) and Proposition \ref{complete} ii).
\end{proof}
 Thus as a corollary to Theorem \ref{complete ideal in the $C^*$-category} we get  
\begin{theorem}[cf. Prop. 5.17 \cite{katsura}]\label{complete bimodule spectral subspace algebra2}
The universal representation of  a $C^*$-corres\-pon\-dence  $X$ over $A$  embeds  $A$  into a relative Cuntz-Pimsner algebra $\OO(J,X)$  as the fixed point algebra for the gauge action if and only if $X$ is a  Hilbert $A$-bimodule  and 
$$
J=(\ker\phi)^\bot \cap J(X).
$$
If this is the case, then $(\ker\phi)^\bot \cap J(X)=\overline{{_A\langle} X , X  \rangle}$ and   $\OO(J,X)$ is canonically isomorphic to the crossed-product $A\rtimes_X\Z$ by Hilbert bimodule $X$ in the sense of \cite{aee}. In particular, each space $X^{\otimes n}$, $n\in \N$, embeds into  $\OO(J,X)$ as the spectral subspace.
\end{theorem} 
 \begin{proof} If $J=(\ker\phi)^\bot \cap J(X)$ and  $X$ is a  Hilbert $A$-bimodule, then by  Theorem \ref{complete ideal in the $C^*$-category} algebra $A$ and spaces $X^{\otimes n}$, $n>0$, embed into $\OO(J,X)$ as spectral subspaces. Conversely, if $A$ embeds into $\OO(J,X)$ as the $0$-spectral space, then the argument from the proof of Theorem \ref{complete ideal in the $C^*$-category} shows that $J=(\ker\phi)^\bot \cap J(X)$ and  $X$  forms a Hilbert bimodule, cf. Proposition \ref{equivalent characteristics of completeness}. 
\end{proof}
By Examples \ref{partial morphism category1.0}, \ref{endomorphism of unitals2}, \ref{paragraf 1}, \ref{paragraf 2} we get
\begin{corollary}\label{complete2424666}
Let $\varphi:A\to M(A_0)$  be a partial morphism. Algebra $A$ embeds into the  crossed product $A\rtimes_{\varphi} \N$ as the $0$-spectral subspace if and only if $\varphi$ arises from a  partial automorphism (in our broader sense, see Example \ref{partial morphism category1.0}). If this is the case, then $A\rtimes_{\varphi} \N$ coincides with partial crossed-product as defined in \cite{exel1} (where instead of the  ideal $J$ one  puts hereditary subalgebra $A_0$).
\end{corollary}
\begin{corollary}\label{complete2424}
Let $\al:A\to A$ be an endomorphism of a unital  $C^*$-algebra $A$. Algebra $A$ embeds into the relative crossed product $C^*(\A,\al;J)$ as the $0$-spectral subspace if and only if $ \al$ admits a complete transfer operator $\LL$ %(equivalently   $\ker\al$  is unital and $\al(A)$ is hereditary subalgebra of $A$) 
and $J=\LL(A)=(\ker\al)^\bot$. If this is the case $C^*(\A,\al;J)$ coincides with the crossed product $A\rtimes_\al  \Z$ introduced in \cite{Ant-Bakht-Leb}, cf. \cite{kwa-leb3}.
\end{corollary}
By Examples \ref{example conerning (E,V)-families},  \ref{paragraf 3} we get
\begin{corollary}\label{complete24245345}
Let $E=(E^0,E^1,r,s)$ be a  directed graph and let $R(E)=\{v\in E^0:0<|s^{-1}(v)|<\infty\}$. Then  $A=C_0(E^0)$  embeds  into the relative graph algebra $C^*(E,V)$ as the $0$-spectral subspace if and only if $V=R(E)$ and $r,s$ are injective. In this case $C^*(E,V)$ coincides with both  the graph algebra $C^*(E)$ and the partial crossed-product
defined by the  partial homeomorphism $s \circ r^{-1}:r(E^1)\to s(E^1)$.
\end{corollary}

The above  statement could be generalized in a natural manner to topological graphs, cf. Example \ref{precategory of a directed graph ex1}. It reflects the fact that as a rule description of the core (the fixed point subalgebra) of a graph algebra in terms of the graph is a nontrivial step.
\section{Structure Theorem for $\OO_\TT(\KK,\JJ)$}
In this section we generalize  the main result of  \cite{fmr}. To this end we introduce  and discuss the notions of invariance and saturation for ideals in right tensor $C^*$-precategories.
\begin{definition}\label{invariant ideals and other things}
We  say that an ideal $\NN$ in a right tensor $C^*$-precategory $\TT$ is \emph{invariant} if 
$\NN \otimes 1 \subset \NN$. For such an ideal there is a \emph{quotient right tensor $C^*$-precategory}  $\TT/\NN$  defined by 
\begin{align*}
(\TT/\NN)(n,m)& :=  \TT(n,m)/\NN(n,m),
\\
\big(a+ \NN(n,m)\big)\otimes 1& := a\otimes 1 + \NN(n,m),
\end{align*}
where $
a\in \TT(n,m)$, $n,m\in \N$, cf. Proposition \ref{quotient C*-category}. More generally, if $\KK$ is an ideal in $\TT$   we say that an ideal $\NN$ is $\KK$-invariant if
$$
(\NN(n,m)\otimes 1) \KK(m+1,l)\subset  \NN(n+1,l), \qquad n,m,l\in\N,
$$
shortly  $(\NN\otimes 1)\KK\subset \NN$.
\end{definition}
%In view of Proposition \ref{another one bites the proposition} we  see that kernels of universal right tensor representations are invariant.  In fact we have the following
%Invariant ideals are "compatible" with right tensor representations.
\begin{notation}
If $\NN$ and $\KK$ are ideals in $\TT$, then  the image of $\KK$ in the quotient $C^*$-precategory $\TT/\NN$ shall be denoted by $\KK/\NN$ (obviously it may be  identified the quotient $C^*$-precategory $\KK/(\NN\cap \KK)$). If $\SSS$ is a sub-$C^*$-precategory of $\TT$ and $\KK$ is an ideal in $\TT$,  we denote by $\SSS+\KK$ a sub-$C^*$-precategory of $\TT$ where  $(\SSS+ \KK)(n,m):=\SSS(n,m)+ \KK(n,m)$, $n,m\in \N$. Similarly, if $\{\SSS_k\}_{k\in\N}$ is a family of ideals in $\TT$, we denote by $\sum_{k=0}^\infty \SSS_k$ the ideal in $\TT$ where  $(\sum_{k=0}^\infty \SSS_k )(n,m):=\clsp\{a\in\SSS_k(n,m):k\in \N\}$. If additionally $\SSS_0\subset \SSS_1\subset ...$, then we write $\underrightarrow{\,\lim\,\,}\SSS_k$ for $\sum_{k=0}^\infty \SSS_k$. If $a$ is a morphism in $\TT(n,m)$, then to say that $a$ is in $\SSS(n,m)$ we  briefly write $a\in  \SSS$.
\end{notation}
\begin{lemma}\label{the smallest invariant ideals}
Let  $\NN$ be a sub-$C^*$-precategory  in a right tensor $C^*$-precategory $\TT$. 
Then, for each $k\in \N$, the spaces 
$$
(\TT(\NN\otimes 1^k)\TT)(n,m):=\clsp\{a (b\otimes 1^k) c: a\in \TT(n,p+k), b\in \NN(p,r), c\in \TT(r+k,m)\}
$$
form an ideal in $\TT$.  In particular,  $\TT\NN\TT$ is the smallest ideal containing $\NN$ and 
$$
\widetilde{\NN}=\sum_{k=0}^\infty\TT(\NN\otimes 1^k)\TT
$$
is the smallest invariant ideal in $\TT$ containing $\NN$. 
\end{lemma}
\begin{proof}
One sees that $\TT(\NN\otimes 1^k)\TT$  is a  well defined $C^*$-precategory and since it is invariant under left and right multiplication, it is an ideal. The ideal  $\widetilde{\NN}=\sum_{k=0}^\infty\TT(\NN\otimes 1^k)\TT$ is a direct limit of the ascending sequence of ideals   $\sum_{k=0}^N\TT(\NN\otimes 1^k)\TT$ where 
$$
(\TT(\NN\otimes 1^k)\TT)\subset \TT(\NN\otimes 1^{k+1})\TT.
$$
Thus $\widetilde{\NN}$ is invariant. The minimality of $\widetilde{\NN}$  and $\TT\NN\TT$ is obvious. 
\end{proof}
\begin{remark}
With analogy to the nomenclature concerning graphs, cf. Example \ref{structure theorem for graphs} below,  one could call  invariant ideals in $\TT$ \emph{hereditary ideals} and then for any ideal $\NN$  one could call $\widetilde{\NN}=\sum_{k=0}^\infty\TT(\NN\otimes 1^k)\TT$ a \emph{hereditation} of $\NN$.
\end{remark}
\begin{proposition}\label{invariant ideals from right tensor representations1}
Let $\KK$ be  an ideal  in a right tensor $C^*$-precategory $\TT$. For  an ideal $\NN$ in $\KK$ the following conditions are equivalent: 
\begin{itemize}
\item[i)] $\NN$ is $\KK$-invariant,
\item[ii)] $\NN=\widetilde{\NN} \cap \KK$ for an invariant ideal $\widetilde{\NN}$ in $\TT$,  
\item[iii)]  $\NN=\ker\pi$ for a  right tensor representation $\pi$ of $\KK$.
\end{itemize}
In particular,  if $\pi$ is a  right tensor representation of $\KK$ and $\widetilde{\NN}$ is an invariant ideal in $\TT$  such that $\ker\pi=\widetilde{\NN} \cap \KK$, then  $\pi$ factors through to the faithful right tensor representation of the  ideal $\KK/\widetilde{\NN}\cong \KK/\ker\pi$ in the quotient right tensor $C^*$-precategory $\TT/\widetilde{\NN}$. 
\end{proposition}
\begin{proof}
i) $\Rightarrow$ ii). Let $
\widetilde{\NN}=\sum_{k=0}^\infty\TT(\NN\otimes 1^k)\TT
$ be the smallest invariant ideal in $\TT$ containing $\NN$. Note that by $\KK$-invariance of $\NN$ we have
$$
\TT(\NN\otimes 1^k)\TT \KK=\TT(\NN\otimes 1^k)\KK \subset \TT \NN =\NN.
$$ 
Therefore $\NN=\widetilde{\NN} \cap \KK$. 
\\
ii) $\Rightarrow$ iii). Consider the image $\KK/\widetilde{\NN}$ of  $\KK$ under the quotient homomorphism $q_{\widetilde{\NN}}:\TT\to \TT/\widetilde{\NN}$. Then one may define $\pi$ to be the composition of $q_{\widetilde{\NN}}$ with any faithful right tensor representation of $\KK/\widetilde{\NN}$ (such a representation exists by Corollary \ref{kernel of universal representation corollary}).
\\
iii) $\Rightarrow$ i). If $\NN=\ker\pi$ for a right tensor representation of $\KK$,  $a \in \NN(n,m)$ and $b \in \KK(m+1,l)$, then $\pi_{n+1,l}((a\otimes 1) b)=\pi_{n,m}(a)\pi_{m+1,l}(b)=0$, that is $(a\otimes 1) b\in \NN(n+1,l)$ and hence $\NN$ is $\KK$-invariant.
\end{proof}

Invariant ideals are closely related with kernels of right tensor representations. We reveal a similar relationship between saturated ideals and ideals of coisometricity for  right tensor representations.
\begin{definition}\label{definition of saturation for categories} Let $\NN$ and $\JJ$ be  ideals in a right tensor $C^*$-precategory $\TT$. We say that $\NN$ is \emph{$\JJ$-saturated} if  $\JJ\cap \otimes 1^{-1}(\NN)\subset \NN$. In general we put 
$$
\SSS_{\JJ}(\NN):=\underrightarrow{\,\lim\,\,}\SSS_k\qquad\,\, \textrm{ where }\,\, \SSS_0:=\NN\,\, \textrm{ and }\,\, \SSS_k:=\JJ\cap \otimes 1^{-1}(\SSS_{k-1})+\SSS_{k-1},\,\, k>0.
$$
Then $\SSS_{\JJ}(\NN)$ is the smallest $\JJ$-saturated ideal containing $\NN$ and we  shall call it  \emph{$\JJ$-saturation} of $\NN$. 
\end{definition}
The saturation works well with invariance.
\begin{lemma}\label{characterization of the saturation lemma}
If $\NN$ and  $\JJ$ are  ideals in a right tensor $C^*$-precategory $\TT$ and $\NN$ is invariant, then 
the  $\JJ$-saturation $S_{\JJ}(\NN)$ of $\NN$ is invariant and $a\in \TT(n,m)$ is in $\SSS_\JJ(\NN)$ if and only if
\begin{equation}\label{characterization of the saturation}
a\otimes 1^k \in \JJ +\NN,\,\,\,\textrm{ for all } k\in \N,\quad\textrm{ and }\quad \lim_{k\to \infty}\| q_\NN(a\otimes 1^k)\|=0.
\end{equation}
More generally, if $\JJ\subset J(\KK)$ and $\NN$ is $\KK$-invariant for an ideal $\KK$ in $\TT$, then 
$$
S_{\JJ}(\NN)=S_{\JJ}(\widetilde{\NN})\cap \KK
$$
for any invariant  ideal $\widetilde{\NN}$ such that $\NN=\widetilde{\NN}\cap \KK$.  In particular, $S_{\JJ}(\NN)$ is $\KK$-invariant.  
\end{lemma}
\begin{proof}
Invariance of $\SSS_\JJ(\NN)$ follow from inclusions  $\SSS_k\otimes 1\subset \SSS_{k-1}$, $k>0$. Passing to  $\TT/\NN$ one sees  that every element in $\TT(n,m)$ satisfying \eqref{characterization of the saturation} may be approximated (arbitrarily closely) by elements $a\in \TT(n,m)$ such that 
\begin{equation}\label{dense subset condition}
 a\otimes 1^k \in \JJ +\NN,\,\textrm{ for }k=0,..,N-1,\,\,\textrm{ and }\,\, a\otimes 1^N \in \NN, \,\, N\in \N. 
 \end{equation}
 We claim that $a$ satisfies  \eqref{dense subset condition} iff $a \in \SSS_N$ and therefore $a$ satisfies \eqref{characterization of the saturation} iff $a \in \SSS_\JJ(\NN)$. Indeed, let $a\in \TT(n,m)$ satisfy  \eqref{dense subset condition}. Then $a\otimes 1^N\in \SSS_0=\NN$ and 
 $$
 a\otimes 1^{N-1}\in \Big(\otimes 1^{-1}(\NN)\cap (\JJ+\NN)\Big)=(\otimes 1^{-1}(\NN)\cap  \JJ)+  \NN= \SSS_1.
 $$
Similarly, for  $k=1,..,N$,  one  sees that 
$$
a\otimes 1^{N-k}\in \SSS_k \quad\textrm{ and }\quad a\otimes 1^{N-k-1}\in \NN+\JJ 
$$
implies that $a\otimes 1^{N-k-1}\in \SSS_{k+1}$. Hence, by induction,  $a\in  \SSS_N$. Conversely, let 
$a\in \SSS_N(n,m)$. Then,  since $\SSS_{\JJ}(\NN)\subset \JJ+\NN$, we have $a\otimes 1^k \in \JJ +\NN$, $k=0,..,N-1$, and  since  $\SSS_0\otimes 1\subset \SSS_0=\NN$ and $\SSS_k\otimes 1\subset \SSS_{k-1}$, $k>0$, we get $a\otimes 1^N\in \SSS_0=\NN$.  
\\
To prove the second part of the statement we denote by $\{\SSS_{k}(\NN)\}_{k\in \N}$ and $\{\SSS_{k}(\widetilde{\NN})\}_{k\in \N}$ the increasing sequences of ideals whose direct limits  yield respectively $S_{\JJ}(\NN)$ and $S_{\JJ}(\widetilde{\NN})$. We shall show our claim by proving that
$$
\SSS_{k}(\NN)=\SSS_{k}(\widetilde{\NN})\cap \KK, \qquad \textrm{ for all }k\in \N.
$$
This relation trivially holds for $k=0$, and if we  suppose it holds for $k=N-1$,  then using inclusions $\JJ\subset J(\KK)\subset \KK$ we get
\begin{align*}
\SSS_{N}(\NN)&=\JJ\cap \otimes 1^{-1}(\SSS_{N-1}(\NN))+\SSS_{N-1}(\NN)
\\
&=\JJ\cap \otimes 1^{-1}(\SSS_{N-1}(\widetilde{\NN})\cap \KK)+\SSS_{N-1}(\widetilde{\NN})\cap \KK
\\
&=\JJ\cap \otimes 1^{-1}(\SSS_{N-1}(\widetilde{\NN}))+\SSS_{N-1}(\widetilde{\NN})\cap \KK
\\
&= \SSS_{N}(\widetilde{\NN})\cap \KK.
\end{align*}

\end{proof}
\begin{proposition}\label{invariant ideals from right tensor representations2}
Let $\JJ$ and  $\KK$ be ideals  in a right tensor $C^*$-precategory $\TT$ such that $\JJ\subset J(\KK)$.  For  an ideal $\NN$ in $\KK$ the following conditions are equivalent: 
\begin{itemize}
\item[i)] $\NN$ is $\KK$-invariant and $\JJ$-saturated,
\item[ii)] $\NN=\widetilde{\NN} \cap \KK$ for an invariant and $\JJ$-saturated ideal $\widetilde{\NN}$ in $\TT$,  
\item[iii)]  $\NN=\ker\pi$ for a  right tensor representation $\pi$ of $\KK$  coisometric on $\JJ$.
\end{itemize}
\end{proposition}
\begin{proof} i) $\Rightarrow$ ii). In view of Proposition \ref{invariant ideals from right tensor representations1} it suffices to note that if $\NN=\widetilde{\NN} \cap \KK$, then $\NN$ is $\JJ$-saturated iff  $\widetilde{\NN}$ is $\JJ$-saturated, which follows from the inclusion  $\JJ\subset J(\KK)$. 
\\
 ii) $\Rightarrow$ iii). Consider the same representation as in the prove of Proposition \ref{invariant ideals from right tensor representations1}, and note that   $\widetilde{\NN}$ is  $\JJ$-saturated iff the ideal $\JJ/\widetilde{\NN}$ in $\TT/\widetilde{\NN}$ is contained in the annihilator of the kernel of the quotient right tensoring on $\TT/\widetilde{\NN}$. Hence, by Corollary \ref{kernel of universal representation corollary}, the considered representation may be chosen to be coisometric on $\JJ$.
   \\
  iii) $\Rightarrow$ i).   $\KK$-invariance of  $\ker \pi$ follows from Proposition \ref{invariant ideals from right tensor representations1}  and  $\JJ$-saturation is a direct consequence of definitions of  $J(\KK)$ and $J$-saturation, Definitions \ref{nie wiem jak to nazwac definition}, \ref{definition of saturation for categories}. 
\end{proof}
Now we are in a position to prove the main result of this section.
\begin{theorem}[Structure Theorem]\label{structure theorem}
Let $\KK$, $\JJ$ and $\NN$ be ideals in a $C^*$-precategory $\TT$ such that $\JJ\subset J(\KK)$ and  $\NN\subset \KK$ is $\KK$-invariant. The subspace 
$$
\OO(\NN)=\clsp\{\iota_{n,m}(a): a\in  \NN(n,m),\,\, n,m\in \N\} \subset \OO_\TT(\KK,\JJ)
$$   
generated by the image of $\NN$ under the universal representation $\iota=\{\iota_{n,m}\}_{n,m\in \N}$ is an ideal in $\OO_\TT(\KK,\JJ)$ and  there are  natural isomorphisms
\begin{equation}\label{isomorphisms of structure theorem}
 \OO(\NN)\cong \OO_\TT(\NN,\JJ\cap\NN),\quad \OO_\TT(\KK,\JJ)/\OO(\NN)\cong \OO_{\TT/\widetilde{\NN}}(\KK/\widetilde{\NN},\JJ/\widetilde{\NN}).
 \end{equation}
 where $\widetilde{\NN}$ is an arbitrary invariant ideal in $\TT$ such that $\NN=\widetilde{\NN}\cap \KK$.
 Moreover, for the $\JJ$-saturation $\SSS_\JJ(\NN)$ of $\NN$ we have 
$$
\SSS_\JJ(\NN)=\iota^{-1}(\OO(\NN)).
$$
Hence   $\OO(\SSS_\JJ(\NN))=\OO(\NN)$ and in the right hand sides of  \eqref{isomorphisms of structure theorem} the ideals  $\NN$ and  $\widetilde{\NN}$ may be replaced by their $\JJ$-saturations  $\SSS_\JJ(\NN)$ and $\SSS_\JJ(\widetilde{\NN})$. 
 \end{theorem} 
 \begin{proof}
 Using operations \eqref{add1} -- \eqref{star1} we see that 
$$
\M_\TT(\NN)=\spane\{\iota_{n,m}(a): a\in \NN(n,m),\,\, n,m\in \N\}
$$
is a two-sided ideal in the algebra $\M_\TT(\KK)$ defined in Subsection \ref{2}. Therefore  $\OO(\NN)=\clsp\{\iota_{n,m}(a): a\in \NN(n,m),\,\, n,m\in \N\}$ is an ideal in  $\OO_\TT(\KK,\JJ)$. To prove that $\OO(\NN)\cong \OO_\TT(\NN,\JJ\cap\NN)$ we show that  the seminorms $\|\cdot \|_\JJ$ and $\|\cdot \|_{\JJ\cap \NN}$  give the same quotients of  $\M_\TT(\NN)$, cf.  Proposition \ref{propostion nie udowodnione}. Let $a\in \M_\TT(\NN)$. Obviously,  $\|a\|_{\JJ\cap \NN}=0$  implies $\|a\|_{\JJ}=0$. Conversely, we may assume that  $a= \sum_{s=0}^{r}
   \iota_{s,s}(a_{s,s})$,  $a_{s,s}\in  \NN(s,s)$, $s=0,...,r$, $r\in \N$,  and then the condition $\|a\|_{\JJ}=0$ is equivalent to
  $$
  \sum_{j=0}^{s}
   a_{j,j}\otimes 1^{s-j}\in \JJ(s,s), \quad s=1,...,r-1,\qquad \lim_{r\to \infty}\sum_{j=0}^{r}
   a_{j,j}\otimes 1^{r-j}=0.
  $$
By $\KK$-invariance of $\NN$ we get  
$$
  \sum_{j=0}^{s}
   a_{j,j}\otimes 1^{s-j}\in \JJ\cap \NN (s,s), \quad s=1,...,r-1,\qquad \lim_{r\to \infty}\sum_{j=0}^{r}
   a_{j,j}\otimes 1^{r-j}=0
  $$  
  which is equivalent to $\|a\|_{\JJ\cap \NN}=0$. Hence $\OO(\NN)\cong \OO( \NN,\JJ\cap\NN)$.
  \\
Let $\widetilde{\NN}$ be an invariant ideal in $\TT$ such that $\NN=\widetilde{\NN}\cap \KK$.  To construct the isomorphism $\OO_\TT(\KK,\JJ)/\OO(\NN)\cong \OO_{\TT/\widetilde{\NN}}(\KK/\widetilde{\NN},\JJ/\widetilde{\NN})$ we consider the right tensor representation  $\pi=\{\pi_{n,m}\}_{n,m\in \N}$    of
$\KK$ in $\OO_{\TT/\widetilde{\NN}}(\KK/\widetilde{\NN},\JJ/\widetilde{\NN})$ given by  $
\pi_{n,m}=\iota_{n,m}\circ q_{\widetilde{\NN}}$, $m,n\in \N$, where $\iota=\{\iota_{n,m}\}_{n,m\in \N}$ denotes the universal representation of $\KK/\widetilde{\NN}$ in $\OO_{\TT/\widetilde{\NN}}(\KK/\widetilde{\NN},\JJ/\widetilde{\NN})$. Since for $a\in \JJ(n,m)$ we have
\begin{align*}
q_{\widetilde{\NN}}(a) \in (\JJ/\widetilde{\NN})(n,m) \,\,\,&\Longleftrightarrow\,\,\, \iota_{n,m}(q_{\widetilde{\NN}}(a))=\iota_{n+1,m+1}(q_{\widetilde{\NN}}(a)\otimes 1) 
\\ & \Longleftrightarrow \,\,\, \pi_{n,m}(a)=\pi_{n+1,m+1}(a\otimes 1),
\end{align*}
it follows that $\pi$ is coisometric on $\JJ$ and thereby induces a homomorphism $\Psi_\pi$ from $\OO_\TT(\KK,\JJ)$ onto $\OO_{\TT/\widetilde{\NN}}(\KK/\widetilde{\NN},\JJ/\widetilde{\NN})$. Plainly,   $\Psi_\pi$ is zero on $\OO(\NN)$ and hence  it factors through to an epimorphism $\Psi_\pi:\OO_\TT(\KK,\JJ)/\OO(\NN) \longrightarrow \OO_{\TT/\widetilde{\NN}}(\KK/\widetilde{\NN},\JJ/\widetilde{\NN})$.
One  proves injectivity of $\Psi_\pi$ by constructing its inverse. Indeed, the formula
$$
\omega_{n,m}(a +\widetilde{\NN}(n,m))=q(\iota_{n,m}(a)),\qquad a\in (\KK/\widetilde{\NN})(n,m), \,\,\, n,m\in \N,
$$
where $q:\OO_\TT(\KK,\JJ)\to \OO_\TT(\KK,\JJ)/\OO(\NN)$ is the quotient map, defines a  right tensor representation $\omega$ of $\KK/\widetilde{\NN}$ in $\OO_\TT(\KK,\JJ)/\OO(\NN)$ which is coisometric on $\JJ/\widetilde{\NN}$. Thus $\omega$ integrates to  a homomorphism from $\OO_{\TT/\widetilde{\NN}}(\KK/\widetilde{\NN},\JJ/\widetilde{\NN})$ to $\OO_\TT(\KK,\JJ)/\OO(\NN)$ which is   inverse to $\Psi_\pi$.% Hence the first part of the  proof is complete.\\
\\
Furthermore,  we have $\ker  \pi_{n,m}=\iota_{n,m}^{-1}(\OO(\NN))$ (where $\pi$ is  the representation defined  above) and thus by Corollary \ref{kernel of universal representation corollary}, $a\in \KK(n,m)$ is  in $\iota_{n,m}^{-1}(\OO(\NN))$ if and only if
$$
 a\otimes 1^k \in \JJ +\widetilde{\NN},\,\,\,\textrm{ for all }\,\, k\in \N,\quad\textrm{ and }\quad \lim_{k\to \infty}\| q_{\widetilde{\NN}}(a\otimes 1^k)\|=0. 
$$
In view of Lemma  \ref{characterization of the saturation lemma} this proves the second part of the theorem.
\end{proof}
The Structure Theorem has a number of important consequences. 
\\
Firstly, $\OO(\NN)$ is a  \emph{gauge invariant ideal} in $\OO_\TT(\KK,\JJ)$, i.e. it is globally invariant under the associated gauge action, and if $\NN$ is $\JJ$-saturated it  is uniquely determined by $\OO(\NN)$. Hence,   denoting by $\Lat_\JJ(\KK)$ the lattice of $\KK$-invariant  and $\JJ$-saturated ideals in $\KK$, and by $\Lat(\OO_\TT(\KK,\JJ))$ the lattice of gauge invariant ideals in   $\OO_\TT(\KK,\JJ)$  we get the natural embedding
$$
\Lat_\JJ(\KK) \hookrightarrow  \Lat(\OO_\TT(\KK,\JJ)).
$$
In general this embedding is not an isomorphism.  However, we  show in Theorem \ref{theorem on lattice structure2} that in certain important cases  we have $\Lat_\JJ(\KK) \cong \Lat(\OO_\TT(\KK,\JJ))$.\par
Secondly,  the zero ideal $\NN=\{0\}$ is invariant and hence dividing $\OO_\TT(\KK,\JJ)$ by the ideal generated by the $\JJ$-saturation $\SSS_\JJ(\{0\})$ of $\{0\}$ we  reduce the relations defining $\OO_\TT(\KK,\JJ)$ without  affecting the algebra $\OO_\TT(\KK,\JJ)$ itself. We formalize this remark as follows.
\begin{definition}\label{reduction ideal definition}
Let $\JJ$ be an ideal in a right tensor $C^*$-precategory $\TT$. We denote the $\JJ$-saturation $\SSS_\JJ(\{0\})$  of the zero ideal by $\RR_\JJ$ and call it a \emph{reduction ideal} associated with  $\JJ$. 
\end{definition}
 \begin{theorem}[Reduction of relations] \label{reduction of relations theorem}
  Let $\KK$ be an ideal in a right tensor $C^*$-category $\TT$ and let $\JJ$ be an ideal in $\JJ(\KK)$. % and $\RR=\SSS_\JJ(\{0\})$ the $\JJ$ saturation of $\{0\}$. 
  Putting
  $$
  \TT_\RR:= \TT/ \RR_\JJ,  \qquad \KK_\RR:=\KK/ \RR_\JJ, \qquad \JJ_\RR:=\JJ/ \RR_\JJ,
  $$
we get  the "reduced" right tensor $C^*$-category $\TT_\RR$, with right tensoring  $\otimes 1_\RR$. The "reduced" ideals $\KK_\RR$, $\JJ_\RR$ are such that   $\JJ_\RR\subset \JJ(\KK_\RR) \cap (\ker\otimes 1_\RR)^\bot$ and there is  a natural isomorphism
  $$
  \OO(\KK,\JJ)\cong \OO(\KK_\RR,\JJ_\RR).
  $$
    \end{theorem}
  \begin{proof}
Clear by   Theorem \ref{structure theorem}. \end{proof}
    \begin{remark}\label{remark reduction}

We may often reduce even  more relations in $\TT$  
in the sense that for  any invariant ideal $\RR'$ in $\TT$
such that $\RR'=\RR_\JJ\cap \KK$ we have
  $$
  \OO_\TT(\KK,\JJ)\cong \OO_{\TT_{\RR'}}(\KK_{\RR'},\JJ_{\RR'}),
  $$
 where  $\KK_{\RR'}:=\KK/ \RR'$ and $\JJ_{\RR'}:=\JJ/\RR'$ coincide with $\KK_\RR$ and  $\JJ_\RR$, but in general $\TT_{\RR'}=\TT/ {\RR'}$ is  "smaller" than $\TT_\RR= \TT/ \RR_\JJ$. 
 \end{remark}

 \subsection{Structure theorem for relative Cuntz-Pimsner algebras}\label{structure theorem for relative Cuntz-Pimsner algebras subsection}
As an application of   Theorem \ref{structure theorem} we  get a version of Structure Theorem for relative Cuntz-Pimsner algebras which improves  \cite[Thm 3.1]{fmr}. For that purpose we
fix a $C^*$-correspondence $X$ over a $C^*$-algebra $A$, and establish the relationships between the relevant ideals in $\KK_X$ and $A$.
\begin{proposition}\label{proposition on quotients of correspondences}
Each of the relations  
\begin{equation}\label{relations on quotients of correspondences}
I=\NN(0,0), \qquad \NN= \KK_X(I)
\end{equation}
establish a one-to-one correspondence  between $X$-invariant ideals $I$ in $A$ and $\KK_X$-invariant ideals $\NN$ in $\KK_X$. Moreover,  we have a natural isomorphism of $C^*$-precategories 
$$
\KK_X/\NN\cong \KK_{X/XI}=\{\KK\left((X/XI)^{\otimes n},(X/XI)^{\otimes m}\right)\}_{n,m\in \N}.
$$
\end{proposition}
\begin{proof}
For an $X$-invariant ideal $I$   in $A$, $\widetilde{\NN}:= \TT_X(I)$ is  invariant and hence $\NN=\widetilde{\NN}\cap \KK=\KK_X(I)$ is $\KK_X$-invariant. Conversely, if $\NN$ is an ideal in $\KK_X$, then by Proposition \ref{Ideals in TT(X)}, $\NN=\KK_X(I)$ where $I:= \NN(0,0)$. If additionally $\NN$ is $\KK_X$-invariant, then  $I$ is $X$-invariant because for $a\in I=\KK_X(I)(0,0)$ and $b=\Theta_{x,y}\in \KK_X(1,1)$ we have
$$
(a\otimes 1)b=\varphi(a)\Theta_{x,y}=\Theta_{\varphi(a)x,y}\in \KK_X(I)(1,1),
$$
that is $\varphi(a)x \in XI$.  The remaining part of the assertion  follows from Lemma \ref{lemma to quotients} and Corollary \ref{corollary to quotients}.
\end{proof}
In the  case when $J=J(X)\cap (\ker\phi)^\bot$, the notion we are about to introduce   coincides with the property called $X$-saturation in \cite[Def. 6.1]{mt} and negative invariance  in \cite[Def. 4.8]{katsura2}, cf. also \cite[Def. 4.14]{katsura2}.
 \begin{definition}\label{definition of saturation for correspondences} 
Let $I$ and $J$ be ideals in $A$. We  say that  $I$ is \emph{$J$-saturated} if 
   $J \cap \phi^{-1}(\LL(XI))\subset I$, or equivalently if
    $$
a \in J\,\,\,\textrm{ and }\,\,\, \varphi(a)X\subset XI\,\, \Longrightarrow \,\, a\in I.
$$ 
In general we put  $
S_{J}(I):=\underrightarrow{\,\lim\,\,}S_k$ where $S_0:=I$ and $S_k:= J\cap \phi ^{-1}(S_{k-1})+S_{k-1}$, $k>0$.
Then $S_{J}(I)$ is the smallest $J$-saturated ideal containing $I$ which we shall call  \emph{$J$-saturation} of $I$.  
  \end{definition}
  \begin{proposition}\label{proposition on saturation of correspondences}
Let $J$ and $I$  be ideals in $A$. The $\KK_X(J)$-saturation  of $\KK_X(I)$ coincides with $\KK_X(S_J(I))$. In particular,  $I$ is $J$-saturated if and only if $\KK_X(I)$ is $\KK_X(J)$-saturated. 
%Relations \eqref{relations on quotients of correspondences} establish a one-to-one correspondence  between $X$-invariant $J$-saturated ideals  in $A$ and ideals of the form $\KK_X\cap \NN$ in $\KK_X$ where   $\NN$ is an invariant $\KK_X(J)$-saturated ideal in $\TT_X$.
\end{proposition}
\begin{proof}
It suffices to  check that under notation of Definitions  \ref{definition of saturation for categories}, \ref{definition of saturation for correspondences} we have $\KK_X(S_n)=\SSS_n$, $n\in \N$, which  is straightforward.  \end{proof}
In view of the  above statement, Corollary \ref{coisometricicity of representations} and  Proposition \ref{invariant ideals from right tensor representations2} we get
  \begin{proposition}\label{invariant and saturated  kernels vs   representations}
Let  $J\subset J(X)$. If  $J$ is an ideal of coisometricity for a  representation  $(\pi,t)$ of  $X$, then  $\ker \pi$ is $J$-saturated and $X$-invariant. Conversely, every  $J$-saturated and $X$-invariant ideal in $A$  is the kernel $\ker \pi$ for a certain representation $(\pi,t)$ of  $X$ whose ideal of coisometricity is $J$. 
\end{proposition}

To state the Structure Theorem for $C^*$-correspondences in its full force we  introduce  algebras that generalize relative Cuntz-Pimsner algebras. 
% Using the notation of subsection \ref{paragraph for TT(X)1}  we could simply write $\OO(J,X)=O_{\TT}(\KK_X,\KK_X(J))$, and in fact  it is useful to consider the following generalization of relative Cuntz-Pimsner algebras $\OO(J,X)$, cf. Theorem \ref{structure theorem2}.
 \begin{definition}\label{definition of strange algebras} 
If  $I$ and $J$  are ideals in $A$ such that $J\subset J(XI)$  we put 
 $$
 \OO_X(I,J):=\OO_{\TT_X}(\KK_X(I),\KK_X(J)).
 $$
 This algebra is well defined because $\KK_X(J) \subset J(\KK_X(I))$ iff $J\subset J(XI)$. In particular,  we have
 $\OO(J,X)=\OO_{X}(A,J)$.
  \end{definition}
  We already noted that ideal structures of  $\KK_{XI}=\{\KK( XI^{\otimes n},XI^{\otimes m})\}_{n,m\in \N}$ and $\KK_X(I)=\{\KK( X^{\otimes n},X^{\otimes m}I)\}_{n,m\in \N}$ are isomorphic, see  Remark \ref{remark on difference between  II_{XI} and II_X(I)}. By the following results, this  could be interpreted as    these $C^*$-precategories are   "Morita equivalent" with the  "equivalence"  established by $L=\{\KK( XI^{\otimes m},X^{\otimes n})\}_{n,m\in \N}$.
 \begin{lemma}\label{lemma on morita equivalence}
Let $I$ be an arbitrary ideal in $A$. Then   
$$
L=\{\KK( XI^{\otimes m},X^{\otimes n})\}_{n,m\in \N}
$$
is a left ideal in a $C^*$-precategory $\TT_X$, and 
 the following relations hold
$$
L^*L=\KK_{XI},\qquad LL^*=\KK_X(I), %\qquad \KK(XI)L\subset \KK(XI),\qquad L\KK(I)\subset \KK(I)
$$
where $L^*L(n,m)=\clsp\{x^*y: x \in L(k,n), y\in L(k,m)\}$ and  $LL^*(n,m)=\clsp\{xy^*: x \in L(n,k), y\in L(m,k)\}$, $n,m\in \N$.
\end{lemma}
\begin{proof} Clearly, $\KK_{XI}\subset L$ and hence $\KK_{XI}=\KK_{XI}^*\KK_{XI}\subset L^*L$. The opposite inclusion follows from the fact that for $x,u \in X^{\otimes k}$, $y\in XI^{\otimes n}$ and $v\in XI^{\otimes m}$
$$
(\Theta_{u,v})^*\Theta_{x,y}=\Theta_{v,y \langle x,u\rangle}\in \KK( XI^{\otimes n},XI^{\otimes m}).
$$
 Thus $L^*L=\KK_{XI}$.  To show that $LL^*=\KK_X(I)$ it suffices to apply  Proposition \ref{Ideals in TT(X)}, since both $LL^*$ and $\KK_X(I)$ are ideals in $\KK_X$ such that  $LL^*(0,0)=\KK_X(I)(0,0)=I$.
\end{proof}
The above "Morita equivalence"  of $C^*$-precategories yields Morita equivalence  of $C^*$-algebras.  
\begin{theorem}[Morita equivalence of $\OO_X(I,I\cap J)$ and $\OO_{XI}(I, J\cap I)$] \label{structure theorem by morita}
Let $I$ be an $X$-invariant ideal in $A$ and let $J$ be an ideal in $J(X)$. The algebras 
$$
 \OO_X(I,J\cap I),\qquad \OO_{XI}(I, J\cap I)= \OO(J\cap I, XI)
$$
may be naturally considered as subalgebras of $\OO_{\TT_X}(\TT_X,\KK_X(J\cap I))$. They coincide whenever $\varphi(I)X=XI$, and in general they are Morita equivalent with an equivalence  established via 
the subspace
$$
\LL =\clsp\{\iota_{n,m}(a): a\in L(n,m)\}\subset \OO_{\TT_X}(\TT_X,\KK_X(J\cap I)).
$$
\end{theorem}
\begin{proof} The first part of the statement is evident, cf. Definition \ref{definition of strange algebras}. To  see the second part consider the  $^*$-algebra $\M_\TT(\TT)$ defined in Subsection \ref{2} for  $\TT=\TT_X$. For any sub-$C^*$-precategory $\SSS$ of $\TT$ we put 
$$
\M_\TT(\SSS)=\spane\{\iota_{n,m}^\M(a): a\in \SSS(n,m)\}.
$$
Then $\M_\TT(\KK(I))$ is a two-sided ideal in $\M_\TT(\TT)$, and $\M_\TT(L)$ is a left ideal in $\M_\TT(\TT)$ such that 
$$
\M_\TT(L)\star \M_\TT(L)^*=\M_\TT(\KK(I)).
$$
Indeed, since $L\subset \KK(I)$ we have $\M_\TT(L)\star \M_\TT(L)^*\subset \M_\TT(\KK(I))$, and the opposite  inclusion follows from Lemma \ref{lemma on morita equivalence}.   Similarly, $\M_\TT(\KK(XI))$ is a $^*$-subalgebra of $\M_\TT(\TT)$ such that  
$$
\M_\TT(L)^*\star \M_\TT(L)=\M_\TT(\KK(XI)).
$$ 
Indeed,    since $\KK(XI)\subset L$ we have  $\M_\TT(\KK(I))\subset\M_\TT(L)^*\star \M_\TT(L)$. To get the opposite inclusion notice that for any $y \in  XI^{\otimes n}$, $v  \in  XI^{\otimes m}$, $u, x_1 \in  X^{\otimes l}$, $x_2 \in  X^{\otimes k}$, by $X$-invariance of $I$, the operator
$$
\big((\Theta_{u,v})^*\otimes 1^k\big)\,\, \Theta_{x_1\otimes x_2,\, y}=\Theta_{v\otimes \varphi(\langle u, x_1 \rangle )x_2,\, y}
$$
is in $\KK( XI^{\otimes n},XI^{\otimes (m+k)})$.
\\
Now it suffices to take the quotients with respect to seminorm $\|\cdot\|_{\KK_X(J\cap I)}$ defined  in Proposition \ref{propostion nie udowodnione}
and  apply the enveloping procedure.
\end{proof}
   We are in a position to formulate   the main result of this subsection. 
\begin{theorem}[Structure Theorem for $C^*$-correspondences]\label{structure theorem2}
Suppose $J$ is an ideal in $J(X)$  and let $\OO(I)$ denote the ideal in $\OO(J,X)$ generated by the image of an $X$-invariant ideal $I$ under the universal representation.  There are  natural isomorphisms
\begin{equation}\label{structure isomorphisms}
  \OO(I)\cong \OO_X(I,J\cap I),\qquad   \OO(J,X)/\OO(I)\cong \OO(J/I,X/XI).
\end{equation}
 Algebras  $\OO(I)$  and $\OO(J\cap I,XI)$ are Morita equivalent and if  $\phi(I)X=XI$, then  simply  $ \OO(I)\cong\OO_X(I,J\cap I)=\OO(J\cap I,XI)$. Moreover, the $J$-saturation $S_J(I)$ of $I$ is $X$-invariant and  
 $$
 S_J(I)=\iota_{0,0}^{-1}(\OO(I)).
 $$
Thus $\OO(I)=\OO(S_J(I))$  and in the right hand sides of  \eqref{structure isomorphisms} the $X$-invariant ideal $I$ may be replaced by the $X$-invariant and $J$-saturated ideal $S_{J}(I)$.
 \end{theorem} 
 \begin{proof}
 One  sees that $\OO(I)=\clsp\{\iota_{n,m}(a): a\in \KK_X(I)(n,m), n,m\in\N\}$ and  hence the assertion follows from   Theorems \ref{structure theorem},  \ref{structure theorem by morita} and Propositions \ref{proposition on quotients of correspondences}, \ref{proposition on saturation of correspondences}.
 \end{proof}
 \begin{remark}\label{structure theorem2 remark}
For an $X$-invariant ideal $I$ in $A$  the  subspace $\phi(I)X$ of $X$ may be considered as a $C^*$-correspondence over $I$. The argument from \cite[Prop. 9.3]{katsura} shows that $\OO(J\cap I,\phi(I)X)$  is Morita equivalent to  $\OO(I)$. Thus we have three $C^*$-algebras with natural embeddings
$$\OO(J\cap I,\phi(I)X)\subset \OO(J\cap I,XI)\subset \OO_X(I,J\cap I)\cong \OO(I)$$
which are all Morita equivalent, and if $\phi(I)X=XI$ then they are actually equal. One of the advantages of our approach is that we have been able to identify $\OO(I)$ precisely as  $\OO_X(I,J\cap I)=\OO_{\TT_X}(\KK_X(I),\KK_X(J\cap I))$ (not only up to Morita equivalence, as it is in \cite{fmr},\cite{katsura}).  
 \end{remark}
The reduction procedure  from  Theorem \ref{reduction of relations theorem} is in full consistency with
 reduction of relations in $C^*$-correspondences presented in \cite{kwa-leb1}. 
 \begin{definition}\label{reduction procedure for X definition} 
 Let $J$ be an ideal in $A$.
The $J$-saturation $S_J(\{0\})$  of the zero ideal will be denoted by $R_J$ and called a \emph{reduction ideal} associated to  $J$. 
\end{definition}
 \begin{theorem}[Reduction of  $C^*$-correspondences] \label{reduction procedure for X}  
 Let   $J$ be an ideal in $J(X)$. The reduction ideal  $R_J$ is $X$-invariant, and putting 
$$
X_R:=X/XR_J,\qquad A_R:=A/R_J, \qquad J_R:=J/R_J,
$$
we get the "reduced"  $C^*$-correspondence $X_R$ over $A_R$    such that $J_R \subset J(X_R)\cap \ker(\phi_R)^\bot$ where $\phi_R$ is the left action on $X_R$, and
 $$
\OO(J,X) 
\cong \OO(J_R,X_R).
$$
   \end{theorem}
\begin{proof}
Clear by Theorem \ref{structure theorem2}. 
\end{proof}\begin{example}\label{structure theorem for partial morphisms}
Let $X=X_\varphi$ be  the $C^*$-correspondence associated with a partial morphism $\varphi:A\to M(A_0)$.  One readily checks that an ideal $I$ in $A$ is  $X$-invariant if and only if  $\varphi(I)A_0\subset I$. For an $X$-invariant ideal $I$ we have a restricted partial morphism $\varphi_I:I \to M(A_0\cap I)$ and a quotient partial morphism  $\varphi^I:A/I \to M(A_0/I)$ where 
$$
\varphi_I=\varphi|_I,\qquad  \varphi^I(a+I)(a_0+I):=\varphi(a)a_0+I,\qquad  a\in A,\,\, a_0\in A_0.
$$
Both $\varphi_I$ and  $\varphi^I$ are well defined as $A_0\cap I$ is a hereditary subalgebra of $I$, $A_0/I$ is a hereditary subalgebra of $A/I$,  $\varphi_I(I)A_0\cap I=A_0\cap I$  and $\varphi^I(A/I)(A_0/I)=A_0/I$.  Furthermore, we may  naturally identify (as $C^*$-correspondences) $X_{\varphi^I}$ with $X/XI$  and $X_{\varphi_I}$ with $\phi(I)X$.  Using our general definition of relative crossed products (Definition \ref{definition of crossed product by partial morphisms}),  by Theorem   \ref{structure theorem2} and Remark \ref{structure theorem2 remark} we get
 $$
 C^*(\varphi;J)/\OO(I)\cong C^*(\varphi^I;J/I),
 $$
 where $\OO(I)$ is Morita equivalent to $C^*(\varphi_I;J\cap I)$, and $J$ is an ideal in $\varphi^{-1}(A_0)$. If we denote by $J_{\infty}$  the ideal in $J$ consisting of those elements $a\in J$  for which the iterates $\varphi^n(a)$, $n\in \N$, make sense and belong to $J$, i.e.
$$
a \in J_\infty \,\,\, \Longleftrightarrow \,\,\, a\in J,\,\,\, \varphi(a)\in J\cap A_0, \,\,\, \varphi^2(a) \in J\cap A_0,\, \, \, ...,\, \,\, \varphi^n(a) \in J\cap A_0,\,\, ...\,\,\,,
$$
then the reduction ideal $R:=R_J$  assumes the following form 
$$
R=\overline{\{a\in J_{\infty}: \exists_{n\in \N} \,\, \varphi^n(a)=0\}},
$$
cf. \cite{kwa-leb1}. For the  quotient partial morphism $\varphi^{R}:A/R\to M(A_0/R)$ we have
$$
C^*(\varphi,J)\cong C^*(\varphi^{R},J/R).
$$
Hence $\varphi^{R}$ may be viewed as  a natural reduction of the partial morphism $\varphi$ relative to the ideal $J$.
\end{example}
\begin{example}\label{structure theorem for graphs}
Suppose that $X=X_{E}$ is the $C^*$-correspondence of a graph $E$, $I$ is an ideal in $A=C_0(E^0)$ and $F\subset E^0$  is a complement of the hull of $I$. Then  
$$
XI=\clsp\{\delta_e:r(e)\in F\},\qquad \phi(I)X=\clsp\{\delta_e:s(e)\in F\}.
$$
It follows that  $I$ is $X$-invariant if and only if $V$  is \emph{hereditary}, that is if $s(e)\in F\Longrightarrow r(e)\in F$, for all $e\in E^1$, cf. \cite{bprs}, \cite{bhrs}, \cite{mt}. When $F$ is hereditary, then (slightly abusing notation) we may consider $E\setminus F:=(E^0\setminus F,r^{-1}(E^0\setminus F),r,s)$ and $F:=(F,s^{-1}(F),r,s)$ as subgraphs of $E$. In this event $X/XI$ is canonically isomorphic to $X_{E\setminus F}$ and $\phi(I)X$ is canonically isomorphic to $X_{F}$, cf. \cite[Ex. 2.4]{fmr}.  We recall that $J(X)= \clsp\{\delta_v:|s^{-1}(v)|< +\infty\}$. Hence Theorem   \ref{structure theorem2} imply that for any hereditary subset $F\subset E^0$ and any  $V \subset \{v \in E^0: |s^{-1}(v)|<\infty\}$  there is an ideal $\OO(I)$ in the relative graph algebra $\OO(E,V)$ such that 
 $$
 \OO(E,V)/\OO(I)\cong \OO(E\setminus F, V\setminus F)
 $$
 and $\OO(I)$ is Morita equivalent to the relative graph algebra  $C^*(F,F\cap V)$.
 We shall say  that a subset $F\subset E^0$ is \emph{$V$-saturated} if every vertex in $V\subset E^0$ which feeds into $F$ and only $F$ is in $F$:
$$
v\in V \,\,\textrm{ and }\,\, \{r(e):s(e)=v\}\subset F \,\,\,\Longrightarrow\,\,\, v \in F.
$$
By a  $V$-\emph{saturation} of a set $F$ we  mean  the smallest saturated subset $S_V(F)$ of $E^0$ containing $F$. In the case $V= \{v \in E^0: 0<|s^{-1}(v)|<\infty\}$ these notions coincide with the ones (without prefix $V$)  defined in \cite{bhrs}, \cite{bprs}. If $J=\clsp\{\delta_v:v\in V\}$, then the reduction ideal $R:=R_J$ is spanned by the point masses of  the $V$-saturation  $S_V(\emptyset)$ of the empty set ($S_V(\emptyset)$ consists of vertices in $V$ that form paths  leading to  sinks). As a consequence we have 
 $$
 \OO(E,V)\cong \OO\left(E\setminus S_V(\emptyset), V\setminus S_V(\emptyset)\right).
 $$
 Hence $E\setminus S_V(\emptyset)$ may be considered as a  natural  reduction of the  graph $E$ relative to the set $V$. 
\end{example}

\section{Ideal structure of $\OO_\TT(\KK,\JJ)$}
In this section we prove gauge invariance theorem for $\OO_\TT(\KK,\JJ)$ and describe the lattice of gauge-invariant ideals in  $\OO_\TT(\KK,\JJ)$. These results generalize the corresponding statements for relative Cuntz-Pimsner obtained in \cite{katsura}, \cite{fmr}, \cite{katsura2} (see the relevant discussions  in Subsection \ref{Ideal structure  for relative Cuntz-Pimsner algebras subsection} and Section \ref{Representations of algebras associated with X}).\par
%\subsection{Gauge-invariance theorem and  gauge-invariant ideals}
Let us fix  ideals $\KK$, $\JJ$   in a $C^*$-precategory $\TT$ such that $\JJ\subset J(\KK)$ and  let $\KK_\JJ$ be the  $C^*$-precategory defined in Theorem \ref{norm formulas on 0-spectral subspace}. Representation $\{\iota_{n,m}^\M\}_{n,m\in\N}$  from Proposition \ref{*-algebra proposition} yields the injective homomorphism of $\KK$ into $\KK_\JJ$:  
\begin{equation}\label{identifying ideals formula}
\KK(n,m)\ni a\longmapsto \iota_{n,m}^\M(a)\in \KK_\JJ(n,m).
\end{equation}
We use it to 
 adopt the identifications
$$
\KK\subset \KK_\JJ,\qquad \DR(\KK_\JJ)=\OO_\TT(\KK,\JJ).
$$
cf. Proposition \ref{roberts and doplicher wear a T-shirt}.
\begin{proposition}\label{extension of right tensor representation is this proposition} 
We have a one-to-one correspondence between the right tensor representations $\pi=\{\pi_{n,m}\}_{n,m\in \N}$ of $\KK$ coisometric on  $\JJ$ and right tensor representations $\widetilde{\pi}=\{\widetilde{\pi}_{n,m}\}_{n,m\in \N}$ of $\KK_\JJ$ coisometric on $\KK_\JJ$. This correspondence is given by \begin{equation}\label{extension of right tensor representation is this} 
\widetilde{\pi}_{r,r+k}\Big(\sum_{j=0, \atop j+k\geq 0}^{r}
   \iota_{j+k,j}^\M(a_{j+k,j})\Big)=\sum_{j=0, \atop j+k\geq 0}^{r}
   \pi_{j+k,j}(a_{j+k,j}).
\end{equation}
%where  $\{\iota_{n,m}\}_{n,m\in \N}$ is a right tensor representation of  $\TT$ in  $\M_\TT(\KK)$, see Theorem \ref{*-algebra proposition}.
\end{proposition} 
\begin{proof}
In view of the definition of $\KK_\JJ$ the assertion  may be   verified directly. One may also deduce it from  Theorem  \ref{universality checking} and  Propositions  \ref{representations  extensions  proposition}, \ref{roberts and doplicher wear a T-shirt}.
\end{proof}
It is well known that if two $C^*$-algebras $A$, $B$ admit circle actions, then a $^*$-homo\-mor\-phism $h:A\to B$ that maps faithfully spectral subspaces of $A$ onto the corresponding spectral subspaces of  $B$ is an isomorphism if and only if it is gauge invariant. Thus the following statement can be thought of as a (stronger) version of what is usually meant by a gauge invariance theorem.
\begin{definition}
A representation $\pi=\{\pi_{n,m}\}_{n,m\in \N}$ of an ideal $\KK$ in a right tensor $C^*$-precategory is said to admit a  gauge action if for every $z\in S^1$ relations 
$$
\beta_z(\pi_{n,m}(a))=z^{n-m}\pi_{n,m}(a),\qquad  a\in \KK(n,m), \, n,m\in \N,
$$
give rise to a well defined $^*$-homomorphism $\beta_z: C^*(\pi) \to C^*(\pi)$ where $C^*(\pi)$ stands for the $C^*$-algebra generated by the spaces $\pi_{n,m}(\KK(n,m))$, $n,m\in \N$.
\end{definition}

\begin{theorem}[Gauge-invariant uniqueness  for $\OO_\TT(\KK,\JJ)$]\label{Gauge invariance theorem for O_T(I,J)1}
Let $\pi=\{\pi_{n,m}\}_{n,m\in \N}$ be a right tensor representation  of $\KK$ coisometric on  $\JJ\subset J(\KK)$ and let $\RR_\JJ$ be the reduction ideal  associated with $\JJ$ (Definition \ref{reduction ideal definition}). The following conditions are equivalent
\begin{itemize}
\item[i)] $\pi$ integrates to a representation  faithful on the  $0$-spectral subspace of $\OO_\TT(\KK,\JJ)$
\item[ii)] $\pi$ integrates to a representation faithful on  all of  the spectral subspaces of $\OO_\TT(\KK,\JJ)$
\item[iii)]  $\ker\widetilde{\pi}$ coincides with $\RR_\JJ$  embedded into $\KK_\JJ$ via \eqref{identifying ideals formula}  
\item[iv)] $\ker\pi=\RR_\JJ$ and for  $n\in \N$  %(the reduction ideal associated to $\JJ$
\begin{equation}\label{necessary condition for faithfulness}
 \JJ(n,n)=\pi^{-1}_{n,n}\left( \sum_{j=1}^{k}\pi_{n+j,n+j}(\KK(n+j,n+j))\right),\qquad\textrm{ for all }k>0 .
\end{equation}
\end{itemize}
In particular, $\pi$ integrates to the  faithful representation of $\OO_\TT(\KK,\JJ)$ if and only if $\pi$ admits a gauge action,
and one of the equivalent conditions i)-iv) holds.\end{theorem}
\begin{proof} Equivalence i) $\Longleftrightarrow$ ii) follows from $C^*$-equality and algebraic relations between spectral subspaces, cf. \cite{exel1}. The representations $\Psi_\pi$ and  $\Psi_{\widetilde{\pi}}$ of $\OO_\TT(\KK,\JJ)$, arising from $\pi$ and $\widetilde{\pi}$, coincide and hence we put $\Psi=\Psi_\pi=\Psi_{\widetilde{\pi}}$.  From the direct limit construction of $\DR(\KK_\JJ)$ we see that  $\Psi$ is  faithful on spectral subspaces iff it is faithful on spaces $\phi_{n,m}(\KK_\JJ(n,m))$, $n,m\in \N$, where $\phi=\{\phi_{n,m}\}_{n,m\in \N}$ is a universal representation of $\KK_\JJ$ in $\DR(\KK_\JJ)$. The latter requirement is equivalent to the equality $\ker\widetilde{\pi}=\ker \phi$ and we claim that $\ker \phi$ coincidences with $\RR_\JJ$.  Indeed,   $
\phi_{r,r}\Big(\sum_{j=0}^{r}
   \iota_{j,j}^\M(a_{j})\Big) =0$ if and only if   $\lim_{k\to\infty} \|\sum_{j=0}^{r} \iota_{j,j}^\M(a_{j})\|_{r+k,r+k}^\JJ=0$, and this is equivalent to the conditions   
$$
\sum_{j=0}^{s}
   a_{j}\otimes 1^{s-j} \in \JJ,\textrm{ for }s=0,..,r-1,\,\,\textrm{ and }\,\,\, b\otimes 1 ^k\in \JJ, \, k\in \N,\, \,\,\,  \lim_{k\to\infty}\|b\otimes 1^k\| =0,
$$
where $b=\sum_{j=0}^{r}
   a_{j}\otimes 1^{r-j}$. In other words,  $\sum_{j=0}^{r}
   \iota_{j,j}^\M(a_{j})$ as an element of $\KK_\JJ(r,r)$ could be identified with $\iota_{r,r}^\M(b)$ and 
  $b \in \RR_\JJ=\SSS_{\JJ}(\{0\})$, 
cf. \eqref{characterization of the saturation}. This proves our claim and hence ii) $\Longleftrightarrow$ iii).
\\
Implication ii) $\Longrightarrow$ iv)  follows from Corollary \ref{kernel of universal representation corollary} and   Proposition \ref{another one bites the proposition}. To obtain iv) $\Longrightarrow$ iii) suppose that $ \sum_{j=0}^{r}
   \iota_{j,j}^\M(a_{j}) \in \ker\widetilde{\pi}_{r,r}$, that is  $\sum_{j=0}^{r}  \pi_{j,j}(a_{j})=0$. %, see \eqref{extension of right tensor representation is this}. 
   Then  
   $
 \pi_{0,0}(a_0)=-\sum_{j=1}^{r}  \pi_{j,j}(a_{j})$ implies (by \eqref{necessary condition for faithfulness}) that $a_0 \in \JJ(0,0).
   $
Hence  $\pi_{0,0}(a_0)=\pi_{1,1}(a_0\otimes 1)$ and consequently $
 \pi_{1,1}(a_0\otimes 1 + a_1)=-\sum_{j=2}^{r}  \pi_{j,j}(a_{j})$ implies that $a_0\otimes 1 + a_1 \in \JJ(0,0)
   $. Proceeding in this way one gets 
   $$
   \sum_{j=0}^{s}
   a_{j}\otimes 1^{s-j} \in \JJ,\,\textrm{ for } \,s=0,..,r-1,\,\,\,\textrm{ and }\,\,\, \pi_{r,r}( \sum_{j=0}^{r}
   a_{j}\otimes 1^{r-j})=0,
   $$
    that is $\sum_{j=0}^{r}
   \iota_{j,j}^\M(a_{j})$ as an element of $\KK_{\JJ}(r,r)$ coincides with 
$  \iota_{r,r}^\M(\sum_{j=0}^{r}
  a_{j}\otimes 1^{r-j})$ and $\sum_{j=0}^{r}
  a_{j}\otimes 1^{r-j}\in \ker\pi=\RR_{\JJ}$.
Hence $\ker\widetilde{\pi}\subset \RR_{\JJ}$. The inclusion $\RR_{\JJ} \subset \ker\widetilde{\pi}$ is obvious.
\end{proof}
\begin{corollary}[Gauge invariance theorem for $\DR(\TT)$] A right tensor representation  $\pi=\{\pi_{n,m}\}_{n,m\in \N}$ of $\TT$ coisometric on $\TT$ integrates to the faithful representation of $\DR(\TT)$ if and only if $\pi$ admits a gauge action and $\ker\pi=\RR_\TT$.
\end{corollary}
\begin{remark}
The condition  \eqref{necessary condition for faithfulness} implies that $\JJ$ is an ideal of coisometricity for $\pi$. In many natural situations the ideal of coisometricity automatically satisfies  \eqref{necessary condition for faithfulness}, cf.  Corollary \ref{takie tam corollary to be or not to be} below.  On the level  of Cuntz-Pimsner algebras, the condition  \eqref{necessary condition for faithfulness}   for $k=1$ implies it   for all $k>0$, cf. Remark \ref{remark for katsura}, and in fact  it reduces to equality \eqref{necessary condition for faithfulness 1}, the role of which  was discovered by Katsura  \cite{katsura} and led him to the notion of a $T$-pair  (see Proposition \ref{proposition on ideals and lattice structure} and Definition \ref{definition of T-pairs} below). 
\end{remark}
Continuing the discussion undertaken below Theorem \ref{structure theorem},  we give a complete description of gauge invariant ideal structure of $\OO_\TT(\KK,\JJ)$.
\begin{theorem}[Lattice structure of gauge invariant ideals in $\OO_\TT(\KK,\JJ)$]\label{theorem on lattice structure}
We have  a one-to-one correspondences  between the following objects  
\begin{itemize}
\item[i)] kernels of right tensor representations $\widetilde{\pi}$ of $\KK_\JJ$ coisometric on $\KK_\JJ$
\item[ii)] invariant and $\KK_\JJ$-saturated ideals $\NN_\JJ$ in $\KK_\JJ$, that is ideals satisfying
$$\NN_\JJ\otimes_\JJ 1\subset \NN_\JJ\quad \textrm{  and } \quad (\otimes_\JJ 1)^{-1}(\NN_\JJ)\subset \NN_\JJ
$$
\item[iii)] gauge invariant ideals $P$ in $\OO_\TT(\KK,\JJ)=\DR(\KK_\JJ) $
\end{itemize}
 These correspondences  preserve inclusions  and are given by 
 $$
\NN_\JJ= \ker\widetilde{\pi},\quad \ker\Psi_{\widetilde{\pi}}=\OO(\NN_\JJ)=P,\quad \NN_\JJ(n,m)=\iota_{n,m}^{-1}(P),\quad n,m\in \N,
 $$
  where $\iota=\{\iota_{n,m}\}_{n,m\in \N}$ denotes the universal right tensor representation of $\KK_\JJ$ in $\DR(\KK_\JJ)$ and $\OO(\NN_\JJ)$ is the closed linear  span of the image of $\NN_\JJ$ in $\DR(\KK_\JJ)$. In particular,  
  $$
  P \cong \DR(\NN_\JJ),\qquad \qquad  \textrm P=\clsp\{\iota_{n,m}(a): a\in \NN_\JJ(n,m)\}\subset \DR(\KK_\JJ),
 $$
 and we have the lattice isomorphism $\Lat_{\KK_\JJ}(\KK_\JJ)\cong \Lat(\OO_\TT(\KK,\JJ))$.
 \end{theorem} 
\begin{proof}
If   $\widetilde{\pi}$ is  a right tensor representation of $\KK_\JJ$, then by Proposition \ref{invariant ideals from right tensor representations2} $\ker\widetilde{\pi}$ is an invariant and $\KK_\JJ$-saturated ideal in $\KK_\JJ$. Conversely, if $\NN_\JJ$ is an invariant and $\KK_\JJ$-saturated ideal in $\KK_\JJ$, then by Theorem \ref{structure theorem} we have a gauge invariant homomorphism $\Psi:\DR(\KK_\JJ)\to \DR(\KK_\JJ/\NN_\JJ)$ whose kernel is  $\OO(\NN_\JJ)$. Hence disintegrating $\Psi$ we get a  right tensor representation $\widetilde{\pi}$ of $\KK_\JJ$ coisometric on $\KK_\JJ$   such that  $\NN_\JJ= \ker\widetilde{\pi}$. This proves the correspondence between the objects in i) and ii).
\\
Let  $P$ be an arbitrary gauge invariant ideal in $\OO_\TT(\KK,\JJ)=\DR(\KK_\JJ)$.   The spaces $\NN_\JJ(n,m):=\iota_{n,m}^{-1}(P)$, $n,m\in \N$, form  an ideal $\NN_\JJ$ in $\KK_\JJ$. By  definition of the algebra  $\DR(\KK_\JJ)$ we have $\iota_{n,m}(a)=\iota_{n+1,m+1}(a\otimes 1)$, for all $a\in \KK_\JJ(n,m)$. Thus  $\NN_\JJ$ is both   invariant and $\KK_\JJ$-saturated. Since  $\OO(\NN_\JJ)\subset P$   the identity map factors through to  a    surjection 
$$
\Psi:\DR(\KK_\JJ)/\OO(\NN_\JJ)\longrightarrow \DR(\KK_\JJ)/P.
$$
As  the ideals $\OO(\NN_\JJ)$ and $P$ are gauge invariant the gauge action on $\DR(\KK_\JJ)$ factors through to gauge actions on  $\DR(\KK_\JJ)/\OO(\NN_\JJ)$ and $\DR(\KK_\JJ) /P$. The epimorphism $\Psi$ intertwines these actions.  In view of Theorem \ref{structure theorem} we may identify  $\DR(\KK_\JJ)/\OO(\NN_\JJ)$ with $\DR(\KK_\JJ/\NN_\JJ)$ and then  $\Psi\Big(\iota_{n,m}(a+\NN_\JJ(n,m))\Big)=\iota_{n,m}(a)+P$, for $a\in \KK_\JJ(n,m)$.  Thus we see that   $\Psi$ is injective on every space $\iota_{n,m}\Big(\KK_{\JJ}(n,m)+\NN_\JJ(n,m)\Big)$. Therefore  $\Psi$ is injective on spectral subspaces of $\DR(\KK_\JJ/\NN_\JJ)=\DR(\KK_\JJ)/\OO(\NN_\JJ)$, and since it is gauge invariant,  it is an isomorphism. Hence $\OO(\NN_\JJ)=P$. This together with Theorem  \ref{structure theorem} proves the correspondence between the objects in ii) and iii).
\end{proof}
\begin{corollary}[Lattice structure of gauge invariant ideals in $\DR(\TT)$]  There is a lattice isomorphism $\Lat_{\TT}(\TT)\cong \Lat(\DR(\TT))$  between the gauge invariant ideals in $\DR(\TT)$ and invariant $\TT$-saturated ideals in $\TT$. Moreover an ideal in $\Lat(\DR(\TT))$ corresponding to an ideal  $\NN$ in  $\Lat_{\TT}(\TT)$ is isomorphic to $\DR(\NN)$.
\end{corollary}
It follows from Propositions \ref{invariant ideals from right tensor representations1}, \ref{extension of right tensor representation is this proposition}  that 
for every invariant  $\KK_\JJ$-saturated  ideal $\NN_\JJ$ in $\KK_\JJ$ the intersection $
\NN:=\KK\cap \NN_\JJ$ (we use  our  identification $\KK\subset \KK_\JJ$ given by \eqref{identifying ideals formula}) is  a $\KK$-invariant and $\JJ$-saturated ideal $\NN$ in $\KK$. Conversely, any $\KK$-invariant and $\JJ$-saturated ideal $\NN$ in $\KK$   give rise to the invariant and $\KK_\JJ$-saturated  ideal $\NN_\JJ[\NN]$ in $\KK_\JJ$ where 
$$
\NN_\JJ[\NN](r,r):=\left\{ \sum_{j=0}^{r}
   \iota_{j,j}^\M(a_{j}): a_{j} \in  \NN(j,j)\right\}.
$$
For an invariant and $\KK_\JJ$-saturated  ideal $\NN_\JJ$ in $\KK_\JJ$ we have  $\NN_\JJ[\KK\cap \NN_\JJ]\subset \NN_\JJ$ but in general  $\NN_\JJ[\KK\cap \NN_\JJ]\neq \NN_\JJ$. This is exactly the reason why the embedding  $\Lat_\JJ(\KK) \hookrightarrow \Lat_{\KK_\JJ}(\KK_\JJ)\cong \Lat(\OO_\TT(\KK,\JJ))$ in general  fails to be an isomorphism. 
 We now introduce conditions under which the aforementioned obstacle vanish.
\begin{theorem}\label{theorem on lattice structure2}
If one of the conditions hold 
\begin{itemize}
\item[i)] $\KK\subset \JJ+\ker\otimes 1$,
\item[ii)] $\KK$ admits a transfer homomorphism and $\JJ=J(\KK)\cap (\ker(\otimes 1))^\bot$,
\end{itemize}
then every ideal in $\Lat_{\KK_\JJ}(\KK_\JJ)$ have the form $\NN_\JJ[\NN]$ where $\NN \in \Lat_\JJ(\KK)$. As a consequence we get the lattice isomorphisms
$$
\Lat_\JJ(\KK)\cong \Lat_{\KK_\JJ}(\KK_\JJ)\cong \Lat(\OO_\TT(\KK,\JJ)).
$$
Moreover, every gauge invariant ideal $P$ in $\OO_\TT(\KK,\JJ)$ is generated by the image of a $\KK$-invariant  ideal $\NN$  in $\KK$ and  $P\cong \OO_\TT(\NN,\JJ\cap \NN)$. 
\end{theorem}
\begin{proof}
Let $\NN_\JJ \in \Lat_{\KK_\JJ}(\KK_\JJ)$ and let $\NN:=\KK\cap \NN_\JJ$.   By Theorem   \ref{diagonal of ideals} (and definition of $\KK_\JJ$) it suffices to show   that an element $a$ in $\NN_\JJ(n,n)$ represented by  $\sum_{j=0}^{n}
   \iota_{j,j}^\M(a_{j})$ where $a_j\in \KK(j,j)$, may also be represented by $\sum_{j=0}^{n}
   \iota_{j,j}^\M(b_{j})$ where $b_j\in \NN(j,j)$, $j=0,..,n$. For abbreviation  we shall write equality between $a$ and its  representatives.
   \par
    i). Suppose that $\KK\subset \JJ+\ker\otimes 1$. Then $a_0=b_0+j_0$ where $b_0\in (\ker\otimes 1)(0,0)$ and $j_0\in \JJ(0,0)$. Plainly $a=\iota_{0,0}^\M(b_0) + \iota_{1,1}^\M (j_0\otimes 1+a_1)+ \sum_{j=2}^{n}
   \iota_{j,j}^\M(a_{j})$. Since  $j_0\otimes 1+a_1$ is in $\KK\subset \JJ+\ker\otimes 1$ we have  $j_0\otimes 1+a_1 =b_1+j_1$ where $b_1\in (\ker\otimes 1)(1,1)$, $j_1\in \JJ(1,1)$ and then $a=\iota_{0,0}^\M(b_0) + \iota_{1,1}^\M (b_1)+ \iota_{1,1}^\M (j_1\otimes 1+a_2)+\sum_{j=3}^{n}
   \iota_{j,j}^\M(a_{j})$. Proceeding in this way one gets  
   $$
   a=\sum_{j=0}^{n}
  \iota_{j,j}^\M(b_j)\,\,\, \textrm{ where }\,\, b_j\in (\ker\otimes 1)(j,j),\,\, j=0,...,n-1,\,\, b_n\in \KK(n,n).
   $$ 
   We denote by $\{\mu_\lambda^{(j)}\}_{\lambda}$ an approximate unit in $\KK(j,j)$, $j=0,...,n$. Using the form of multiplication in $\KK_\JJ$ we obtain 
   $$
   \iota_{n,n}^\M(b_n)=\lim_{\lambda}\,\, \iota_{n,n}^\M(b_n) \star  \iota_{n,n}^\M(\mu_\lambda^{(n)})=\lim_{\lambda}\,\, a \star  \iota_{n,n}^\M(\mu_\lambda^{(n)})\in \NN_{\JJ}(n,n)
,$$
 that is  
        $b_n\in \NN(n,n)$. Similar computations show that $\iota_{n-1,n-1}^\M(b_{n-1})$ equals to
        $$
        \lim_{\lambda}\,\, \left(a \star  \iota_{n-1,n-1}^\M(\mu_\lambda^{(n-1)})- \iota_{n,n}^\M(b_n (\mu_\lambda^{(n-1)}\otimes 1))\right)\in \NN_{\JJ}(n,n).
        $$
Hence, since $(\otimes_\JJ 1)^{-1}(\NN_\JJ)\subset \NN_\JJ$, we get $\iota_{n-1,n-1}^\M(b_{n-1})\in \NN_{\JJ}(n-1,n-1)$, that is  
        $b_{n-1}\in  \NN(n-1,n-1)$. Proceeding in this way one obtains $b_{j}\in \NN(j,j)$, $j=0,..,n$.
       \par
    ii). Suppose that  $\JJ=J(\KK)\cap (\ker(\otimes 1))^\bot$ and $\KK$ admits a transfer homomorphism  $\LL$. 
One  sees that
$$
a=  \iota_{0,0}^\M(b_0) \,\,\, \textrm{ where }\,\,\, b_0=\sum_{j=0}^{n}\LL^{n-j} (a_j)\in \KK(0,0).
$$ 
Since $(\otimes_\JJ 1)^{-1}(\NN_\JJ)\subset \NN_\JJ$ one gets $\iota_{0,0}^\M(b_0)\in \NN_{\JJ}(0,0)$, that is  
        $b_{0}\in  \NN(0,0)$ and the proof is complete. 
        \end{proof}
        \begin{corollary}\label{takie tam corollary to be or not to be}
        If  $\KK\subset \JJ+\ker\otimes 1$, then a  right tensor representation $\pi$  of $\KK$ coisometric on  $\JJ$ integrates to the faithful representation of $\OO_\TT(\KK,\JJ)$ if and only if  $\pi$ is gauge invariant, $\JJ$ is the ideal of coisometricity for $\pi$ and  $\ker\pi=\RR_\JJ$.
        \end{corollary}
   
    \begin{corollary}\label{takie tam corollary to be or not to be2}
        If $\KK$ admits a transfer homomorphism and $\JJ=J(\KK)\cap (\ker(\otimes 1))^\bot$, then a   right tensor representation $\pi$  of $\KK$ coisometric on  $\JJ$ integrates to the faithful representation of $\OO_\TT(\KK,\JJ)$ if and only if  $\pi$ is faithful and gauge invariant.
        \end{corollary}
 \subsection{Ideal structure  of relative Cuntz-Pimsner algebras}\label{Ideal structure  for relative Cuntz-Pimsner algebras subsection}
Now we show   how a complete description of the  gauge invariant ideal structure of  relative Cuntz-Pimsner algebras $\OO(J,X)$  can be deduced  from Theorems  \ref{theorem on lattice structure}, \ref{theorem on lattice structure2}. Hopefully it will shed more light on the    results of   \cite{katsura2}.
 We start with a  useful statement  which follows from Theorem \ref{theorem on lattice structure2}  and    contains, as particular cases, \cite[Cor. 8.7, Thm. 10.6]{katsura2}, \cite[Thm. 6.4]{mt}.

     \begin{theorem}\label{sometimes all ideals are of simple form} Let $X$ be a $C^*$-correspondence over $A$ and let $J$ be an ideal in $J(X)$. If one of the conditions hold 
\begin{itemize}
\item[i)] $A=J+\ker\phi$,
\item[ii)] $X$ is a Hilbert bimodule  and $J=J(X)\cap (\ker \phi)^\bot$,
\end{itemize}
then we  have an  isomorphism between the lattice of gauge invariant ideals in $\OO(J,X)$ and lattice of $X$-invariant  $J$-saturated ideals in $A$: 
$$
\OO(J,X)\,\, \triangleright \,\,P \,\,\,\mapsto \,\,\, \iota_{0,0}^{-1}(P)\,\, \triangleleft \,\, A.
$$
 In particular,  every gauge invariant ideal $P$ in $\OO(J,X)$ is generated by the image of an $X$-invariant ideal $I$ in $A$ and then  $P$ is isomorphic to $\OO_X(I,J\cap I)$ and Morita equivalent to $\OO(J\cap I,XI)$.  If the condition ii) holds we actually have $P\cong\OO(J\cap I,XI)$. 
\end{theorem}
\begin{proof}
If  $A=J+\ker\phi$, then  on the level of the $C^*$-precategory $\TT_X$ we have $\KK_X\subset \KK_X(J)+\ker\otimes 1$.  If $X$ is a Hilbert bimodule and $J=J(X)\cap (\ker \phi)^\bot$, then  by Proposition \ref{complete bimodule spectral subspace algebra8}, $\KK_X$ admits a transfer homomorphism and $\KK_X(J)=J(\KK_X)\cap  (\ker \otimes 1)^\bot$. Moreover, for every  $X$-invariant $J$-saturated ideal $I$ in $A$  we have $\phi(I)X=XI$, see \cite[Prop. 10.2]{katsura2}. Thus it suffices to apply Theorems  \ref{theorem on lattice structure2}, \ref{structure theorem2}.
\end{proof}
A general result (without additional assumptions) requires  description of  invariant and $\KK_\JJ$-saturated ideals in  the right tensor $C^*$-precategory $\KK_\JJ$ where 
$$\KK:=\KK_X \,\,\textrm{ and }\,\,  \JJ:=\KK_X(J).$$
Equivalently, by Theorem  \ref{theorem on lattice structure}, we need to describe all the kernels of right tensor representations $\widetilde{\pi}$ of $\KK_\JJ$ coisometric on $\KK_\JJ$.  We recall that each such representation $\widetilde{\pi}$ is uniquely determined by a  representation  $(\pi,t)$ of $X$ coisometric on $J$, see Proposition \ref{proposition:rho} and Proposition \ref{extension of right tensor representation is this proposition}.

\begin{proposition}\label{proposition on ideals and lattice structure0}
Let  $(\pi,t)$ be  a faithful representation of $X$ coisometric on $J$ and let $I'$ be the ideal of coisometricity for $(\pi,t)$ (then automatically $J\subset I'\subset (\ker\phi)^\bot\cap J(X)$). 
The kernel of the corresponding representation $\widetilde{\pi}$ of $\KK_\JJ$  is uniquely determined by $I'$. Namely, $\ker \widetilde{\pi}(n,n)$, $n\in \N$, consists of elements that can be represented by $\sum_{k=0}^n \iota_{k,k}^\M(a_k)$,  $a_k\in \KK(X^{\otimes k})$, such that 
\begin{equation}\label{relations to be proved}
\sum_{k=0}^j a_k\otimes 1^{j-k}\in \KK(X^{\otimes j}I'),\,\, j=0,...,n-1,\,\qquad \sum_{k=0}^n a_k\otimes 1^{n-k}=0.
\end{equation}
\end{proposition}
\begin{proof}
Let   $a=\sum_{k=0}^n \iota_{k,k}^\M(a_k)\in \KK_\JJ(n,n)$,  $a_k\in \KK(X^{\otimes k})$, $k=0,...,n$, $n\in \N$, be such that $\widetilde{\pi}(a)=0$. We consider $(\pi,t)$, and thereby also $\widetilde{\pi}$, as a representation on a Hilbert space $H$ and take  advantage of the    family $\{P_m\}_{m\in\N}$ of decreasing projections   onto the essential subspaces of $\pi_{m,m}$,  cf.  page \pageref{remarks concerning projections}. Within the notation of  Proposition \ref{representation of TT(X) proposition}, we get
$$
0  =\|\widetilde{\pi}(a)\|=\max \left\{\max_{j=0,1,...,n-1}\big\{ \|\widetilde{\pi}(a)(P_j-P_{j+1})\|\big\},\, \big\|\widetilde{\pi}(a)P_n\big\| \right\}
$$
and thus  
\begin{equation}\label{relations to be proved2}
\overline{\pi}_{j,j}(\sum_{k=0}^{j}
   a_{k}\otimes 1^{j-k})(P_j-P_{j+1})=0, \,\, j=1,...,n-1,
      \quad \overline{\pi}_{n,n}(\sum_{k=0}^{n}
   a_{k}\otimes 1^{n-k})=0.
 \end{equation}
Since  $\pi$ is faithful, $\{\overline{\pi}_{i,j}\}_{i,j\in\N}$ is faithful.  Hence   \eqref{relations to be proved2} implies that $\sum_{k=0}^{n}
   a_{k}\otimes 1^{n-k} = 0$ and   each element $\overline{\pi}_{j,j}(\sum_{k=0}^{j}
   a_{k}\otimes 1^{j-k})$ is supported on $P_{j+1}H$, $j=0,...,n-1$. In particular, 
   $
   \left(\sum_{k=0}^{n-1}
   a_{k}\otimes 1^{n-1-k}\right)\otimes 1=\sum_{k=0}^{n}
   a_{k}\otimes 1^{n-k} -a_n \in \KK(X^{\otimes n})
   $ and thus by Lemma \ref{lemma almost destroying} ii), $\sum_{k=0}^{n-1}
   a_{k}\otimes 1^{n-1-k} \in \KK(X^{\otimes (n-1)}I')$. Similarly,
   if we assume that $\sum_{k=0}^{n-m}
   a_{k}\otimes 1^{n-m-k}\in \KK(X^{\otimes (n-m)}I')$, for  certain $m=0,...,n-1$, then
   $$
    \left(\sum_{k=0}^{n-m-1}
   a_{k}\otimes 1^{n-m-1-k}\right)\otimes 1=\sum_{k=0}^{n-m}
   a_{k}\otimes 1^{n-m-k}- a_{n-m}\in \KK(X^{\otimes (n-m-1)})
   $$ 
   and hence by Lemma \ref{lemma almost destroying} ii),  $\sum_{k=0}^{n-m-1}
   a_{k}\otimes 1^{n-m-1-k}\in \KK(X^{\otimes (n-m-1)}I')$.
   Thus by induction relations \eqref{relations to be proved} are satisfied.
\end{proof}
By passing to  quotients we may use the above proposition to  get a description of the kernel of $\widetilde{\pi}$ in a general situation. To this end  we use the following lemma. 
 \begin{lemma}\label{T-pairs lemma}
 Let  $(\pi,t)$ be  a  representation of $X$ and let 
\begin{equation}\label{T-pair of ideals from representation}
I=\ker\pi ,\qquad I'=\{a\in A: \pi(a)\in \pi_{1,1}(\KK(X))\}.
\end{equation}
Then  representation $(\pi,t)$ factors through to the faithful representation of $X/XI$ for which the ideal of coisomtericity is $I'/I$. In particular, the following relations hold
\begin{equation}\label{T-pairs condition}
I\,\, \textrm{ is }\,\,X\textrm{-invariant,}\qquad I\subset I'\subset q^{-1}(J(X/XI)\cap(\ker\phi^I)^\bot). 
\end{equation}
where  $q:A\to A/I$ is the quotient map and $\phi^I$ is the left action on $X/XI$.
\end{lemma}
\begin{proof}
See  \cite[Lem. 5.10]{katsura2}  and Corollary  \ref{Proposition by muhly and solel}. 
\end{proof}

 \begin{proposition}\label{proposition on ideals and lattice structure}
If $(\pi,t)$ is a representation of $X$ coisometric on $J$, then  the kernel of the corresponding right tensor representations $\widetilde{\pi}$ of $\KK_\JJ$ is uniquely  determined  by the ideals 
$$
I=\ker\pi \quad \textrm{ and } \quad I'=\{a\in A: \pi(a)\in \pi_{1,1}(\KK(X))\}.
$$
Namely, $\ker \widetilde{\pi}(n,n)$, $n\in \N$, consists of those elements that can be represented in a form $\sum_{k=0}^n \iota_{k,k}^\M(a_k)$, $a_k\in \KK(X^{\otimes k})$, $k=0,...,n$, where 
$$
\sum_{k=0}^j a_k\otimes 1^{j-k}\in \KK(X^{\otimes j}I'),\,\, j=0,...,n-1,\,\,\, \sum_{k=0}^n a_k\otimes 1^{n-k}\in \KK(X^{\otimes n}I).
$$
% where $$ I=\ker\pi ,\qquad I'=\pi^{-1}(\pi_{11}(\KK(X)). $$
\end{proposition}
\begin{proof}
Pass from $(\pi,t)$ to the (faithful) quotient representation $(\dot{\pi},\dot{t})$ of the quotient $C^*$-correspondence $X/XI$ over $A/I$. By   Proposition \ref{proposition on quotients of correspondences} this corresponds to passing from the right tensor representation $\{\pi_{n,m}\}_{n,m\in \N}$ of $\KK$ to the (faithful) right tensor representation $\{\dot{\pi}_{n,m}\}_{n,m\in \N}$ of $\KK/\TT_X(I)=\KK_{X/XI}$.  Applying Proposition \ref{proposition on ideals and lattice structure0} to representation $(\dot{\pi},\dot{t})$ we get the assertion  by  Lemma \ref{T-pairs lemma}.
\end{proof}
\begin{remark}\label{remark for katsura}
Using the above statement one may deduce that for a representation $[\pi,t]=\{\pi_{n,m}\}_{n,m\in \N}$ of $\KK_{X}$,  all $n\in \N$ and $k>0$ we have
$$
\pi^{-1}_{n,n}\left( \sum_{j=1}^{k}\pi_{n+j,n+j}(\KK_X(n+j,n+j))\right)=\pi^{-1}_{n,n}(\pi_{n+1,n+1}(\KK_X(n+1,n+1))),
$$
and this space is equal to $\KK(X^{\otimes n}I')$. Hence for right tensor representations of $\KK_X$ it suffices to check the condition \eqref{necessary condition for faithfulness} only for $k=1$.
\end{remark}
We adapt the notion of coisometricity to the notion of a $T$-pair \cite[Def. 5.6]{katsura2}.
\begin{definition}\label{definition of T-pairs}
A pair  $(I,I')$ of ideals in $A$ satisfying \eqref{T-pairs condition} is called a \emph{$T$-pair} of $X$.  We shall say that a $T$-pair $(I,I')$ is \emph{coisometric }on an ideal $J$  in $A$,  if $J\subset I'$.  In particular, an $O$-pair introduced  in \cite[Def. 5.21]{katsura2} is simply a $T$-pair coisometric on $J=(\ker\phi)^\bot \cap J(X)$.
\end{definition}
Let us note  that if $(I,I')$ is a $T$-pair coisometric on  $J$, then  $I$ is automatically $J$-saturated. Indeed,  in view of Lemma \ref{T-pairs lemma} we have  
$$
J\subset q^{-1}(J(X/XI)\cap(\ker\phi^I)^\bot) \quad \Longrightarrow \quad J\cap\varphi^{-1}(\KK(XI))\subset I.
$$ 
Furthermore, $T$-pairs form a lattice with the natural order induced by inclusion.   
\begin{theorem}[Lattice structure of gauge invariant ideals in $\OO(J,X)$, cf.  \cite{katsura2}]\label{Lattice structure of gauge invariant ideals in O(J,X) theorem}
We have  lattice isomorphisms   between the following objects  
\begin{itemize}
\item[i)] $T$-pairs $(I,I')$ of $X$ coisometric on $J$, 
\item[ii)] invariant and $\KK_\JJ$-saturated  ideals $\NN_\JJ$ in $\KK_\JJ$,
\item[iii)] gauge invariant ideals $P$ in $\OO(J,X)$.
\end{itemize}
 The correspondences between the objects in ii) and iii), and  i) and ii) are respectively  given by the equality $
\OO(\NN_\JJ)=P$, and the equivalence: 
$a\in \KK_\JJ(n,n)$ is in $\NN_\JJ(n,n)$, $n\in \N$, iff it  may be represented by $\sum_{k=0}^n \iota_{k,k}(a_k)$ where
$$
\sum_{k=0}^j a_k\otimes 1^{j-k}\in \KK(X^{\otimes j}I'),\,\, j=0,...,n-1,\,\,\,\,\,\,\, \sum_{k=0}^n a_k\otimes 1^{n-k}\in \KK(X^{\otimes n}I).
$$
Moreover, we have
$$
P\cong \DR(\NN_\JJ),\qquad \OO(J,X)/P\cong \OO(I'/I,X/XI).
$$ 
 \end{theorem} 
 \begin{proof}
In view of Theorem \ref{theorem on lattice structure} and Proposition \ref{proposition on ideals and lattice structure} it suffices to  show that the pair of ideals $(I,I')$ from item i) define (via the described equivalence) invariant and $\KK_\JJ$-saturated ideal $\NN_\JJ$ in $\KK_\JJ$. To this end note that by  \eqref{T-pairs condition} and   Corollary \ref{Kernel of universal representation correspondence} the universal representation of $X/XI$ in $\OO(I'/I,X/XI)$ is faithful. Composing this representation with   quotient maps $X\to X/XI$, $A\to A/I$ one gets    the representation $(\pi,t)$ of $X$ such that that  relations \eqref{T-pair of ideals from representation} are satisfied. The pair $(\pi,t)$ give rise to  representation $\widetilde{\pi}$ of $\KK_\JJ$ whose kernel  (by Proposition \ref{proposition on ideals and lattice structure}) is the desired ideal $\NN_\JJ$. Using Theorem \ref{Gauge invariance theorem for O_T(I,J)1} we get  $\OO(J,X)/P\cong \OO(I'/I,X/XI)$. \end{proof}

\begin{remark}The ideals  we identified in Theorem \ref{structure theorem2} as algebras of the form $\OO_X(I,J\cap I)$ are exactly these gauge invariant ideals in  $\OO(J,X)$ that correspond to $T$-pairs $(I,I')$ where $I'=I+J$. As we noticed in Theorem \ref{sometimes all ideals are of simple form} in many situations all gauge invariant ideals are of this form. 
\end{remark}

\section{Embedding conditions for $\OO_\TT(\KK,\JJ)$}\label{Embedding theorems for O_T(I,J)}
We fix a right tensor $C^*$-precategory $\TT$ and exhibit    conditions for algebras of  type $\OO_\TT(\KK,\JJ)$ to be embedded into one another via  universal representations.   These results have similar motivation as  \cite[Prop. 3.2]{dpz}, \cite[Prop. 6.3]{fmr} and shall be applied in Section \ref{Representations of algebras associated with X} to algebras associated with $C^*$-correspondences. They also may be viewed as a description of certain gauge invariant subalgebras of $\OO_\TT(\KK,\JJ)$. We start with
\begin{proposition}[Necessary conditions]\label{Necessary conditions proposition}
For $j=1,2$, let $\KK_{\,j}$ and $\JJ_{j}$ be ideals in $\TT$ such that 
$ \JJ_{j}\subset J(\KK_{\,j})$ and $\KK_{\,1}\subset \KK_{\,2}$. Denote by $\{\iota_{n,m}^{(j)}\}_{m,n\in \N}$  the universal representations of $\KK_{j}$  in $\OO_\TT(\KK_{\,j},\JJ_{j})$, $j=1,2$.   The natural homomorphism $\Psi:\OO_\TT(\KK_{\,1},\JJ_{1})\longmapsto \OO_\TT(\KK_{\,2},\JJ_2)$
\begin{equation}\label{quasi embedding 1}
\OO_\TT(\KK_{\,1},\JJ_{1})\ni  \iota_{n,m}^{(1)}(a)\stackrel{\Psi}{\longmapsto} \iota_{n,m}^{(2)}(a)\in \OO_\TT(\KK_{\,2},\JJ_2)
\end{equation}
is well define if and only if $
\JJ_1\subset \JJ_2.$
Moreover, if $\Psi$ is well defined and injective, then 
$$
\JJ_1=\JJ_2\cap J(\KK_{\,1}), \qquad \RR_{\JJ_2}\cap \KK_{\,1}=\RR_{\JJ_1}  
$$
 where $\RR_{\JJ_j}$  is  the reduction ideal associated with $\JJ_j$, $j=1,2$, see Definition \ref{reduction ideal definition}.
 \end{proposition}
 \begin{proof}
If  $
\JJ_1\subset \JJ_2
$,  then $\Psi$ is well defined by the construction of  norm  in $\OO_\TT(\KK_{\,j},\JJ_{j})$, $j=1,2$, see Proposition \ref{propostion nie udowodnione}. Conversely, if $\Psi$ is well defined, then for $a\in J(\KK_{\,1})(n,m)$ we have
$$
  \iota_{n,m}^{(1)}(a)=\iota_{n+1,m+1}^{(1)}(a\otimes 1) \,\, \,\Longrightarrow\,\,\, 
\iota_{n,m}^{(2)}(a)=\iota_{n+1,m+1}^{(2)}(a\otimes 1),
$$
that is $\JJ_1 \subset \JJ_2$. In the event  $\Psi$ is injective the above implication is an equivalence and thus we have $\JJ_1 = \JJ_2\cap J(\KK_{\,1})$. Moreover, representations   $\{\iota_{n,m}^{(1)}\}_{m,n\in \N}$  and  $\{ \iota_{n,m}^{(2)}|_{\KK_1}\}_{m,n\in \N}=\{\Psi \circ \iota_{n,m}^{(1)}\}_{m,n\in \N}$ have the same kernels  and hence  we get $\RR_{\JJ_2}\cap \KK_{\,1}=\RR_{\JJ_1} $.
 \end{proof}
 By the above statement  we may narrow our attention down to algebras  $\OO_\TT(\KK_{\,1},\JJ\cap J(\KK_{1}))$ and $\OO_\TT(\KK_{\,2},\JJ)$ where $\KK_1\subset \KK_2$ and $\JJ\subset J(\KK_2)$. 
 \begin{proposition}[Sufficient conditions on ideals $\KK_1$ and $\JJ$] \label{proposition of hierarchy of algebras2} 
 Let $\KK_1$, $\KK_2$ and $\JJ$ be ideals in $\TT$ such that $\KK_1\subset \KK_2$ and
$ \JJ\subset J(\KK_2)$. The  condition
\begin{equation}\label{sufficient conditions on ideals equation1}
\KK_1\cap J(\JJ) \subset J(\KK_1)
\end{equation}
(which holds e.g. whenever  $\JJ\subset \KK_1$,  $\JJ\subset \KK_1^\bot$ or $\KK_1$ is invariant) implies that there is a natural embedding 
$$
\OO_\TT(\KK_1,\JJ\cap J(\KK_1))\subset  \OO_\TT(\KK_2,\JJ).
$$
\end{proposition}
\begin{proof}
Notice that $\KK_1\cap J(\JJ) \subset J(\KK_1)$ is equivalent to $
\JJ\cap \KK_1 \cap (\otimes 1)^{-1}(\JJ+\KK_1)\subset J(\KK_1)
$. Indeed, on the one hand  we have  $\KK_1\cap J(\JJ) \subset \JJ\cap \KK_1 \cap (\otimes 1)^{-1}(\JJ+\KK_1)$, and on the other if $a \in   \JJ\cap \KK_1 \cap (\otimes 1)^{-1}(\JJ+\KK_1)(n,m)$, then either $a\in J(\KK_1)$ or $a\in \KK_1\cap J(\JJ)$: for $a\in \TT(n,m)$, $m,n\in \N$, we have
\begin{equation}\label{Sufficient condition on ideals I_1 and J}
a\in (\KK_1\cap \JJ) \textrm{ and }a\otimes 1 \in (\KK_1 + \JJ)\,\, \Longrightarrow a\in J(\KK_1).
\end{equation}
We show that $\Psi$ given by \eqref{quasi embedding 1} is faithful on the core (and hence on all of the spectral subspaces) of    $\OO_\TT(\KK_1,\JJ\cap J(\KK_1))$. To this end let   $a\in \OO_\TT(\KK_1,\JJ\cap J(\KK_1))$ be of the form $a= \sum_{j\in \N}
   \iota_{j,j}(a_{j})$,  $a_{j}\in \KK_1(j,j)$, $j\in \N$, and suppose that $\|\Psi(a)\|=0$. Then  by  Theorem \ref{proposition of norm formulas}  we have
 $$
  \sum_{j=0}^{s}
   a_{j}\otimes 1^{s-j}\in \JJ(s,s), \quad s\in \N,\qquad \lim_{r\to \infty}\sum_{j=0}^{r}
   a_{j}\otimes 1^{r-j}=0. 
$$
Since $a_{0}\in (\KK_1\cap \JJ)(0,0)$ and $a_{0}\otimes 1= (a_{0}\otimes 1 + a_{1}) -a_{1}\in  (\KK_1+\JJ)(1,1)$, by \eqref{Sufficient condition on ideals I_1 and J} we get $a_{0}\in J(\KK_1) (1,1)$ and consequently  $a_{0}\otimes 1+a_{1}\in (\JJ\cap \KK_1)(1,1)$. Suppose  now that  $\sum_{j=0}^{s}
   a_{j}\otimes 1^{s-j}\in (\JJ\cap \KK_1)(s,s)$ for certain $s\in \N$. Since   
   $$
  \left(\sum_{j=0}^{s}
   a_{j}\otimes 1^{s-j}\right)\otimes 1 =\sum_{j=0}^{s+1}
   a_{j}\otimes 1^{s+1-j}- a_{s+1} \in (\JJ+\KK)(s+1,s+1),$$
    by \eqref{Sufficient condition on ideals I_1 and J} we then  have $\sum_{j=0}^{s}
   a_{j}\otimes 1^{s-j}\in J(\KK_1)(s+1,s+1)$ and  consequently   $\sum_{j=0}^{s+1}
   a_{j}\otimes 1^{s+1-j}\in (\JJ\cap \KK_1)(s+1,s+1)$. Thus by   induction we get 
   $$
  \sum_{j=0}^{s}
   a_{j}\otimes 1^{s-j}\in (\JJ\cap J(\KK_1))(s,s), \quad s\in \N,\qquad \lim_{r\to \infty}\sum_{j=0}^{r}
   a_{j}\otimes 1^{r-j}=0 
$$ 
which   in view of  Theorem \ref{proposition of norm formulas} is  equivalent to $\|a\|=0$. \\
Clearly, $\Psi:\OO_\TT(\KK_1,\JJ\cap J(\KK_1))\to  \OO_\TT(\KK_2,\JJ)$ preserves the gauge actions and hence injectivity of $ \Psi$  on spectral subspaces implies the injectivity of $ \Psi$ on the whole algebra $\OO_\TT(\KK_1,\JJ\cap J(\KK_1))$.
     \end{proof}
\begin{corollary}
For any ideals $J\subset J(X)$  and $I$ in $A$ such that 
\begin{equation}\label{sufficient conditions on ideals equation2}
I\cap J(XJ) \subset J(XI) 
\end{equation}
(which holds e.g. whenever $J\subset I$,  $J\subset I^\bot$ or $I$ is $X$-invariant) we have the  natural embedding 
$$
\OO_X(I,J\cap J(XI))\subset \OO(J,X).
$$
\end{corollary}
\begin{proof}
Apply Proposition \ref{proposition of hierarchy of algebras2}  to
 $\KK_2=\KK_X$, $\KK_1=\KK_X(I)$ and $\JJ=\KK_X(J)$.
\end{proof}
We notice that if $I$ is $X$-invariant, then  for  $J\subset J(X)$  we have $J\cap J(XI)=J\cap I$   and thereby the foregoing statement implies the inclusion $\OO_X(I,J\cap I)\subset \OO(J,X)$ from Theorem \ref{structure theorem2}. 
As the next example shows the condition \eqref{sufficient conditions on ideals equation1}, or more precisely its special case \eqref{sufficient conditions on ideals equation2},  is essential.
\begin{example}
Suppose $A=A_0\oplus A_0$ where $A_0$ is a unital $C^*$-algebra. Consider the ideals $J=A_0\oplus A_0$,  $I= A_0 \oplus\{0\}$   and  the $C^*$-correspondence $X=X_\al$  associated with  the endomorphism $\al:A\to A$ given  by the formula $\al(a,b)=(0,a)$. Then   we have
$$
I\cap J(XJ)=A_0 \oplus\{0\}\nsubseteq J(XI)= \{0\}\oplus A_0,
$$
so the inclusion \eqref{sufficient conditions on ideals equation2} fails.
On the other hand the algebra $\OO_X(I,J\cap J(X I))$ can not be embedded into  $
C^*(J,X)$ as we have   $
C^*(J,X)=\{0\}$ and $\OO_X(I,J\cap J(X I))= I$
(the latter relation may be checked using,  for instance,  Theorem \ref{proposition of norm formulas}).
\end{example}
Let us now exploit the  condition $\RR_\JJ\cap \KK=\RR_{\JJ\cap J(\KK_{\,1})}$ introduced in  Proposition \ref{Necessary conditions proposition}. It is evident that $\JJ\subset \KK_1$ implies  $\RR_\JJ=\RR_{\JJ\cap J(\KK_{\,1})}$ (it  follows  from Propositions \ref{Necessary conditions proposition}, \ref{proposition of hierarchy of algebras2}). Moreover,  we have
$$
\JJ\subset (\ker \otimes 1)^{\bot} \,\, \,\Longrightarrow\,\,\, \RR_\JJ=\RR_{\JJ\cap J(\KK_{\,1})}=\{0\}.
$$ 
Thus for $C^*$-algebras  generalizing Katsura's algebras of $C^*$-correspondences (cf. Remark  \ref{remark on definition of main object})
the  condition  $\RR_\JJ=\RR_{\JJ\cap J(\KK_{\,1})}$ is trivially  satisfied. 
\begin{proposition}[Sufficient conditions on ideal $\KK_1$] \label{proposition of hierarchy of algebras} 
 Let $\KK_1$  be an ideal in $\TT$ such that 
\begin{equation}\label{condition sufficient for embeddings}
\otimes^{-1}(\KK_1)\cap (\ker\otimes 1)^\bot \subset \KK_1.
\end{equation}
Then for every ideal  $\KK_2$ in $\TT$ and every ideal $\JJ$ in $J(\KK_2)$ the natural homomorphism \eqref{quasi embedding 1} establish the embedding 
$$
\OO_\TT(\KK_1,\JJ\cap J(\KK_1))\subset  \OO_\TT(\KK_2,\JJ)
$$
if and only if $ \RR_{\JJ}\cap \KK_1=\RR_{\JJ\cap J(\KK_1)}$.
\end{proposition}
\begin{proof}
 The "only if" part follows from  Proposition \ref{Necessary conditions proposition}. Let us   assume that  $\RR_\JJ\cap \KK_1=\RR_{\JJ\cap J(\KK_1)}$.  By  the reduction procedure described in Theorem \ref{reduction of relations theorem}  applied to $\JJ$ we may actually  assume that $\RR_\JJ=\{0\}$ and 
$\JJ \subset (\ker\otimes 1)^\bot $. %(it is clear that if $\KK_1$ is invariant, then so is $\KK_1/\RR_\JJ$, and one readily checks that condition i) implies analogous relation for $\KK_1/\RR_\JJ$). 
% Now, similarly as in the proof of Proposition \ref{proposition of hierarchy of algebras2}, 
It suffices to prove that if $a\in \OO_\TT(\KK_1,\JJ\cap J(\KK_1))$ is such that  $a= \sum_{s=0}^r
   \iota_{s,s}(a_{s})$,  $a_{s}\in \KK_1(s,s)$, $s=0,...,r$, and  $\|\Psi(a)\|=0$, then $\|a\|=0$. In view of   Theorem \ref{proposition of norm formulas}  the requirement  $\|\Psi(a)\|=0$ is equivalent to
\begin{equation}\label{equation pomocniczy1}
  \sum_{j=0}^{s}
   a_{j}\otimes 1^{s-j}\in \JJ(s,s), \quad s=0,...,r-1,\qquad \sum_{j=0}^{r}
   a_{j}\otimes 1^{r-j}=0. 
\end{equation}
To show that $\|a\|=0$ we need to check whether
\begin{equation}\label{equation pomocniczy}
\sum_{j=0}^{s}
   a_{j}\otimes 1^{s-j}\in (\JJ\cap J(\KK_1))(s,s), \quad s=0,...,r-1,\qquad \sum_{j=0}^{r}
   a_{j}\otimes 1^{r-j}=0.
\end{equation}
However, since  
$$
\Big(\sum_{j=0}^{r-1}
   a_{j}\otimes 1^{r-1-j}\Big)\otimes 1 =\sum_{j=0}^{r}
   a_{j}\otimes 1^{r-j}- a_{r}= - a_{r}\in  \KK_1(r,r),
   $$
it follows that $\sum_{j=0}^{r-1}
   a_{j}\otimes 1^{r-1-j}\in \Big(\otimes^{-1}(\KK_1)\cap (\ker\otimes 1)^\bot\Big)(r-1,r-1) \subset J(\KK_1)(r-1,r-1)$. 
Analogously, if $\sum_{j=0}^{s}
   a_{j}\otimes 1^{s-j}\in  J(\KK_1)(s,s)$, $s=0,...,r-1$, then
   $$
\Big(\sum_{j=0}^{s-1}
   a_{j}\otimes 1^{s-1-j}\Big)\otimes 1 =\sum_{j=0}^{s}
   a_{j}\otimes 1^{s-j}- a_{s,s}\in  \KK_1(r,r)
   $$
   and  $\sum_{j=0}^{s-1}
   a_{j}\otimes 1^{s-1-j}\in \Big(\otimes^{-1}(\KK_1)\cap (\ker\otimes 1)^\bot\Big)(s-1,s-1) \subset J(\KK_1)(s-1,s-1)$.
   Hence  by induction \eqref{equation pomocniczy} holds.
\end{proof}
 Item iii) of Lemma \ref{lemma embedding of tensors}  states that  the ideal  $\KK_X$  in  $\TT_X$ satisfies \eqref{condition sufficient for embeddings},  and thus we get
\begin{proposition}\label{corollary about correspondences}
Let $\TT=\TT_X$ be a right tensor $C^*$-precategory of a $C^*$-correspon\-dence $X$ and  let 
  $\SSS$ and $\JJ$ be  ideals in $\TT$ such that $\JJ\subset J(\KK_2)$  and $\KK_X\subset \KK_2$. If we  put  $J:=\JJ(0,0)$, then the natural homomorphism  is an embedding 
$$
\OO(J\cap J(X),X)\subset  \OO_{\TT}(\KK_2,\JJ),
$$
 if and only if $R_J=R_{J\cap J(X)}$ where $R_J$ and $R_{J\cap J(X)}$ are reducing ideals  associated to $J$ and $J\cap J(X)$ respectively, see Definition \ref{reduction procedure for X definition}.  
\end{proposition}
In view of Proposition \ref{roberts and doplicher wear a T-shirt}, the algebra  $\OO_{\TT_X}(\TT_X, \TT_X)$ coincides with the \emph{Doplicher-Roberts algebra of a $C^*$-correspondence} investigated in \cite{fmr}, \cite{dpz}. It is natural to consider the following "relative version" of such algebras.
\begin{definition}\label{relative Doplicher-Roberts algebra definition} Suppose $X$ is a $C^*$-correspondence over $A$ and  $J$ is an arbitrary ideal in $A$. We shall call the $C^*$-algebra  
$$
\DR(J, X):= \OO_{\TT_X}(\TT_X, \TT_X(J))
$$
 a \emph{relative Doplicher-Roberts algebra of $X$} relative to $J$.  Within this notation the algebra considered in \cite{fmr}, \cite{dpz} is $\DR(A, X)$.
 \end{definition}
Now Proposition \ref{corollary about correspondences} can be interpreted  as the following generalization of   \cite[Prop. 6.3]{fmr}, \cite[Prop. 3.2]{dpz}.
 \begin{corollary}\label{corollary about correspondences2}
 The natural homomorphism  is an embedding 
$$
\OO(J\cap J(X),X)\subset \DR(J, X).
$$
 if and only if $R_J=R_{J\cap J(X)}$.  In particular, $\OO(J\cap J(X),X)$ embeds into $\DR(J, X)$   whenever
 $
 J\subset (\ker \phi)^{\bot}$ or $J\subset J(X)$.
\end{corollary}
\begin{example}
Let $A_0$ be a non-unital  $C^*$-algebra and $A_0^+$ its minimal unitization. Let us consider the $C^*$-correspondence   $X=A_0\oplus A_0^+$ over $A=A_0^+\oplus A_0^+$ where 
$
\langle x,y \rangle_A :=x^*y$, $x \cdot a= xa$ and  $a\cdot x =\al(a)x$ where $\al(a_0\oplus b_0)=0\oplus a_0$. Then  $J(X)=A_0^+\oplus A_0^+$ and for the ideal $J=A$ we have 
$$
R_J=A_0^+\oplus A_0^+ \neq A_0\oplus A_0^+=R_{J\cap J(X)}.
$$
On the other hand   $\OO(J\cap J(X),X)\cong \C$ and  $\DR(J, X)=\{0\}$.
\end{example}

\section{Application to Doplicher-Roberts algebras associated with  $C^*$-cor\-res\-pon\-den\-ces}\label{Representations of algebras associated with X}
 Our study in this section   is  motivated by \cite{fmr}, \cite{dpz} and our aim is to generalize \cite[Thm. 6.6]{fmr}, \cite[Thm. 4.1]{dpz}.
We recall 
that there is a  one-to-one correspondence, 
 given by relations
$$
\pi=\Psi\circ \iota_{0,0}|_A,\qquad t=\Psi\circ \iota_{1,0}|_X, 
$$
between representations $\Psi$ of the algebra $\OO(J,X)$ and   representations $(\pi,t)$ of $X$ coisometric on $J$.
The inverse of this correspondence  is given by  the equality $\Psi=\Psi_{[\pi,t]}$ where $[\pi,t]$ is the right tensor  representation of $\KK_X$ defined in Proposition \ref{proposition:rho}. For shortening we denote the representation of $\OO(J,X)$  corresponding to $(\pi,t)$  by $\pi\times_J t$. We arrive at the following statement, cf.  \cite[Cor. 11.7]{katsura2}, \cite[Thm. 5.1]{mt}, \cite[6.4]{katsura}. 
\begin{theorem}[Gauge invariance theorem for $\OO(J,X)$] \label{Gauge invariance theorem for O(J,X) hey}
Let $(\pi,t)$ be a representation of a $C^*$-correspondence $X$ coisometric on an ideal $J \subset J(X)$.  Then 
 $\pi\times_J t$ is  a faithful representation of $\OO(J,X)$ if and only if 
 $(\pi,t)$ admits a gauge action,  $\ker\pi=R_J$ %(the reduction ideal associated to $\JJ$
 and 
 \begin{equation}\label{necessary condition for faithfulness 1}
 J=\{a\in A: \pi(a) \in \pi_{1,1}(\KK(X))\}.
\end{equation}
In particular, 
\begin{itemize}
\item[i)] if $J\subset (\ker\varphi)^{\bot}$, then  $\pi\times_J t$ is faithful  if and only if $(\pi,t)$ is faithful, admits a gauge action and $J$ is an ideal of coisometricity for $(\pi,t)$.
\item[ii)] if $J+ (\ker\varphi)^{\bot}=A$, then $\pi\times_J t$ is faithful if and only if $(\pi,t)$ admits a gauge action, $\ker\pi=R_J$  and $J$ is an ideal of coisometricity for $(\pi,t)$.
\end{itemize}
\end{theorem}
\begin{proof}
For the first part use Theorem \ref{Gauge invariance theorem for O_T(I,J)1} and Proposition \ref{proposition on ideals and lattice structure}. To see items i), ii) apply    Proposition \ref{proposition on ideals and lattice structure0} and Theorem  \ref{sometimes all ideals are of simple form}.
\end{proof}
\begin{remark}
One could state the foregoing  theorem in a somewhat more general form similar to Theorem \ref{Gauge invariance theorem for O_T(I,J)1}. Namely, if in the above statement one drops the requirement of admitting a gauge action, one gets necessary and sufficient conditions  for $\pi\times_J t$ to be  faithful on spectral subspaces of $\OO(J,X)$.
\end{remark}
If $(\pi,t)$ is a representation   on a Hilbert space $H$, then   $[\pi,t]=\{\pi_{n,m}\}_{m,n\in\N}$  extends to a  right tensor representation $\overline{[\pi,t]}=\{\overline{\pi}_{n,m}\}_{m,n\in\N}$ of $\TT_X$, see Proposition \ref{extensions of representations on Hilbert spaces}.  
Therefore, for appropriately chosen ideal  $J$,  $\overline{[\pi,t]}$ integrates to a representation of relative Doplicher-Roberts algebra $\DR(J ,X)$ that we shall denote by $\overline{\pi\times_J t}$. By Proposition \ref{representation of TT(X) proposition} we have
\begin{proposition}\label{combining of two representations}  
Let $J$ be an ideal  in $A$.  A representation $\Psi$ of  $\DR(J, X)$ on a Hilbert space $H$ is of the form $\overline{t\times_J\pi}$ for a representation $(\pi,t)$ of $X$ 
if and only if
$$
\Psi(\iota_{n,n}(\LL( X^{\otimes n}))H=\Psi(\iota_{n,n}(\KK( X^{\otimes n}))H,\qquad n\in \N,
$$ 
where  $\{\iota_{n,m}\}_{m,n\in \N}$ is the universal representations of $\TT_X$ in $\DR(J, X)$.
%and $J\subset \{a\in A: \Psi(i_{0,0}(a))=\Psi(i_{1,1}(\phi(a))\}$.
If this is the case and additionally  $R_J=R_{J_0}$ where $J_0=J\cap J(X)$, then
  $\overline{t\times_J\pi}$ is an extension of $ t\times_{J_0}\pi$ (as we then have  $\OO(J_0, X) \subset \DR(J, X)$, cf. Corollary \ref{corollary about correspondences2}). 
\end{proposition}
%\begin{proof}Apply Corollary \ref{representations extension core}\end{proof}
We note that 
\begin{align*}
 \{a\in A: \overline{\pi}_{1,1}(\phi(a))=\pi(a)\}&=\{a\in A: \pi(a)P_1=\pi(a)\}
 \\
 =\{a\in A: \pi(a) \in \overline{\pi}_{1,1}(\LL(X))\} &=\{a\in A: \pi(a)H\subset \overline{t(X)H}\},
\end{align*}
 cf. Lemma \ref{lemma almost destroying} i) or \cite[Lem. 1.9]{fmr}, so  the forthcoming results may be stated without a use of the mapping $\overline{\pi}_{1,1}$. We shall however keep using it as it indicates the relationship with coisometricity for $C^*$-correspondences. 
The following  statement generalizes \cite[Thm. 6.6]{fmr}.
\begin{theorem}[Extending  representations from $\OO(J_0,X)$  to  $\DR(J, X)$]\label{main theorem of this section} Let  $(\pi,t)$ be a representation of $X$ on a Hilbert space $H$ and let  $J$ be an ideal in $A$ such that
$$
 J\subset \{a\in A: \overline{\pi}_{1,1}(\phi(a))=\pi(a)\},\quad \textrm{and}\quad R_J=R_{J_0}\quad \textrm{where}\quad J_0=J\cap J(X).
$$
Then $\OO(J_0, X) \subset \DR(J, X)$ and the representation $\overline{\pi\times_J t}$ of $\DR(J, X)$ is an extension of the representation  $\pi\times_{J_0} t$  of $\OO(J_0,X)$
 % then representation $\overline{(t\times_J\pi)}$ of $\DR(J, X)$ is an extension of representation  $t\times_{J\cap J(X)}\pi$  of $\OO(J \cap J(X),X)\subset \DR(J, X) $ and the following condition holds
 and 
 \begin{itemize}
 \item[i)]  $\overline{\pi\times_J t}$ is faithful on spectral subspaces of $\DR(J, X)$ if and only if  $\pi\times_{J_0} t$ is faithful on spectral subspaces of $\OO(J_0,X)$ and $J= \{a\in A: \overline{\pi}_{1,1}(\phi(a))=\pi(a)\}$.
 \item[ii)]  $\overline{\pi\times_J t}$ is faithful  if and only if   $\pi\times_{J_0} t$ is faithful and $J= \{a\in A: \overline{\pi}_{1,1}(\phi(a))=\pi(a)\}$.
\end{itemize}

\end{theorem}
\begin{proof}
i). Since $\pi\times_{J_0} t$ is considered as a restriction of  $\overline{\pi\times_J t}$,  faithfulness of $\overline{\pi\times_J t}$  on  spectral subspaces, implies such a faithfulness  of  $\pi\times_{J_0} t$, and then  we get $J= \{a\in A: \overline{\pi}_{1,1}(\phi(a))=\pi(a)\}$ by Proposition \ref{another one bites the proposition}. Conversely, suppose  that $J= \{a\in A: \overline{\pi}_{1,1}(\phi(a))=\pi(a)\}$ and  $\pi\times_{J_0} t$ is faithful on spectral subspaces of  $\OO(J_0,X)$. In this event  $\ker\pi=R_J=R_{J_0}$ and hence by passing to quotients  (dividing  all the associated $C^*$-precategories by $\NN=\TT_X(\ker \pi)$) we may  assume that  $\ker\pi=0$. Then, by  Proposition \ref{representation of TT(X) proposition}, representation $\{\overline{\pi}_{i,j}\}_{i,j\in\N}$ is faithful. To show that   $\overline{\pi\times_J t}$ is faithful let $a= \sum_{s=0}^{r}
   \iota_{s,s}(a_{s})$,  where $a_{s}\in \LL(X^{\otimes s})$, $s=0,...,r$, $r\in \N$, be such that  $\overline{\pi\times_J t}(a)=0$.  Similarly as in the proof of Proposition \ref{proposition on ideals and lattice structure0} we get
%\begin{equation}
$$
\overline{\pi}_{s,s}(\sum_{j=0}^{s}
   a_{j}\otimes 1^{s-j})(P_s-P_{s+1})=0, \,\,\,\, s=1,...,r-1,
      \qquad \overline{\pi}_{r,r}(\sum_{j=0}^{r}
   a_{j}\otimes 1^{r-j})=0.
 $$
 %\end{equation}
In particular,  $\sum_{j=0}^{r}
   a_{j}\otimes 1^{r-j} = 0$ and  each operator $\overline{\pi}_{s,s}(\sum_{j=0}^{s}
   a_{j}\otimes 1^{s-j})$ is supported on $P_{s+1}H$, $s=0,...,r-1$. Thus applying inductively  Lemma \ref{lemma almost destroying} i) we get
   $$  \sum_{j=0}^{s}
   a_{j}\otimes 1^{s-j}\in \LL(X^{\otimes s},X^{\otimes s}J), \quad s=1,...,r-1,\qquad \sum_{j=0}^{r}
   a_{j}\otimes 1^{r-j}=0,
$$
which is equivalent to $\|a\|=0$.\\
ii). If $\overline{\pi\times_J t}$ is faithful, then    $\pi\times_{J_0}t$ is faithful  and $J= \{a\in A: \overline{\pi}_{1,1}(\phi(a))=\pi(a)\}$ by Proposition \ref{another one bites the proposition}.
Conversely, if   $J= \{a\in A: \overline{\pi}_{1,1}(\phi(a))=\pi(a)\}$ and $\pi \times_{J_0} t$ is faithful, then  $\overline{\pi\times_J t}$ is faithful on spectral subspaces of $\DR(J,X)$ by item i). Thus  it  suffices to show  that
\begin{equation}\label{inequality * to be proved}
 \|\overline{\pi\times_J t}(E(a))\|\leq \|\overline{\pi\times_J t}(a)\|,\qquad a\in \DR(J,X),
\end{equation}
where $E$ is the  conditional expectation onto the $0$-spectral subspace of $\DR(J,X)$, see \cite[Thm 4.2]{dr}. For that purpose we  note three facts. 
Firstly, we have %  
\begin{equation}\label{inequality * to be used}
 \|\overline{\pi\times_J t}(E(a))\|\leq \|\overline{\pi\times_J t}(a)\|,\qquad a\in \OO(J_0,X)\subset \DR(J,X),\end{equation}
(a wicker version of \eqref{inequality * to be proved}) because  $\overline{\pi\times_J t}$ is an extension of $\pi \times_{J_0} t$ and $\pi \times_{J_0} t$ is  faithful.
Secondly,  for any $a\in \DR(J,X)$ and any $r\in \N$ we have 
\begin{equation}\label{maximum norms to help}
\|\overline{\pi\times_J t}(a)\|=\max \left\{\max_{s=0,1,...,r-1}\big\{ \|\overline{\pi\times_J t}(a)(P_s-P_{s+1})\|\big\},\, \big\|\overline{\pi\times_J t}(a)P_r\big\| \right\},
\end{equation}
where  $\{P_m\}_{m\in\N}$ is the family of decreasing projections    described  on  page \pageref{remarks concerning projections}. Thirdly, it is enough to prove \eqref{inequality * to be proved} for  elements $a\in \DR(J,X)$ of  the form  
\begin{equation}\label{form of elements to be mealed with}
a=\sum_{k=-\infty}^{\infty} \sum_{j=0, \atop j+k\geq 0}^{r}
   \iota_{j+k,j}(a_{j+k,j}),\qquad  a_{j+k,j}\in \LL(X^{\otimes j},X^{\otimes j+k}), \,\, r\in \N,
 \end{equation}
 as they  form a dense subspace of $\DR(J,X)$.\\
Hence we fix an element $a$ of the form \eqref{form of elements to be mealed with},
  $\varepsilon >0$ and $\xi \in (P_s-P_{s+1})H$ where $s=0,1,...,r-1$. Since $\pi_{s,s}$ restricted to $P_sH$ is nondegenerate, by  Hewitt-Cohen Factorization Theorem
  there exist   $c\in \KK(X^{\otimes s})$ and $\eta \in H$ such that  $\xi=\pi_{s,s}(c)\eta$  and   $$
  \| \pi_{s,s}(c)\|\cdot  \|\eta\|\leq (1+\varepsilon) \cdot  \|\xi\|.
   $$
We have 
\begin{align*}
   \| \overline{\pi\times_J t}(E(a))\xi\| &=  \| \overline{\pi\times_J t}(E(a))(P_s-P_{s+1})\xi\|=
   \|
   \overline{\pi}_{s,s}\big(\sum_{j=0}^{s}
   a_{j,j}\otimes 1^{s-j}\big)\xi\|
 \\
   &=       \| \overline{\pi}_{s,s}\big(\sum_{j=0}^{s}
   a_{j,j}\otimes 1^{s-j}\big)\pi_{s,s}(c)\eta\|
   =
    \|\pi_{s,s}\Big(\big(\sum_{j=0}^{s}
   a_{j,j}\otimes 1^{s-j}\big) c\Big)\eta\|
  \\
    & \leq \|\overline{\pi\times_J t}\Big(\iota_{s,s}\big(\sum_{j=0}^{s}
   a_{j,j}\otimes 1^{s-j} c\big)\Big)\| \cdot \|\eta\|
   \\
  &  = \|\overline{\pi\times_J t}\Big(E\big(\sum_{k=-\infty}^{+\infty}\iota_{s+k,s}(\sum_{j=0, \atop j+k\geq 0}^{s}
   a_{j+k,j}\otimes 1^{s-j} c)\big)\Big)\| \cdot \|\eta\|.
   \end{align*}
Applying  \eqref{inequality * to be used} to $\sum_{k=-\infty}^{+\infty}\iota_{s+k,s}(\sum_{j=0, \atop j+k\geq 0}^{s}
   a_{j+k,j}\otimes 1^{s-j} c) \in \OO(J_0,X)$ we get
\begin{align*}
  \| \overline{\pi\times_J t}(E(a))\xi\| & \leq \|\overline{\pi\times_J t}\Big(\sum_{k=-\infty}^{+\infty}\iota_{s+k,s}(\sum_{j=0, \atop j+k\geq 0}^{s}
   a_{j+k,j}\otimes 1^{s-j} c)\Big)\| \cdot \|\eta\|
  \\
  & \leq \|\sum_{k=-\infty}^{+\infty}\overline{\pi}_{s+k,s}\big(\sum_{j=0, \atop j+k\geq 0}^{s}
   a_{j+k,j}\otimes 1^{s-j}\big)\|\cdot \|\pi_{s,s}(c)\| \cdot \|\eta\|
   \\
   &\leq \|\sum_{k=-\infty}^{+\infty}\sum_{j=0, \atop j+k\geq 0}^{r}
   \overline{\pi}_{j+k,j}(a_{j+k,j})(P_s-P_{s+1})\|\cdot(1+\varepsilon) \cdot \|\xi\|
   \\
  &=\|\overline{\pi\times_J t}(a)(P_s-P_{s+1})\|\cdot(1+\varepsilon) \cdot \|\xi\|.
  \end{align*}    
Therefore by arbitrariness of $\varepsilon$ and $\xi$ we have 
 $$\| \overline{\pi\times_J t}(E(a))(P_s-P_{s+1})\|\leq  \|\overline{\pi\times_J t}(a)(P_s-P_{s+1})\|\quad \textrm{ for}\quad s=0,...,r-1.
 $$
  Similarly one shows that $\| \overline{\pi\times_J t}(E(a))P_r\|\leq  \|\overline{\pi\times_J t}(a)P_r\|$.
  Hence by \eqref{maximum norms to help} inequality \eqref{inequality * to be proved} holds and the proof is complete.
     
\end{proof}
The requirement  $R_J=R_{J_0}$ is automatically satisfied whenever one deals with faithful representations and in particular every faithful representation of $\OO(J_0,X)$ extends to a faithful representation of $\DR(J,X)$ for an appropriate ideal $J$.
\begin{theorem} Let   $(\pi,t)$ be a representation of $X$ on Hilbert space $H$. Put
$$
J=\{a\in A: \overline{\pi}_{1,1}(\phi(a))=\pi(a)\}\qquad  \textrm{ and }\qquad J_0=J\cap J(X).
$$
If  $\pi\times_{J_0} t$   is faithful on  $\OO(J_0,X)$ (resp. on spectral subspaces of $\OO(J_0,X)$), then $R_J=R_{J_0}$ and representation $\overline{\pi\times_J t}$ is faithful on $\DR(J,X)$ (resp. on spectral subspaces of $\DR(J,X)$).
\end{theorem}
\begin{proof}
Representations $\pi\times_{J_0} t $ and $\overline{\pi\times_J t}$ are intertwined  by the natural homomorphism $\Psi$, see Proposition \ref{Necessary conditions proposition}: 
$
(\pi\times_{J_0} t)(a)=\overline{\pi\times_J t}\,\,(\Psi(a))$, $a\in \OO(J_0,X).
$
Hence, if  $\pi\times_{J_0} t$   is faithful on spectral subspaces of $\OO(J_0,X)$, then so is $\Psi$. As faithfulness of $\Psi$ on spectral subspaces of  $\OO(J_0,X)$ implies the equality $R_J=R_{J_0}$ the assertion follows from  Theorem \ref{main theorem of this section}.
\end{proof}

\begin{corollary}\label{if this is  true theorem}
If  $\pi\times_{J_0} t$   is a faithful representation of $\OO(J_0,X)$ on a Hilbert space $H$, then
$$
J_0=J\cap J(X) \quad  \textrm{ and }\quad R_J=R_{J_0} \quad  \textrm{ where }\quad J=\{a\in A: \overline{\pi}_{1,1}(\phi(a))=\pi(a)\},
$$
and  $\pi\times_{J_0} t$  extends to a   faithful  representation $\overline{\pi\times_J t}$ of $\DR(J,X)$.
\end{corollary}

 \bibliographystyle{amsplain}

\end{document}